\documentclass[a4paper,oneside,reqno]{amsart}
\usepackage[utf8]{inputenc}
\usepackage[T1]{fontenc}
\usepackage[british]{babel}

\usepackage{regexpatch}
\usepackage{amsmath}

\usepackage[foot]{amsaddr}

\usepackage{stix}
\usepackage{accents}
\usepackage{enumitem}

\usepackage{thmtools}

\usepackage{graphicx}

\usepackage[dvipsnames]{xcolor}
\usepackage{tikz}
\usepackage{tikz-3dplot}
\usepackage{pgfplots}
\pgfplotsset{compat=newest}
\usepackage{float}

\usepackage[labelformat=simple]{subcaption}

\usepackage{bookmark}
\usepackage[noabbrev,capitalise]{cleveref}

\numberwithin{equation}{section}

\newcommand\numberthis{\addtocounter{equation}{1}\tag{\theequation}}

\hypersetup{colorlinks,linkcolor={Cerulean},citecolor={ForestGreen}}

\renewcommand{\subset}{\subseteq}

\newcommand{\R}{\mathbb R}
\newcommand{\Z}{\mathbb Z}
\newcommand{\N}{\mathbb N}
\newcommand{\E}{\mathcal E}

\newcommand{\eins}{\mathbb 1}
\newcommand{\univ}{\mathbb e}

\newcommand{\cb}{C_{\textnormal{b}}}
\newcommand{\cc}{C_\textnormal{c}}

\newcommand{\cic}{C^1_{\textnormal{c}}}
\newcommand{\mplus}{\mathcal M_+}
\newcommand{\mb}{\mathcal M_{\textnormal{b}}}

\newcommand{\fcb}{\mathscr F\cb^\infty}

\newcommand{\im}{\textnormal{Im}}

\newcommand{\supp}{\textnormal{supp}}

\newcommand{\dom}{\mathcal D}

\newcommand{\de}{\mathop{}\!\mathrm{d}}

\newtheorem{thm}{Theorem}[section]

\newtheorem{prop}[thm]{Proposition}

\theoremstyle{definition}

\newtheorem{cond}[thm]{Condition}

\theoremstyle{remark}
\newtheorem{lem}[thm]{Lemma}
\newtheorem{rem}[thm]{Remark}

\setenumerate[1]{label={(\roman*)}}

\crefname{cond}{Condition}{Conditions}
\crefname{lem}{Lemma}{Lemmas}

\crefname{equation}{}{}
\Crefname{equation}{}{}

\crefname{enumi}{}{}
\Crefname{enumi}{}{}

\newlist{thmlist}{enumerate}{1}
\setlist[thmlist]{label=(\roman{thmlisti}),
	ref=\thethm\,(\roman{thmlisti}),
	noitemsep}

\Crefname{listthm}{Theorem}{Theorems}
\Crefname{listlem}{Lemma}{Lemma}
\Crefname{listrem}{Remark}{Remarks}
\Crefname{listcond}{Condition}{Conditions}
\Crefname{listcor}{Corollary}{Corollaries}

\addtotheorempostheadhook[thm]{\crefalias{thmlisti}{listthm}}
\addtotheorempostheadhook[lem]{\crefalias{thmlisti}{listlem}}
\addtotheorempostheadhook[rem]{\crefalias{thmlisti}{listrem}}
\addtotheorempostheadhook[cond]{\crefalias{thmlisti}{listcond}}
\addtotheorempostheadhook[cor]{\crefalias{thmlisti}{listcor}}

\allowdisplaybreaks

\begin{document}

\title
[Mosco convergence of gradient forms]{Mosco convergence of gradient forms with\\non-convex interaction potential}
\date{24.06.2024}

\author{Martin Grothaus\textsuperscript{ \lowercase{a},$\dagger$}}
\address{\textsuperscript{\lowercase{a}}University of Kaiserslautern-Landau, Department of Mathematics, Kaiserslautern, Germany.}
\email{\textsuperscript{$\dagger$}grothaus@mathematik.uni-kl.de}

\author{Simon Wittmann\textsuperscript{ \lowercase{b},$\ddag$,$\ast$}}
\address{\textsuperscript{\lowercase{b}}The Hong Kong Polytechnic University, Department of Applied Mathematics, Hong Kong, China.}
\email{\textsuperscript{$\ddag$}simon.wittmann@polyu.edu.hk}
\thanks{\textsuperscript{$\ast$}Corresponding author}

\makeatletter
\@namedef{subjclassname@2020}{%
	\textup{2020} Mathematics Subject Classification}
\makeatother



\begin{abstract}
	This article provides a new approach to address Mosco convergence of gra\-di\-ent-type Dirichlet forms, $\mathcal E^N$ on $L^2(E,\mu_N)$ for $N\in\mathbb N$, 
	in the framework of converging Hilbert spaces by K.~Kuwae and T.~Shioya \cite{kuwae}. 
	The basic assumption is weak measure convergence of the family
	${(\mu_N)}_{N}$ on the state space $E$ - either a separable Hilbert space or a locally convex topological vector space.
	Apart from that, the conditions on ${(\mu_N)}_{N}$ try to impose as little restrictions as possible.
	The problem has fully been solved if the family ${(\mu_N)}_{N}$
	contain only log-concave measures, due to L.~Ambrosio, G.~Savaré and L.~Zambotti \cite{zambros}.
	However for a large class of convergence problems the assumption of log-concavity fails. The article suggests a way to overcome this hindrance,
	as it presents a new approach. Combining the theory of Dirichlet forms 
	with methods from numerical analysis we find abstract criteria for Mosco convergence of standard gradient forms with varying reference measures.
	These include cases in which the measures are not log-concave.
	To demonstrate the accessibility of our abstract theory we discuss a first application,
	generalizing an approximation result by S.~K.~Bounebache and L.~Zambotti
	from \cite{boune}.
\end{abstract}

\maketitle

\section{Introduction}

\subsection{Mosco convergence of gradient forms}

The abstract framework presented by {Ka\-zu\-hi\-ro} Kuwae and Takashi Shioya in \cite{kuwae}
takes up the functional analytic ideas of Umberto Mosco, who in \cite{umberto} investigates the convergence of spectral structures on a Hilbert space,
and fits it into a setting of varying Hilbert spaces. Their method has found application in partial differential equations, see e.g.~\cite{lancia},
and in probability theory, see e.g.~\cite{boune,nvili,andres}. The probabilistic disputes are often motivated by problems from statistical mechanics
involving the scaling limit of a dynamical system. There, one typically starts by looking at a
statistically distributed ensemble of interacting particles or sites in a finite volume. 
Sites are the interacting entities, which replace the physical particles, in phenomenological or effective models.
Technically, a finite volume marks a subset $E_N$ in the collection of all states $E$, which is characterized by a limited number of degrees of freedom.
That number increases as the index $N\in\N$ increases. Descriptively, the limit of $N\to\infty$ represents a transition from a micro- or mesoscopic 
understanding of the problem to a  macroscopic point of view. 
On $E_N$ a natural reference measure is provided by the Lebesgue measure. At each point
there is a natural tangent space which isomorphic to the Euclidean space.
A probability $\mu_N$ with a density proportional to $\exp(-V_N)$ describes a system in its thermal equilibrium.
The function $V_N:E_N\to\R$ is called potential, or Hamiltonian, assigned to a microscopic state. Once the weak measure convergence of $\mu_N$ for $N\to\infty$ is known,
the closest question related to a dynamical result is concerned with the fluctuations around the equilibrium. For each $N$ such a dynamic should
admit $\mu_N$ as a reversible measure and heuristically behave according to the
stochastic differential equation 
\begin{equation*}
	\de X_t=-\nabla V_N\de t + \sqrt{2}\,\de W_t.
\end{equation*}
Convergence of the finite-dimensional distributions of the laws on $E$ under the scaling limit $N\to\infty$ is equivalent to the Mosco convergence of the gradient-type Dirichlet forms  
\begin{equation}\label{eqn:helloworld}
	\E^N(u,v)\,=\,\int\limits_{E_N}{\langle \nabla u, \nabla v\rangle}_{E_N} \de\mu_N,\qquad u,v\in\dom(\E^N).
\end{equation}
The elements of $\dom(\E^N)$ are contained in a local Sobolev space $H^{1,1}_\textnormal{loc}$ over $E_N$.
The problem becomes more involved the less regularity is assumed for $V_N$. 
The applications we consider, do not require the continuity of $V_N$, for example.

Given the weak convergence of the invariant measures, the problem of identifying the asymptotic Dirichlet form becomes an interesting 
topic on its own right, as it stands at the beginning of a further discussion on the probabilistic side.
Gradient forms appear as standard examples in the books of \cite{maro,fuku94}. If the state space $E$ is Polish and $m$ is a probability measure 
on its Borel $\sigma$-algebra $\mathcal B(E)$, then
the family of local, quasi-regular, conservative and symmetric Dirichlet forms on $L^2(E,m)$ are in 1:1 correspondence 
with the family of conservative $m$ -symmetric diffusion processes
on $(E,\mathcal B(E))$ (up to equivalence). 
A conservative diffusion process $X=(\Omega,\mathscr F,{(X_t)}_{t\geq 0},{(P_x)_{x\in E}})$ is a Hunt process with path space $C([0,\infty),E)$. 
The transition function 
$p_t(x,A):=P_x(\{X_t\in A\}),$ $x\in E,$ $A\in\mathcal B(E),$ $t\geq 0,$
is $m$ -symmetric and hence $m$ is an invariant measure.
Extending the linear operator 
\begin{equation*}
	\tilde p_t :u\,\mapsto \,\int\limits_Eu(y)\de p_t(\cdot  ,\de y),
\end{equation*}
which acts on the bounded, measurable functions on $E$, to a symmetric contraction operator $T_t$ on $L^2(E,m)$ for $t\geq 0$, the relation of $X$ and $\E$ is given by the equations
\begin{gather}
	\dom(\E)\,=\,\Big\{\,u\in L^2(E,m)\,\Big|\, \sup_{t>0}\,\frac{1}{t}\int\limits_{E}u\,(u-T_tu)\de m\,<\infty\,\Big\}\nonumber\\
	\textnormal{and}\quad\E(u,v)\,=\,\lim_{t\to 0}\,\frac{1}{t}\int\limits_{E}u\,(v-T_tv)\de m.\label{eqn:associate}
\end{gather}
The family ${(T_t)}_{t\geq 0}$ forms a strongly continuous contraction semigroup on $L^2(E,m)$.
Given a family of diffusion processes 
$\{$ $X^N=(\Omega_N$, $\mathscr F^N$, ${(X^N_t)}_{t\geq 0}$, ${(P^N_x)}_{x\in E})$,
$X=(\Omega$, $\mathscr F$, ${(X_t)}_{t\geq 0}$, ${(P_x)}_{x\in E})$ $\}$,  
where $X^N$ is $m_N$ -symmetric and $X$ is $m$ -symmetric, we now write $\widetilde P_N(B):=\int_E P^N_x(B)\de m_N(x)$ for $B\in\mathscr F^N$, $N\in\N$, 
and $\widetilde P(B):=\int_E P_x(B)\de m(x)$ for $B\in\mathscr F$.
Convergence of the finite-dimensional distributions of equilibrium fluctuations, which reads
\begin{multline*}
	\lim_{N\to\infty}\,\int\limits_{\Omega_N} f_1(X^N_{t_1})\cdot f_2(X^N_{t_1+t_2})\cdot \dots \cdot f_k(X^N_{t_1+t_2+\dots+t_k})\de \widetilde P_N\\
	\begin{aligned}
		&=\,\lim_{N\to\infty}\,\int\limits_ET_{t_1}^N(f_1\cdot T_{t_{2}}^N(\dots T_{t_{k-1}}^N (f_{k-1}\cdot T^N_{t_k}f_k)\dots))\de m_N(x)\\
		&=\,\int\limits_ET_{t_1}(f_1\cdot T_{t_{2}}(\dots T_{t_{k-1}} (f_{k-1}\cdot T_{t_k}f_k)\dots))\de m(x)
	\end{aligned}\\
	=\,\int\limits_{\Omega} f_1(X_{t_1})\cdot f_2(X_{t_1+t_2})\cdot \dots \cdot f_k(X_{t_1+t_2+\dots+t_k})\de \widetilde P
\end{multline*}
with $f_1,\dots,f_k\in\cb (E)$, $t_1,\dots,t_k\in[0,\infty)$, $k\in\N$, 
is equivalent to Mosco convergence of the corresponding sequence of Dirichlet forms towards the corresponding asymptotic form.
This is due to the theorem of Mosco-Kuwae-Shioya, as stated in \cite[Theorem 2.4]{kuwae}.
Mosco convergence is formulated in terms of two conditions, (a) of \cite[Definition 2.1]{umberto} respectively 
(F1') of \cite[Definition 2.11]{kuwae}), and (b) of \cite[Definition 2.1]{umberto} respectively 
(F2) of \cite[Definition 2.11]{kuwae}).
In this text we call them {(M1)} and {(M2)}.

The exact domain of the asymptotic form plays a crucial role.
Identifying a Mosco limit includes making a statement concerning the scope of its domain.
This is reflected in the contrasting interplay between the two conditions when they are looked at independently.
If the sequence ${(\E^N)}_N$ satisfies {(M1)} w.r.t.~the asymptotic form $\E^*$ and
simultaneously satisfies {(M2)} w.r.t.~another asymptotic form $\E^{**}$, then $\dom(\E^{**})\subset\dom (\E^*)$ and the quadratic 
form of $\E^{*}$ is dominated by that of $\E^{**}$, i.e.~$\E^{*}(u,u)\leq\E^{**}(u,u)$ for $u\in\dom(\E^{**})$.
To show Mosco convergence we thus have to see why the `smallest' asymptotic form for which \textit{(M2)} holds and the `biggest' asymptotic form
for which {(M1)} holds coincide. \newcommand{\amind}{\cref{eqn:am1nd}}

Another important aspect, when it comes to verifying the two conditions {(M1)} and {(M2)}, is that, given {(M2)}, one can always retreat to prove a slightly modified version of
{(M1)} instead. The modified version is easier to check in such cases, where one has additional knowledge on the resolvents associated with ${(\E^N)}_{N\in\N}$.
The modified version of {(M1)} is formulated in Condition (a) of \amind{}. We can benefit from it due to the fact that the resolvent associated to a Dirichlet form is sub-Markovian.
It is an essential ingredient in the proof of \cref{thm:m1}, a key result in this article.
To our best knowledge, however, this useful modification of the conditions for Mosco convergence has not been stated explicitly in the literature.

Despite the increasing number of applications and its usage in various fields,
schematic guides to deal with Mosco convergence and related results in a general setting are rather rare to find.
With \cite{koles,kolesni,belhadj,toll,puga,suzuki} we would like to name some sophisticated works, which fall into this category. 
This article endeavours to find new methods and tools in the topic of Mosco convergence. 
The idea for our approach is based on an observation in a finite-dimensional vector space $V$.
The properties {(M1)} and {(M2)} are equivalent to each other, if the term of Mosco convergence refers 
to a sequence of symmetric, non-negative definite bilinear forms on $V$.
To benefit from this, we are inspired by a method which is used in numerics and better known under the name of Finite Elements.
This transfer presents the most significant innovation of this article, as we think.
Our motivation to derive and present the abstract theory in this survey is to provide a suitable groundwork 
in the field of Dirichlet forms to address problems form statistical mechanics.
This intention manifests itself in the type of potential functions which are considered and in the way the conditions are formulated.
The characteristic feature of our approach is, that it tries to use as little information as possible on the 
asymptotic invariant measure $\mu$. Instead we formulate the assumptions in terms of the Radon-Nikodym derivatives of the approximating measures ${(\mu_N)}_{N\in\N}$
- more precisely, on the densities of suitable disintegrations.   
We apply our results to a problem, whose relevance originates from the discussion of \cite{boune} dealing with scaling limits of skew reflected stochastic interface models.

\subsection{Interface models with skew reflection}\label{sec:intm}

Even in the comfortable case, in which $E$ is a Hilbert space and the limit $\mu$ of ${(\mu_N)}_N$ admits a density w.r.t.~a Gaussian measure, 
the task of proving Mosco convergence for gradient-type forms can be challenging, depending on the nature of the density.
In \cite{boune}, Said Karim Bounebache and Lorenzo Zambotti
investigate the instance, where $E=L^2((0,1),\de s)$ and $E_N$ is the linear span of indicator functions $\eins_{[2^{-N}(i-1),2^{-N}i)}$, $i=1,\dots, 2^N$.
Their study motivates the application we discuss in \cref{sec:gauss} and we shortly summarize the aspects of their work which are relevant to this article in the next paragraph.
The law of a standard Brownian bridge, i.e.~a Brownian bridge starting from $0$ at time $0$ and ending in $0$ at time $1$, defines a Gaussian measure $\tilde\mu$ on $E$.
Bounebache and Zambotti show Mosco convergence for the sequence of gradient forms ${(\E^{N})}_N$ as in \cref{eqn:helloworld}. The respective reference measure of the gradient form is chosen as
the probability with
\begin{equation}\label{eqn:nonlog}
	\de \mu_N(h)\,\propto\,\exp(-V(h))\,\de\tilde\mu_N(h),\qquad V:E\ni h\,\mapsto\, \int\limits_0^1 f(h(s))\de s,
\end{equation}
where $f:\R\to\R$ is of bounded variation 
and $\tilde \mu_N$ denotes the image measure under the orthogonal projection $E\to E_N$ of $\tilde\mu$.
In this case, the domain of $\E^N$ coincides with the Sobolev space $H^{1,2}(E_N,\tilde\mu_N)$ of the Gaussian measure $\tilde\mu_N$.
The asymptotic form is a perturbed version of the standard gradient form on $E$ in the Gaussian case:
\begin{equation}\label{eqn:insteadof}
	\E(u,v)\,=\,\int\limits_{E}{\langle\nabla u,\nabla v\rangle}_E \,\exp(-V)/Z\,\de\tilde\mu,\
	\quad u,v\in\dom(\E).
\end{equation}
The domain of $\E$ coincides with the Sobolev space $H^{1,2}(E,\tilde\mu)$ and $Z:=\int\exp(-V)\de\tilde\mu$.
With their result for convergence of ${(\E^N)}_{N\in\N}$ towards $\E$ in the sense of Mosco, \cite[Thm. 5.6]{boune}, the authors provide an approximation statement for a skew-reflected stochastic heat equation with Dirichlet boundary conditions on $[0,1]$.
The corresponding stochastic partial differential equation reads
\begin{align}\label{eqn:heat}
	&\frac{\partial u}{\partial t}\,=\,\frac{1}{2}\frac{\partial^2 u}{\partial \theta^2}+\frac{\alpha}{2}\frac{\partial}{\partial\theta}l_{t,\theta}^0+\dot W,\nonumber\\
	&u(t,0)\,=\,u(t,1)\,=\,0,\nonumber\\
	&u(0,\theta)\,=\,u_0(\theta),\quad \theta\in[0,1],
\end{align}
where $\{{l_{t,\theta}^0}|\theta\in[0,1]\}$ is the family of local times at $0$ accumulated over $[0,1]$ by the process $\{{u(t,\theta)}|\theta\in[0,1]\}$,
$\dot W$ is space-time white noise and $u_0\in E$.
The Dirichlet form $\E$ is quasi-regular. 
For the choice $f=\alpha\eins_{(-\infty,0]}$ the diffusion process $X=(\Omega,\mathscr F,{(X_t)}_{t\geq 0},{(P_x)_{x\in E}})$ with state space $E$ which is associated with 
$\frac{1}{2}\E$
yields a weak solution to \cref{eqn:heat} for all starting values in $E$ except for a set of capacity zero.
The probability 
\begin{equation*}
	\mathbb P_\mu(B):=\int\limits_E P_x\big({(X_t)_{t\geq 0}}\in B\big)\de\mu(x)
\end{equation*}
for $B$ from the Borel $\sigma$-algebra of $C([0,\infty),E)$ defines the stationary law of $X$. 
The probabilistic equivalence of Mosco convergence, ${(\E^N)}_{N\in\N}$ towards $\E$, is the weak convergence of marginals of a sequence of laws ${(\mathbb P_N)}_{N\in\N}$, 
towards the marginals of $\mathbb P_\mu$.
The law ${\mathbb P_N}$ arises from the stochastic differential equation describing the equilibrium fluctuations of a microscopic interface model.
For an index $N\in\N$ the dynamic of the microscopic interface $u_N(t,\theta)$, $t\geq 0$, $\theta\in[0,1]$, is given by a diffusion process with state space $E_N$. 
Let ${(u_t^{N,i})}_{t\geq 0}$ denote its coordinate processes, i.e.~we write
\begin{align*}
	u_N(t,\theta)&\,=\,2^{N/2}\sum_{i=1}^{2^N} \eins_{[2^{-N}(i-1),2^{-N}i)}(\theta)\,u_t^{N,i},\\
	\text{where}\quad u^{N,i}_t&\,=\,2^{N/2}\,\langle \eins_{[2^{-N}(i-1),2^{-N}i)},u_N(t,\,\cdot\,)\rangle_E\quad \text{for }t\geq0,\,i=1,\dots, 2^N.
\end{align*}
The coordinate processes follow the system of stochastic differential equations
\begin{equation}\label{eqn:microheat}
	\de u^{N,i}_t\,=\,-A_N\big(u^{N,1}_t,\dots,u^{N,2^N}_t\big)^\textnormal{T}\de t+\beta_N\de l_t^{N,i,0} +\de W_t^{N,i}, \\
\end{equation}
for $i=1,\dots, 2^N$, where 
\begin{equation*}
	\beta_N:=\frac{1-e^{-\alpha 2^{-N}}}{1+e^{-\alpha 2^{-N}}},
\end{equation*}
${(l_t^{N,i,0})}_{t\geq 0}$ is the central local time of ${(u^{N,i}_t)}_{t\geq 0}$ at $0$, 
\sloppy ${(W_t^{N,1},\dots, W_t^{N,2^N})}_{t\geq 0}$ is a $2^N$-di\-men\-sion\-al Brownian motion and
$A_N\in\R^{2^N\times 2^N}$ is the inverse matrix of 
\begin{equation}\label{eqn:covmat}
	\Big[\,\sum_{k=1}^\infty\frac{2^{N+2}}{\pi^2 k^2}\int\limits_{2^{-N}(i-1)}^{2^{-N}i}\sin(\pi k s)\de s\int\limits_{2^{-N}(j-1)}^{2^{-N}j}\sin(\pi k s)\de s\,\Big]_{i,j=1}^{2^N}.
\end{equation}
The matrix of \cref{eqn:covmat} equals
\begin{equation*}
	2\Big[\,\sum_{k=1}^\infty \lambda_k\, \big\langle \varphi_k, 2^{N/2} \eins_{[2^{-N}(i-1),2^{-N}i)}\big\rangle_E\,\big\langle \varphi_k,2^{N/2} \eins_{[2^{-N}(j-1),2^{-N}j)}\big\rangle_E\,\Big]_{i,j=1}^{2^N}
\end{equation*}
with  ${(\varphi_k)}_{k\in\N}$ being an orthonormal basis of eigenvectors for the corresponding sequence of eigenvalues ${(\lambda_k)}_{k\in\N}$ 
of the covariance operator of $\tilde\mu$, i.e.~the Laplace operator $L$ with Dirichlet boundary conditions on $(0,1)$,
\begin{equation*}
	Lh=h'',\quad \dom(L)=W^{2,2}((0,1))\cap W^{1,2}_0((0,1)).
\end{equation*}
The invariant measure of \cref{eqn:microheat} is the probability with 
\begin{equation}\label{eqn:nonlog2}
	\de m_N(x)\,\propto\,\exp\Big(-x^\textnormal{T}A_Nx-\alpha2^{-N}\sum_{i=1}^{2^N}\eins_{(-\infty,0]}(x_i)\Big)
\end{equation}
and is the image measure of $\mu_N$ under the bijection
\begin{equation*}
	E_N\ni h\,\mapsto\, 2^{N/2} \sum_{i=1}^{2^N}\langle h,\eins_{[2^{-N}(i-1),2^{-N}i)}\rangle_E\,\univ_i\,\in\R^{2^N}.
\end{equation*}
Likewise, the image form of $\E^N$ under this map is the gradient Dirichlet form $\E ^{N,\text{micro}}$ on $L^2(\R^{2^N},m_N)$ with reference measure $m_N$, i.e.~the closure of
\begin{equation*}
	\E ^{N,\text{micro}}=\sum_{i=1}^{2^N}\,\int\limits_{\R^{2^N}}\frac{\partial u}{\partial x_i}\,\frac{\partial v}{\partial x_i}\de m_N
\end{equation*}
with pre-domain $C_b^1(\R^{2^N})$ on $L^2(\R^{2^N},m_N)$.
The diffusion process $X^N=(\Omega_N,$ $\mathscr F^N,$ ${(X^N_t)}_{t\geq 0},$ ${(P^N_x)_{x\in \R^{2^N}}})$ with state space $\R^{2^N}$
which is associated to the form $\frac{1}{2}\E ^{N,\text{micro}}$ yields weak solutions to \cref{eqn:microheat} for every starting point in $\R^{2^N}$ except for a set of zero capacity.
With the bijection 
\begin{equation}\label{eqn:micromacro}
	C\big([0,\infty),\R^{2^N}\big)\ni  {\omega(\,\cdot\,)}\,\mapsto\,2^{N/2}\sum_{i=1}^{2^N} \eins_{[2^{-N}(i-1),2^{-N}i)}\,\omega_i(\,\cdot\,)\,\in C([0,\infty),E_N)
\end{equation}
the diffusion process of the coordinates of the microscopic interface model is interpreted as an approximation for the macroscopic dynamic.
The image measure of the equilibrium law 
\begin{equation*}
	B\mapsto \int\limits_{\R^{2^N}} P^N_x\big({(X_t)_{t\geq 0}}\in B\big)\de m_N(x),
\end{equation*}
where $B$ is from the Borel $\sigma$-algebra of $C([0,\infty),\R^{2^N})$, under the map of \cref{eqn:micromacro}
is denoted by $\mathbb P_N$.
Via Mosco convergence of Dirichlet forms the weak measure convergence of the finite-di\-men\-sion\-al distributions of ${(\mathbb P_N)}_{N\in\N}$
towards those of $\mathbb P_\mu$ is shown.
The difficulty in the proof for Mosco convergence 
lies in the fact that the measure of \cref{eqn:nonlog2} is not log-concave, due to the non-convexity of the perturbing potential, as the authors point out.
A standard approach to prove Mosco convergence, which uses gradient bounds for the semigroups, is thus not available in this case.
Gradient bounds have been used to prove convergence of diffusion processes in \cite{Zambotti} (in combination with integration by parts formulae) or in \cite{zambros}.
Such can be obtained either probabilistically (see e.g.~\cite{priola, biel}) or analytically (see e.g.~\cite{baud, luigi}).
However, the typical requirement for the convexity of the potential $V$, i.e.~positivity of the Hessian $\nabla\nabla V$ or equivalently the log-concavity of $\exp(-V)$, 
unites the different approaches to gradient bounds and prohibits their applicability in this context. 
To prove the convergence result of \cite[Thm. 5.6]{boune} despite the lack of such tools, the authors exploit an integration by parts formula for the limiting invariant measure $\mu$.
Estimates for the local time of the Brownian bridge at $0$, which appears in the integration by parts formula for $\mu$, from the abundant literature on that particular measure are used.
The methods of \cite[Thm. 5.6]{boune} do not generalize to skew-reflected interface models with Gaussian reference measures which are different from the Brownian bridge in a straightforward way.
From the physical point of view, however, it is strongly desirable to be able to take also different Gaussian measures into consideration. The most prominent reason is that
a physical interface is observed in a $(2+1)$-di\-men\-sion\-al space. A more fitting choice for the state space $E$ 
is therefore the $L^2$-space over $[0,1]^2$, for example. Moreover, the covariance of the Gaussian measure mimics the intrinsic physical nature of the interface, i.e.~its stiffness or surface tension. 
The covariance operator determines the linear part of the drift in the stochastic dynamics. In the case of \cref{eqn:heat} it is the Laplacian on the interval $(0,1)$.
To have more options for that would increase the physical scope of a probabilistic model. 
In the next paragraph we describe a probabilistic interpretation of this article's results on Mosco convergence, \cref{thm:m1},  \cref{cor:mosco}, and \cref{thm:bougen} in particular, 
which provide broader options 
concerning the state space $E$ and the Gaussian reference measure.

Generalizing the statement of \cite[Thm. 5.6]{boune} for
a broader class of reference measures and state spaces,
we consider a generic finite measure $\lambda$ on a $\sigma$-algebra over a set $\Omega$ and define $E:=L^2(\Omega,\lambda)$.
Our assumptions allow to treat non-convex potentials, as those of the type considered in \cite{boune}.
The reference measure on $E$ we consider reads
\begin{equation*}
	\de\mu(h)\,:=\,\frac{1}{Z}\exp(-V(h))\de\tilde\mu(h)\quad\textnormal{with}\quad V:E\ni h\,\mapsto\,\int\limits_\Omega f(h(\omega))\de \lambda(\omega),
\end{equation*}
where $\tilde\mu$ is mean-zero, non-degenerate Gaussian, $f:\R\to\R$ is a function of bounded variation and $Z:=\int\exp(-V)\de\tilde\mu$.
For $f$ it is necessary to assume that for every point of discontinuity the corresponding level sets of $\tilde\mu$ in $\Omega$ are $\lambda$ -nullsets almost surely (see \cref{cond:lvl}).
We look at an increasing sequence of exhausting, finite-dimensional subspaces $E_N\nearrow E$ and the corresponding sequence of orthogonal projections ${(\pi_N)}_N$. 
Using the probabilities
\begin{equation*}
	\de\mu_N(h)\,\propto\,\exp(-V(h))\de (\tilde\mu\circ\pi_N^{-1})(h)
\end{equation*}
we define the gradient Dirichlet forms $\E^N$ as in \cref{eqn:helloworld}. The domain of $\E^N$ coincides with the Sobolev space $H^{1,2}(E_N,\tilde\mu\circ\pi_N^{-1})$
of the Gaussian measure $\tilde\mu\circ\pi_N^{-1}$. We show
Mosco convergence of ${(\E^{N})}_N$ towards the perturbed Gaussian gradient form $\E$ with
\begin{equation*}
	\E(u,v)\,=\,\int\limits_{E}{\langle\nabla u,\nabla v\rangle}_E \,\de\mu,\
	\quad u,v\in\dom(\E)=H^{1,2}(E,\tilde\mu),
\end{equation*}
on $L^2(E,\mu)$.
Hence, our convergence result accommodates the physically more relevant cases of $(2+1)$-dimensional interface models.
Gaussian measures on $E:=L^2((0,1)^2)$ can be constructed by means of the Fourier transform (see e.g.~\cite[Theorem 2.3.1]{boga98}).
For example, we can consider the Gaussian measure $\tilde\mu$ on $E$ whose characteristic function is given via the inverse operator $Q$ of 
the squared Laplacian with Dirichlet boundary conditions
\begin{equation*}
	Lh=\big(\frac{\partial^2}{\partial^2 x} +\frac{\partial^2}{\partial^2 y}\big)\big(\frac{\partial^2}{\partial^2 x} +\frac{\partial^2}{\partial^2 y}\big)h,\quad
	\dom(L)=W^{4,2}((0,1)^2)\cap W^{3,2}_0((0,1)^2).
\end{equation*}
Then, $Q$ is a symmetric, non-negative, nuclear operator on $E$ and $\tilde\mu$ is the unique Gaussian measure with
\begin{equation*}
	\int\limits_E\exp{\big(i\langle h_1,h_2\rangle_E\big)}\de\tilde\mu(h_1)=\exp\big(-\frac{1}{2}\langle h_2,Qh_2\rangle_E\big),\quad h_2\in E.
\end{equation*}
Arranging all elements from the orthonormal family of eigenvectors
\begin{equation*}
	(0,1)^2\ni(x,y)\,\mapsto\,2\sin(l\pi x)\sin(m\pi y)
\end{equation*}
and corresponding eigenvalues $\pi^{-4}(l^2+m^2)^{-2}$ with indices $l,m\in\N$ into a sequence indexed by one parameter $k\in\N$, we obtain
an orthonormal basis ${(\varphi_k)}_{k\in\N}$ of eigenvectors of $Q$ with the corresponding sequence of eigenvalues ${(\varphi_k)}_{k\in\N}$.
For each $N\in\N$ we choose a finite-dimensional subspace $E_N$ of $E$ such that $E_{N-1}\subseteq E_N$ for $N\geq 2$ and moreover the linear span of $\bigcup_{N\in\N}E_N$ is dense in $E$.
Let $d_N\in\N$ denote the dimension of the space $E_N$ and $\xi^{N,1},\dots,\xi^{N,d_N}$ be an orthonormal basis for $E_N$.
One possible choice, to which we refer in the following discussion, is given by
\begin{multline}\label{eqn:basisfu}
	E_N:=\textnormal{span}\Big(\Big\{u:[0,1]^2\to[0,1]\,\Big|\,\text{there exist }i,j=1,\dots, 2^N\text{ with }\\
	u(x,y)=\frac{1}{d_N}\, \eins_{[2^{-N}(i-1),2^{-N}i)}(x)\, \eins_{[2^{-N}(j-1),2^{-N}j)}(y)\,\Big\}\Big),
\end{multline}
where $d_N=2^{2N}$ and $\xi^{N,1},\dots,\xi^{N,d_N}$ is an arbitrary ordering of the basis functions in \cref{eqn:basisfu}.
If we set-up the non-log-convex perturbing potential with $f=\alpha\eins_{(-\infty,0]}$, as regarded in the previous paragraph, the probabilistic interpretation
of our result on Mosco convergence, \cref{thm:bougen}, is as follows.
Again, for each index $N$ a diffusion on $E_N$ describes the dynamical microscopic interface $u_N(t,\theta)$, $t\geq 0$, $\theta\in[0,1]^2$.
Its coordinate processes are denoted by ${(u_t^{N,i})}_{t\geq 0}$, i.e.~we write
\begin{align*}
	u_N(t,\theta)&\,=\,\sum_{i=1}^{d_N} \xi^{N,i}(\theta)\,u_t^{N,i},\\
	\text{where}\quad u^{N,i}_t&\,=\,\langle \xi^{N,i},u_N(t,\,\cdot\,)\rangle_E\quad \text{for }t\geq0,\,i=1,\dots, d_N.
\end{align*}
They are subject to the system of stochastic differential equations
\begin{equation}\label{eqn:microheat2}
	\de u^{N,i}_t\,=\,-A_N\big(u^{N,1}_t,\dots,u^{N,d_N}_t\big)^\textnormal{T}\de t+\beta_N\de l_t^{N,i,0} +\de W_t^{N,i}, \\
\end{equation}
for $i=1,\dots, d_N$, where 
\begin{equation*}
	\beta_N:=\frac{1-e^{-\alpha /d_N}}{1+e^{-\alpha /d_N}},
\end{equation*}
${(l_t^{N,i,0})}_{t\geq 0}$ is the central local time of ${(u^{N,i}_t)}_{t\geq 0}$ at $0$, 
\sloppy ${(W_t^{N,1},\dots, W_t^{N,d_N})}_{t\geq 0}$ is a $d_N$-di\-men\-sion\-al Brownian motion and
$A_N\in\R^{d_N\times d_N}$ is the inverse matrix of 
\begin{equation}\label{eqn:defA}
	\Big[\,\sum_{k=1}^\infty 2\lambda_k\,\langle \varphi_k, \xi^{N,i}\rangle_E\,\langle \varphi_k, \xi^{N,j}\rangle_E\,\Big]_{i,j=1}^{d_N}.
\end{equation}
The invariant measure of \cref{eqn:microheat2} is the probability with 
\begin{equation*}
	\de m_N(x)\,\propto\,\exp\Big(-x^\textnormal{T}A_Nx-\frac{\alpha}{d_N}\sum_{i=1}^{d_N}\eins_{(-\infty,0]}(x_i)\Big).
\end{equation*}
Under the bijection
\begin{equation*}
	E_N\ni h\,\mapsto\,  \sum_{i=1}^{d_N}\langle h,\xi^{N,i}\rangle_E\,\univ_i\,\in\R^{d_N},
\end{equation*}
$m_N$ is the image measure of $\mu_N$ and the gradient-type Dirichlet form on $L^2(\R^{d_N},m_N)$,
\begin{equation*}
	\E ^{N,\text{micro}}=\sum_{i=1}^{d_N}\,\int\limits_{\R^{d_N}}\frac{\partial u}{\partial x_i}\,\frac{\partial v}{\partial x_i}\de m_N
\end{equation*}
with pre-domain $C_b^1(\R^{d_N})$, is the image form of $\E^N$.
Hence, Mosco converge of ${(\E^N)}_{N\in\N}$ towards $\E$ provides a probabilistic interpretation connected to the large-scale asymptotic of the stochastic differential equation \cref{eqn:microheat2}.
The diffusion process $X^N=(\Omega_N,$ $\mathscr F^N,$ ${(X^N_t)}_{t\geq 0},{(P^N_x)_{x\in \R^{d_N}}})$ with state space $\R^{d_N}$
which is associated to the form $\frac{1}{2}\E ^{N,\text{micro}}$ yields weak solutions to \cref{eqn:microheat2} for every starting point in $\R^{d_N}$
except for a set of zero capacity.
The bijection 
\begin{equation}\label{eqn:micromacro2}
	C\big([0,\infty),\R^{d_N}\big)\ni  {\omega(\,\cdot\,)}\,\mapsto\,\sum_{i=1}^{d_N} \xi^{N,i}\,\omega_i(\,\cdot\,)\,\in C([0,\infty),E_N)
\end{equation}
transforms the diffusion process of the coordinates of the microscopic interface model into an approximation for the macroscopic dynamic.
Analogously to the previous paragraph, the statement of Mosco convergence concerns the image measure $\mathbb P_N$ of the equilibrium law
\begin{equation*}
	B\mapsto \int\limits_{\R^{d_N}} P^N_x\big({(X_t)_{t\geq 0}}\in B\big)\de m_N(x), 
\end{equation*}
where $B$ is from the Borel $\sigma$-algebra of $C([0,\infty),\R^{d_N})$, under the map of \cref{eqn:micromacro2}.
The weak measure convergence of the finite dimensional distributions of ${(\mathbb P_N)}_{N\in\N}$ towards those of $\mathbb P_\mu$ is the probabilistic consequence of our analytic result.
The Gaussian measure $\tilde\mu$ on $E$ can be chosen arbitrarily in our case. So, the statement on the weak convergence of the marginals of $P_N$ for $N\to\infty$
remains vaild if we define the drift $A_N$ in \cref{eqn:microheat2} as the inverse of the matrix in \cref{eqn:defA} w.r.t.~any orthonormal basis 
${(\varphi_k)}_{k\in\N}$ and values ${(\lambda_k)}_{k\in\N}$ such that 
\begin{equation*}
	Q:E\ni h\,\mapsto\, \sum_{i=1}^\infty\lambda_k\,\langle h,\varphi_k\rangle_E\,\varphi_k
\end{equation*}
is the covariance operator of a Gaussian measure on $E$. With \cref{thm:bougen}, where the Mosco convergence of ${(\E ^N)}_N$
in the context of skew-reflected interface models is stated, we intend to give a first application generalizing the statement of \cite[Thm. 5.6]{boune}. It aims
to demonstrate the utility and accessibility of the abstract convergence results of \cref{sec:lowsii,sec:stagra}.

Further applications of the abstract results are the subject of ongoing studies. In the next lines we shortly describe what is about to come as a follow-up of this article. The
methods of \cref{sec:lowsii,sec:stagra} are apt to show Mosco convergence in such cases, in which the Gaussian measures as well as the perturbations vary with the index $N$.
Instead of looking at the sequence of images under orthogonal projections of one particular Gaussian measure, we can initiate a skew-reflected interface 
model  with a sequence ${(\tilde\mu_N)}_{N}$, where $\tilde\mu_N$ is a Gaussian measure on $E_N$ and the only condition on ${(\tilde\mu_N)}_{N}$
is the weak measure convergence. 
Regarding the above example with the Brownian bridge, the criteria for Mosco convergence of Dirichlet forms developed in this article would alternatively allow to set up
an explicit microscopic model with quadratic nearest neighbour interaction, as is a typical choice in the literature on stochastic interface models.
That changes the drift $A_N$ in \cref{eqn:microheat}. Instead of the implicit definition as the inverse of the projected covariance matrix as given above, $A_N$ then is explicitly given as the tridiagonal $\R^{2^N}\times\R^{2^N}$-matrix
\begin{equation*} 
	A_N:=2^{2N}
	\begin{array}{cccccc}
		\begin{bmatrix}
			2&-1&&&&\\
			-1&2&-1&&&\\
			&-1&2&-1&&\\
			&&&&&&\\
			&&&\ddots&&\\
			\\
			&&&-1&2&-1\\
			&&&&-1&2
		\end{bmatrix}
	\end{array}.
\end{equation*}
The linear operation of $A_N$ corresponds to a discretization of the Laplacian on $(0,1)$.
A similar modification can be done in the above example with the Gaussian measure on $L^2((0,1)^2)$ whose covariance operator is given by the squared Laplacian.
$A_N\in\R^{d_N\times d_N}$ could be explicitly defined as a banded matrix arising from a natural microscopic interface model
and corresponding to a discretization of the squared Laplacian on $(0,1)^2$.
These extensions concerning the applicability of the main results in \cref{sec:lowsii,sec:stagra},
along with a characterization of the limiting macroscopic diffusion process for
a $(2+1)$-di\-men\-sion\-al interface model with skew-reflection, will be discussed in a publication to be completed shortly.
The main purpose of this article is to explain the concept and development of the suitable Dirichlet form techniques in detail.

\subsection{Outline}

We start by constructing function spaces on $\R^d$ for $d\in\N$, which are spanned by a finite selection of elementary functions. These are
obtained by the rescaling and shifting of a compactly supported, archetype function.
Hence the basis functions carry two indices - one connected to the spacial shift and the other linked to a scaling parameter or grid size.
The first part of \cref{sec:diri} sets up a particular scheme of finite elements. These accommodate the class of piecewise linear functions on $\R^d$
w.r.t.~an equidistant triangulation, called the Coxeter-Freudenthal-Kuhn triangulation.
For an element of the resulting function space, the calculation of the weak gradient and its squared norm becomes a particularly easy expression
in terms of the basis. This is stated in \cref{thm:partun}.
The second part of \cref{sec:diri} introduces quantities, which we call the residuum and the perturbation of a given probability density $\varrho$ on $\R^d$.
Their interpretation as linear functionals on $L^2(\R^d,\varrho\de x)$, respectively on $L^1(\R^d,\varrho\de x)$, is the foundation for \cref{lem:polygon}.
In terms of the operator norm we precisely express in this lemma how well a general function of finite energy can be approximated by the finite elements,
where we estimate the $L^2$- and the energy norm. \cref{lem:polygon} builds the bridge between 
the analysis of \cref{sec:tria} and the convergence theory of \cref{sec:lowse,sec:appli}.
In \cref{sec:mokushi} we first recall the essential terminology of \cite{kuwae}.
The introduction to the theory of Mosco-Kuwae-Shioya is written in a self-contained way.
Beyond the reader's comfort there are two other reasons which motivate this procedure.
Firstly, the elaborate notation in the original paper of Kuwae and Shioya, which is more focussed on the topological aspects of the theory,
surpasses the needs of this text. So, we would like to have a more basic notation.
However, there is no generally agreed custom how to initiate the concepts with a suitably simplified yet precise notation.
Secondly, the validity of the version of the theorem of Mosco-Kuwae-Shioya which is presented in our text may be known to experts of the theory, yet it
is not directly evident from the original formulation of the theorem.
To avoid any obscurity we give the proof (analogous as in \cite[Proof of Theorem 2.4]{kuwae} and \cite[Proof of Theorem 2.4]{umberto} to a large extent) in detail, 
and the version written in this article, \cref{thm:mokush}, becomes apparent.
Our main results are then stated and proven in \cref{sec:lowsii}.
The analysis of finite elements is done on $\R^d$. We initiate the abstract theory on a state space $S\times\R^d$
with a Polish space $S$. The family of reference measures ${(\mu_N)}_{N\in\N}$ on $S\times\R^d$ and their weak limit $\mu$
disintegrate into the respective conditional distributions ${(m_s^N)}_N$ and $m_s$ on $\R^d$, 
given that the canonical projection $\pi_1:S\times\R^d\to S$ takes the value $s$. Accordingly, we write
$\mu_N(A)=\int_S\int_{\R^d}\eins_A(s,x)\de m_s^N(x)\de\nu_N(s)$ for $A\in\mathcal B(S\times\R^d)$, where $\nu_N$ denotes the image measure $\mu_N\circ\pi_1^{-1}$
for $N\in\N$. \cref{thm:m1} manifests an asymptotic result for the superposition of $d$-di\-men\-sional gradient Dirichlet forms, defined
on $L^2(\R^d,m_s^N)$ respectively for $s\in S$ and $N\in\N$, with varying mixing measures $\de \nu_N(s)$.
We would like to point out that we do not assume the weak convergence of the disintegration measures ${(m_s^N)}_N$ in a pointwise sense on $S$.
This question might not even make sense since the support of the mixing measure $\nu_N$ might be a nullset w.r.t.~the asymptotic mixing measure $\mu\circ\pi_1^{-1}$.
The fact that we consider varying mixing measures requires a more delicate analysis than would be needed in the case of a fixed mixing measure.
The section closes with a discussion on the stability of the underlying assumptions of \cref{sec:lowsii}, listed in \cref{cond:mucken}, under certain perturbations.
\cref{sec:stagra} explains the relevance of \cref{thm:m1} for an effectively infinite-dimensional setting, where
the state space $E$ is a Fr{\'e}chet space and  a densely embedded Hilbert space $H$ takes the role of a tangent space to define a gradient
on the cylindrical smooth functions.
An abstract convergence result for minimal gradient forms on $E$ (as defined and analysed in \cite{rocla, albe}) 
with varying reference measures is obtained 
by applying the methods of \cref{sec:lowsii} on suitable component forms. 
The assumption concerning the domain of the asymptotic form, which \cref{cor:mosco}, the central result of \cref{sec:stagra}, requires,
is closely related to the question of Markov uniqueness and is the subject of the discussion in \cite{rocla}.
\cref{sec:gauss} then presents a Hilbert space setting in which the required characterization of the form domain is known. 
It displays the application of our abstract results in the context of the scaling limit of skew-reflected stochastic interface models (as described in \cref{sec:intm} above).

Summarizing the outline we list the central accomplishments of this article:
\begin{itemize}
	\item The version of the theorem of Mosco-Kuwae-Shioya formulated and proven in this article provides an amendment to the original version, 
	in which {(M1)} is replaced by the alternative condition (a) of \amind{}. In some instances the new condition is easier to check, since it allows to exploit knowledge on the resolvents,
	for example the sub-Markovianity.
	\item  \cref{thm:m1} ensures Mosco convergence for a sequence of superposed $d$-di\-men\-sional gradient-type Dirichlet forms on $S\times\R^d$ with a Polish space $S$. 
	The respective dis\-integration- and mixing measures vary.
	\item \cref{cor:mosco} uses the statement provided by \cref{thm:m1} to 
	derive a result on Mosco convergence for a sequence of minimal gradient forms with varying reference measures in an effectively infinite-dimensional setting.
	\item For the proof of \cref{thm:m1} we recall the method of Finite Elements, which is used in numerical analysis.
	Starting with the Coxeter-Freudenthal-Kuhn triangulation of $\R^d$ we set up a particular scheme of finite elements.
	The relevant properties, which make them useful in the theory of Dirichlet forms, are proven in \cref{thm:partun} and \cref{lem:polygon}.
	\item A first application in the context of a non-log-concave reference measure on a general state space $E=L^2(\Omega,\lambda)$ is presented in \cref{thm:bougen}.
	We consider the images of a Gaussian measure under orthogonal projections and a perturbing density $\exp(-\int_\Omega f\circ h\de\lambda)$, $h\in E$,
	for a function $f:\R\to\R$ with bounded variation.
\end{itemize}

\section{Finite Elements}\label{sec:tria}\label{sec:diri}
\newcounter{b}
\setcounter{b}{\value{section}}
\subsection{Triangulation and tent functions}\label{sec:kuhn}
We first give some notation.
A positive integer $d\geq 2$ indicating the dimension is fixed throughout this section.
In the following $\univ_k$ denotes the $k$-th unit vector of $\R^d$ for $k=1,\dots,d$. Their sum $\univ:=\univ_1+\dots+\univ_d$ is the vector whose components are constantly $1$.
For a point $x\in\R^d$ we write $[x]\in\Z^d$ for the component-wise floor of $x$, i.e.~$[x]$ is the unique element in $\Z^d$ such that $x-[x]\in[0,1)^d$.
Let $M$ be a set and $A$ be a family of maps from $M$ into $\R^d$. The family $b\,A+c$ is defined as $\{b\,a(\cdot)+c:M\to\R^d|a\in A\}$ for $b\in\R$, $c\in\R^d$.
Occasionally it is convenient to abbreviate `$x\in\im(a)$' by `$x\in a$' for a map $a:M\to\R^d$. 
Furthermore, $\eins_K:M\to\{0,1\}$ is the indicator function of an arbitrary subset $K\subset M$.
The topological support of the measure $|f(x)|\de x$ is denoted by $\supp[f]$ for a Borel measurable function on $f:\R^d\to\R$. 
We call a measurable function $\varphi:\R^d\to[0,1]$ \textbf{primal} if \vspace{1ex}
\begin{center}
	\begin{tabular}{l}
		$\bullet$ $\supp[\varphi]\subset [-2,2]^d$,\vspace{1ex}\\
		$\bullet$ $\displaystyle\int\limits_{\R^d}\varphi(x)\de x\,=\,1$,\vspace{1ex}\\
		$\bullet$ $\displaystyle\sum_{\alpha\in \Z^d}\varphi(x-\alpha)=1$\quad for $x\in\R^d$.\vspace{1ex}				
	\end{tabular}
\end{center} 
So, the last condition says that the family 
$\{\R^d\ni x\mapsto \varphi(x-\alpha)\,|\,\alpha\in\Z^d\}$ form a partition of unity.
The set of primal functions is denoted by $\mathscr C$.
For a scaling parameter $r\in(0,\infty)$ and $\alpha\in r\,\Z^d$ define $\varphi_r^\alpha(x)=\varphi((x-\alpha)/r)$, $x\in\R^d$, $\varphi\in\mathscr C$.
In this section a family $\chi_r^\alpha:\R^d\to[0,1]$, called the \textbf{tent functions},
with index $\alpha\in r\,\Z^d$ and $r\in(0,\infty)$ are constructed, which contains the element $\chi_1^0\in\mathscr C$.
This particular primal function is a piecewise linear interpolation
of the sample points $(z,\eins_{\{0\}}(z))$ over all nodes $z$ from the lattice $\Z^d$. 
The construction of $\chi_r^\alpha$ is explained step by step in the following text as their family turns out particularly useful for our purpose. Then, 
\cref{thm:partun} sums up all 
their properties which are relevant to the part following after it.
The functions' construction has a stand-alone status among the other sections and
the reader who quickly wants to get into the matter of Mosco convergence gets all the necessary preparation
for the subsequent part simply by taking note of the statements of \cref{thm:partun}.

We start by giving a triangulation of the unit $d$-cube. Its appearance traces back to \cite{cox,freu,kuhn}.
The reader can find a helpful outline of that matter in \cite{moore}.
The set $\mathscr T_1^0$ contains the shortest paths which start in $0\in\R^d$, end in $\univ\in\R^d$ and only walk along the edges of the unit $d$-cube.
An element $T\in \mathscr T_1^0$ visits exactly $d+1$ points of the set $X=\{0,1\}^{d}$ of corners. 
It reaches each of those $d+1$ points exactly once, say in an order $C_0,\dots,C_d$, where $C_0=0\in\R^d$ and $C_d=\univ$. 
We identify $T$ with an injection $\{0,\dots, d\}\to X$ writing $T(i)=C_i$.
A good way to characterize the set $\mathscr T_1^0$ exploits its one-to-one correspondence with the symmetric group $\mathscr S_d$. Let $T\in\mathscr T_1^0$. 
To find the corresponding permutation from $\mathscr S_d$ we choose
$\sigma_T(i)\in\{1,\dots,d\}$ for $i=1,\dots,d$ such that $\univ_{\sigma_T(i)}$ is the direction parallel to the edge which connects $T(i-1)$ and $T(i)$, i.e.
\begin{equation}\label{eqn:tirect}
	\univ_{\sigma_T(i)}\,=\,T(i)-T(i-1).
\end{equation}
Then the map $i\mapsto \sigma_T(i)$ is a permutation on $\{1,\dots,d\}$ indeed. The $k$-th component of the starting point ${T(0)}_k$ equals $0$
and the $k$-th component of the end point ${T(d)}_k$ equals $1$. 
So, for each $k\in\{1,\dots,d\}$ there has to be an edge of the unit cube 
parallel to $\univ_k$ along which the path of $T$ runs.
This means that $\sigma_T$ is surjective.
Moreover, since the number of edges along which $T$ runs equals $d$, the map $\sigma_T$ is also injective. By induction w.r.t.~$i$ it follows from \cref{eqn:tirect} that
\begin{equation}\label{eqn:suniv}
	T(i)\,=\,\sum_{j=1}^i \univ_{\sigma_T(j)}
\end{equation}
for $i=1,\dots,d$.
The convex hull 
\begin{align*}
	C_T\,:&=\,\Big\{\,\sum_{i=0}^d\lambda_i\,T(i)\,\Big|\,0\leq \lambda_0,\dots,\lambda_d\leq 1\textnormal{ and }\sum_{i=0}^d\lambda_i\,=\,1\,\Big\}\\
	&=\,\Big\{\,x\in\R^d\,\Big|\,0\leq x_{\sigma_T(d)}\leq x_{\sigma_T(d-1)}\leq\dots\leq x_{\sigma_T(1)}\leq 1\,\Big\}
\end{align*}
of $\{T(0),\dots,T(d)\}$ defines a polyhedron. If we choose, for given $x\in[0,1]^d$, a permutation $\sigma\in S_d$ such that
\begin{equation}\label{eqn:sigref}
	0\leq x_{\sigma(d)}\leq x_{\sigma(d-1)}\leq\dots\leq x_{\sigma(1)}\leq 1
\end{equation}
and then choose $T'$ as the unique element from $\mathscr T_1^0$ with $\sigma_{T'}=\sigma$, then it holds $x\in C_{T'}$.
Of course, the element $T'$ with $x\in C_{T'}$ is not unique as there might be more than one element in $\mathscr S_d$ for which \cref{eqn:sigref} is satisfied.
The family $\{C_T|T\in \mathscr T_1^0\}$ are called the \textbf{Coxeter-Freudenthal-Kuhn triangulation} of the unit cube.

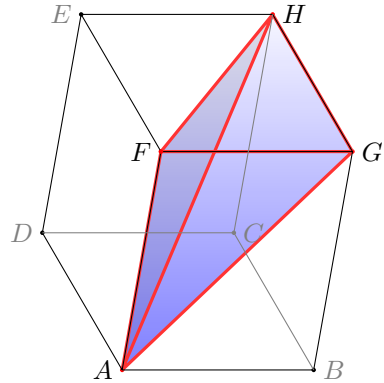
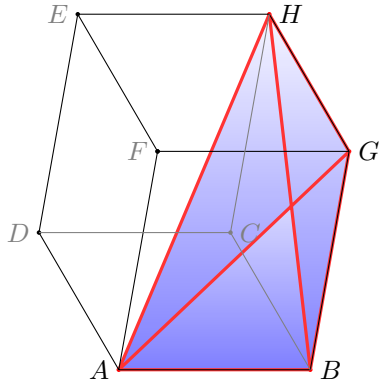
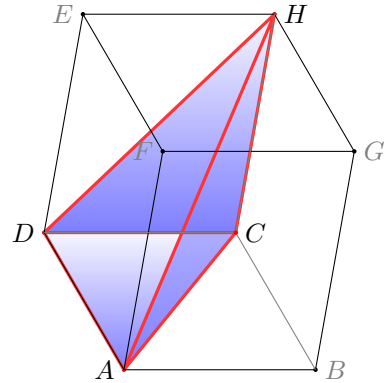
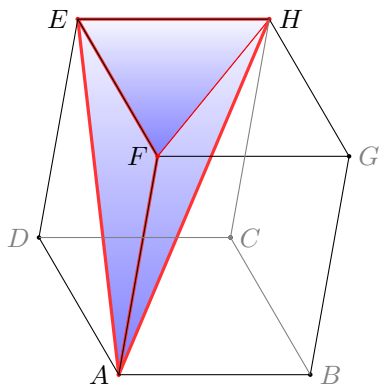
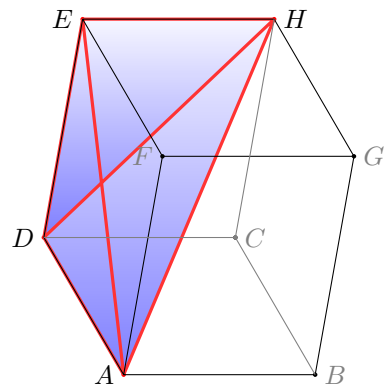
\begin{figure}
	\begin{subfigure}[t]{0.4\textwidth}
		\centering
		\begin{tikzpicture}[scale=0.42] 
	
	\def \tta{ -10.00 } 
	\def \k{    -3.00 } 
	\def \l{     6.00 } 
	\def \d{     5.00 } 
	\def \h{     7.00 } 
	
	\coordinate (A) at (0,0); 
	\coordinate (B) at ({-\h*sin(\tta)},{\h*cos(\tta)}); 
	\coordinate (C) at ({-\h*sin(\tta)-\d*sin(\k*\tta)},
	{\h*cos(\tta)+\d*cos(\k*\tta)}); 
	\coordinate (D) at ({-\d*sin(\k*\tta)},{\d*cos(\k*\tta)}); 
	
	\coordinate (Ap) at (\l,0); 
	\coordinate (Bp) at ({\l-\h*sin(\tta)},{\h*cos(\tta)}); 
	\coordinate (Cp) at ({\l-\h*sin(\tta)-\d*sin(\k*\tta)},
	{\h*cos(\tta)+\d*cos(\k*\tta)}); 
	\coordinate (Dp) at ({\l-\d*sin(\k*\tta)},{\d*cos(\k*\tta)}); 
	
	\fill[red]    (A) circle [radius=2pt]; 
	\fill[black]  (B) circle [radius=2pt]; 
	\fill[black]    (C) circle [radius=2pt]; 
	\fill[black]    (D) circle [radius=2pt]; 
	\fill[red]  (Ap) circle [radius=2pt]; 
	\fill[black] (Bp) circle [radius=2pt]; 
	\fill[red]   (Cp) circle [radius=2pt]; 
	\fill[red] (Dp) circle [radius=2pt]; 
	
	\filldraw[draw=red,bottom color=blue!50!white, top color=blue!3!white]
	(A) -- (Ap) -- (Dp);
	\filldraw[draw=red,bottom color=blue!50!white, top color=blue!3!white]
	(A) -- (Dp) -- (Cp);
	\filldraw[draw=red,bottom color=blue!50!white, top color=blue!3!white]
	(Ap) -- (Dp)  -- (Cp);

	\draw[red!80!white,-,very thick] (A) -- (Cp)
	(A) -- (Ap)
	(A) --  (Dp)
	(Ap) -- (Cp)
	(Ap)  -- (Dp)
	(Dp) --  (Cp);
	
	\draw [-,thin] (B)  --  (A)
	(Ap) -- (Bp)
	(B)  --  (C)
	(D)  --  (C)
	(A)  --  (D)
	(Ap) --  (A)
	(Cp) --  (C)
	(Bp) --  (B)
	(Bp) -- (Cp);
	
	\draw [gray,-,thin] (Dp) -- (Cp)
	(D) --  (Dp)
	(Ap) -- (Dp);

	\draw (Ap) node [right]    {$B$}
	(Bp) node [right,gray]     {$G$}
	(Cp) node [right]          {$H$}
	(C)  node [left,gray]      {$E$}
	(D)  node [left,gray]      {$D$}
	(A)  node [left]           {$A$}
	(B)  node [left,gray]      {$F$}
	(Dp) node [right]          {$C$};
	
	\fill[black]   (B) circle [radius=2pt]; 
	\fill[red]    (Dp) circle [radius=2pt]; 
	
	
\end{tikzpicture}
		\caption{The vertexes $A,B,C,H$ of $T_1$ as the corner points of the tetrahedron $C_{T_1}$.}
		\label{fig:tetra1}
	\end{subfigure}
	\hfill
	\begin{subfigure}[t]{0.4\textwidth}
		\centering
		\begin{tikzpicture}[scale=0.42] 
	
	\def \tta{ -10.00 } 
	\def \k{    -3.00 } 
	\def \l{     6.00 } 
	\def \d{     5.00 } 
	\def \h{     7.00 } 
	
	\coordinate (A) at (0,0); 
	\coordinate (B) at ({-\h*sin(\tta)},{\h*cos(\tta)}); 
	\coordinate (C) at ({-\h*sin(\tta)-\d*sin(\k*\tta)},
	{\h*cos(\tta)+\d*cos(\k*\tta)}); 
	\coordinate (D) at ({-\d*sin(\k*\tta)},{\d*cos(\k*\tta)}); 
	
	\coordinate (Ap) at (\l,0); 
	\coordinate (Bp) at ({\l-\h*sin(\tta)},{\h*cos(\tta)}); 
	\coordinate (Cp) at ({\l-\h*sin(\tta)-\d*sin(\k*\tta)},
	{\h*cos(\tta)+\d*cos(\k*\tta)}); 
	\coordinate (Dp) at ({\l-\d*sin(\k*\tta)},{\d*cos(\k*\tta)}); 
	
	\fill[red]  (A) circle [radius=2pt]; 
	\fill[red]    (B) circle [radius=2pt]; 
	\fill[black]  (C) circle [radius=2pt]; 
	\fill[black]    (D) circle [radius=2pt]; 
	\fill[black]   (Ap) circle [radius=2pt]; 
	\fill[red] (Bp) circle [radius=2pt]; 
	\fill[red]   (Cp) circle [radius=2pt]; 
	\fill[black] (Dp) circle [radius=2pt]; 
	
	
	\filldraw[draw=red,bottom color=blue!50!white, top color=blue!3!white]
	(B) -- (Bp)  -- (Cp);	
	\filldraw[draw=red,bottom color=blue!50!white, top color=blue!3!white]
	(A) -- (Cp) -- (Bp);
	\filldraw[draw=red,bottom color=blue!50!white, top color=blue!3!white, opacity=0.9]
	(A) -- (B) -- (Cp);

	\draw[red!80!white,-,very thick] (A) --  (Cp)
	(A) --  (B)
	(A) -- (Bp)
	(B)  --  (Cp)
	(B) --  (Bp)
	(Bp)  -- (Cp);
	
	\draw [-,thin] (B)  --  (A)
	(Ap) -- (Bp)
	(B)  --  (C)
	(D)  --  (C)
	(A)  --  (D)
	(Ap) --  (A)
	(Cp) --  (C)
	(Bp) --  (B)
	(Bp) -- (Cp);
	
	\draw [gray,-,thin] (Dp) -- (Cp)
	(D) --  (Dp)
	(Ap) -- (Dp);

	\draw (Ap) node [right,gray] {$B$}
	(Bp) node [right]           {$G$}
	(Cp) node [right]           {$H$}
	(C)  node [left,gray]       {$E$}
	(D)  node [left,gray]      {$D$}
	(A)  node [left]           {$A$}
	(B)  node [left]            {$F$}
	(Dp) node [right,gray]      {$C$};
	
	\fill[red]   (B) circle [radius=2pt]; 
	\fill[gray] (Dp) circle [radius=2pt]; 
	
	
\end{tikzpicture}
		\caption{The vertexes $A,F,G,H$ of $T_2$ as the corner points of the tetrahedron $C_{T_2}$.}
		\label{fig:tetra2}
	\end{subfigure}\vspace{1cm}\\
	\begin{subfigure}[t]{0.45\textwidth}
		\centering
		\begin{tikzpicture}[scale=0.42] 
	
	\def \tta{ -10.00 } 
	\def \k{    -3.00 } 
	\def \l{     6.00 } 
	\def \d{     5.00 } 
	\def \h{     7.00 } 
	
	\coordinate (A) at (0,0); 
	\coordinate (B) at ({-\h*sin(\tta)},{\h*cos(\tta)}); 
	\coordinate (C) at ({-\h*sin(\tta)-\d*sin(\k*\tta)},
	{\h*cos(\tta)+\d*cos(\k*\tta)}); 
	\coordinate (D) at ({-\d*sin(\k*\tta)},{\d*cos(\k*\tta)}); 
	
	\coordinate (Ap) at (\l,0); 
	\coordinate (Bp) at ({\l-\h*sin(\tta)},{\h*cos(\tta)}); 
	\coordinate (Cp) at ({\l-\h*sin(\tta)-\d*sin(\k*\tta)},
	{\h*cos(\tta)+\d*cos(\k*\tta)}); 
	\coordinate (Dp) at ({\l-\d*sin(\k*\tta)},{\d*cos(\k*\tta)}); 
	
	\fill[red]    (A) circle [radius=2pt]; 
	\fill[black]  (B) circle [radius=2pt]; 
	\fill[black]    (C) circle [radius=2pt]; 
	\fill[black]    (D) circle [radius=2pt]; 
	\fill[red]  (Ap) circle [radius=2pt]; 
	\fill[red] (Bp) circle [radius=2pt]; 
	\fill[red]   (Cp) circle [radius=2pt]; 
	\fill[black] (Dp) circle [radius=2pt]; 
	
	\filldraw[draw=red,bottom color=blue!50!white, top color=blue!3!white]
	(A) -- (Ap)  -- (Bp);
	\filldraw[draw=red,bottom color=blue!50!white, top color=blue!3!white]
	(A) -- (Bp) -- (Cp);	
	\filldraw[draw=red,bottom color=blue!50!white, top color=blue!3!white, opacity=1]
	(Ap) -- (Bp) -- (Cp);
	
	\draw[red!80!white,-,very thick] (A) -- (Cp)
	(A) -- (Bp)
	(A) --  (Ap)
	(Ap) -- (Cp)
	(Bp)  -- (Cp)
	(Ap) --  (Bp);
	
	\draw [-,thin] (B)  --  (A)
	(Ap) -- (Bp)
	(B)  --  (C)
	(D)  --  (C)
	(A)  --  (D)
	(Ap) --  (A)
	(Cp) --  (C)
	(Bp) --  (B)
	(Bp) -- (Cp);
	
	\draw [gray,-,thin] (Dp) -- (Cp)
	(D) --  (Dp)
	(Ap) -- (Dp);

	\draw (Ap) node [right] {$B$}
	(Bp) node [right]     {$G$}
	(Cp) node [right]          {$H$}
	(C)  node [left,gray]          {$E$}
	(D)  node [left,gray]         {$D$}
	(A)  node [left]           {$A$}
	(B)  node [left,gray]      {$F$}
	(Dp) node [right,gray]      {$C$};
	
	\fill[black]   (B) circle [radius=2pt]; 
	\fill[gray] (Dp) circle [radius=2pt]; 
	
	
\end{tikzpicture}
		\caption{The vertexes $A,B,G,H$ of $T_3$ as the corner points of the tetrahedron $C_{T_3}$.}
		\label{fig:tetra3}
	\end{subfigure}
	\hfill
	\begin{subfigure}[t]{0.45\textwidth}
		\centering
		\begin{tikzpicture}[scale=0.42] 
	
	\def \tta{ -10.00 } 
	\def \k{    -3.00 } 
	\def \l{     6.00 } 
	\def \d{     5.00 } 
	\def \h{     7.00 } 
	
	\coordinate (A) at (0,0); 
	\coordinate (B) at ({-\h*sin(\tta)},{\h*cos(\tta)}); 
	\coordinate (C) at ({-\h*sin(\tta)-\d*sin(\k*\tta)},
	{\h*cos(\tta)+\d*cos(\k*\tta)}); 
	\coordinate (D) at ({-\d*sin(\k*\tta)},{\d*cos(\k*\tta)}); 
	
	\coordinate (Ap) at (\l,0); 
	\coordinate (Bp) at ({\l-\h*sin(\tta)},{\h*cos(\tta)}); 
	\coordinate (Cp) at ({\l-\h*sin(\tta)-\d*sin(\k*\tta)},
	{\h*cos(\tta)+\d*cos(\k*\tta)}); 
	\coordinate (Dp) at ({\l-\d*sin(\k*\tta)},{\d*cos(\k*\tta)}); 
	
	\fill[red]    (A) circle [radius=2pt]; 
	\fill[black]  (B) circle [radius=2pt]; 
	\fill[black]    (C) circle [radius=2pt]; 
	\fill[red]    (D) circle [radius=2pt]; 
	\fill[black]  (Ap) circle [radius=2pt]; 
	\fill[black] (Bp) circle [radius=2pt]; 
	\fill[red]   (Cp) circle [radius=2pt]; 
	\fill[red] (Dp) circle [radius=2pt]; 
	
	\filldraw[draw=red,bottom color=blue!50!white, top color=blue!3!white]
	(A) -- (D) -- (Dp);
	\filldraw[draw=red,bottom color=blue!50!white, top color=blue!3!white]
	(A) -- (Dp) -- (Cp);
	\filldraw[draw=red,bottom color=blue!50!white, top color=blue!3!white]
	(D) -- (Dp)  -- (Cp);
	
	\draw[red!80!white,-,very thick] (A) -- (Cp)
	(A) -- (D)
	(A) --  (Dp)
	(D) -- (Cp)
	(Dp) --  (Cp)
	(D)  -- (Dp);
	
	\draw [-,thin] (B)  --  (A)
	(Ap) -- (Bp)
	(B)  --  (C)
	(D)  --  (C)
	(A)  --  (D)
	(Ap) --  (A)
	(Cp) --  (C)
	(Bp) --  (B)
	(Bp) -- (Cp);
	
	\draw [gray,-,thin] (Dp) -- (Cp)
	(D) --  (Dp)
	(Ap) -- (Dp);

	\draw (Ap) node [right,gray]    {$B$}
	(Bp) node [right,gray]     {$G$}
	(Cp) node [right]          {$H$}
	(C)  node [left,gray]      {$E$}
	(D)  node [left]      {$D$}
	(A)  node [left]           {$A$}
	(B)  node [left,gray]      {$F$}
	(Dp) node [right]          {$C$};
	
	\fill[black]   (B) circle [radius=2pt]; 
	\fill[red]    (Dp) circle [radius=2pt]; 
	
	
\end{tikzpicture}
		\caption{The vertexes $A,D,C,H$ of $T_4$ as the corner points of the tetrahedron $C_{T_4}$.}
		\label{fig:tetra4}
	\end{subfigure}
	\vspace{1cm}\\
	\begin{subfigure}[t]{0.45\textwidth}
		\centering
		\begin{tikzpicture}[scale=0.42] 
	
	\def \tta{ -10.00 } 
	\def \k{    -3.00 } 
	\def \l{     6.00 } 
	\def \d{     5.00 } 
	\def \h{     7.00 } 
	
	\coordinate (A) at (0,0); 
	\coordinate (B) at ({-\h*sin(\tta)},{\h*cos(\tta)}); 
	\coordinate (C) at ({-\h*sin(\tta)-\d*sin(\k*\tta)},
	{\h*cos(\tta)+\d*cos(\k*\tta)}); 
	\coordinate (D) at ({-\d*sin(\k*\tta)},{\d*cos(\k*\tta)}); 
	
	\coordinate (Ap) at (\l,0); 
	\coordinate (Bp) at ({\l-\h*sin(\tta)},{\h*cos(\tta)}); 
	\coordinate (Cp) at ({\l-\h*sin(\tta)-\d*sin(\k*\tta)},
	{\h*cos(\tta)+\d*cos(\k*\tta)}); 
	\coordinate (Dp) at ({\l-\d*sin(\k*\tta)},{\d*cos(\k*\tta)}); 
	
	\fill[red]  (A) circle [radius=2pt]; 
	\fill[red]    (B) circle [radius=2pt]; 
	\fill[red]  (C) circle [radius=2pt]; 
	\fill[black]    (D) circle [radius=2pt]; 
	\fill[black]   (Ap) circle [radius=2pt]; 
	\fill[black] (Bp) circle [radius=2pt]; 
	\fill[red]   (Cp) circle [radius=2pt]; 
	\fill[black] (Dp) circle [radius=2pt]; 
	
	\filldraw[draw=red,bottom color=blue!50!white, top color=blue!3!white]
	(A) -- (B) -- (C);
	\filldraw[draw=red,bottom color=blue!50!white, top color=blue!3!white]
	(A) -- (B) -- (Cp);
	\filldraw[draw=red,bottom color=blue!50!white, top color=blue!3!white]
	(C) -- (B)  -- (Cp);
	\draw[red!80!white,-,very thick] (A) --  (Cp)
	(A) --  (B)
	(A) -- (C)
	(C)  --  (Cp)
	(B) --  (C);
	
	\draw [-,thin] (B)  --  (A)
	(Ap) -- (Bp)
	(B)  --  (C)
	(D)  --  (C)
	(A)  --  (D)
	(Ap) --  (A)
	(Cp) --  (C)
	(Bp) --  (B)
	(Bp) -- (Cp);
	
	\draw [gray,-,thin] (Dp) -- (Cp)
	(D) --  (Dp)
	(Ap) -- (Dp);

	\draw (Ap) node [right,gray] {$B$}
	(Bp) node [right, gray]           {$G$}
	(Cp) node [right]           {$H$}
	(C)  node [left]       {$E$}
	(D)  node [left,gray]      {$D$}
	(A)  node [left]           {$A$}
	(B)  node [left]            {$F$}
	(Dp) node [right,gray]      {$C$};
	
	\fill[red]   (B) circle [radius=2pt]; 
	\fill[gray] (Dp) circle [radius=2pt]; 
	
	
\end{tikzpicture}
		\caption{The vertexes $A,F,E,H$ of $T_5$ as the corner points of the tetrahedron $C_{T_5}$.}
		\label{fig:tetra5}
	\end{subfigure}
	\hfill
	\begin{subfigure}[t]{0.45\textwidth}
		\centering
		\begin{tikzpicture}[scale=0.42] 
	
	\def \tta{ -10.00 } 
	\def \k{    -3.00 } 
	\def \l{     6.00 } 
	\def \d{     5.00 } 
	\def \h{     7.00 } 
	
	\coordinate (A) at (0,0); 
	\coordinate (B) at ({-\h*sin(\tta)},{\h*cos(\tta)}); 
	\coordinate (C) at ({-\h*sin(\tta)-\d*sin(\k*\tta)},
	{\h*cos(\tta)+\d*cos(\k*\tta)}); 
	\coordinate (D) at ({-\d*sin(\k*\tta)},{\d*cos(\k*\tta)}); 
	
	\coordinate (Ap) at (\l,0); 
	\coordinate (Bp) at ({\l-\h*sin(\tta)},{\h*cos(\tta)}); 
	\coordinate (Cp) at ({\l-\h*sin(\tta)-\d*sin(\k*\tta)},
	{\h*cos(\tta)+\d*cos(\k*\tta)}); 
	\coordinate (Dp) at ({\l-\d*sin(\k*\tta)},{\d*cos(\k*\tta)}); 
	
	\fill[red]    (A) circle [radius=2pt]; 
	\fill[black]  (B) circle [radius=2pt]; 
	\fill[red]    (C) circle [radius=2pt]; 
	\fill[red]    (D) circle [radius=2pt]; 
	\fill[black]  (Ap) circle [radius=2pt]; 
	\fill[black] (Bp) circle [radius=2pt]; 
	\fill[red]   (Cp) circle [radius=2pt]; 
	\fill[black] (Dp) circle [radius=2pt]; 
	
	\filldraw[draw=red,bottom color=blue!50!white, top color=blue!3!white]
	(A) -- (D)  -- (C);
	\filldraw[draw=red,bottom color=blue!50!white, top color=blue!3!white]
	(A) -- (D) -- (Cp);
	\filldraw[draw=red,bottom color=blue!50!white, top color=blue!3!white]
	(D) -- (C) -- (Cp);
	
	\draw[red!80!white,-,very thick] (A) -- (Cp)
	(A) -- (C)
	(A) --  (D)
	(D) -- (C)
	(C)  -- (Cp)
	(D) --  (Cp);
	
	\draw [-,thin] (B)  --  (A)
	(Ap) -- (Bp)
	(B)  --  (C)
	(D)  --  (C)
	(A)  --  (D)
	(Ap) --  (A)
	(Cp) --  (C)
	(Bp) --  (B)
	(Bp) -- (Cp);
	
	\draw [gray,-,thin] (Dp) -- (Cp)
	(D) --  (Dp)
	(Ap) -- (Dp);

	\draw (Ap) node [right,gray] {$B$}
	(Bp) node [right,gray]     {$G$}
	(Cp) node [right]          {$H$}
	(C)  node [left]          {$E$}
	(D)  node [left]         {$D$}
	(A)  node [left]           {$A$}
	(B)  node [left,gray]      {$F$}
	(Dp) node [right,gray]      {$C$};
	
	\fill[black]   (B) circle [radius=2pt]; 
	\fill[gray] (Dp) circle [radius=2pt]; 
	
	
\end{tikzpicture}
		\caption{The vertexes $A,D,E,H$ of $T_6$ as the corner points of the tetrahedron $C_{T_6}$.}
		\label{fig:tetra6}
	\end{subfigure}
	\caption{The Coxeter-Freudenthal-Kuhn triangulation of the unit cube for $d=3$.}
	\label{fig:tetra}
\end{figure}

We now illustrate the constructive idea behind the Coxeter-Freudenthal-Kuhn triangulation of the unit cube for the case $d=3$.
We denote the corner points by
\begin{gather*}
	A\,=\,(0,0,0),\qquad B\,=\,(1,0,0),\qquad C\,=\,(1,1,0),\qquad D\,=\,(0,1,0),\\
	E\,=\,(0,1,1),\qquad F\,=\,(0,0,1),\qquad G\,=\,(1,0,1),\qquad H\,=\,(1,1,1).
\end{gather*}
For an element $T\in\mathcal T_1^0$ we write $T:T(0)\to T(1)\to T(2)\to T(3)$ and
for a permutation $\sigma\in\mathscr S_3$ we write $\sigma=(\sigma(1),\sigma(2),\sigma(3))$.
There are six elements in the symmetric group $\mathscr S_3$, hence six elements in $\mathscr T_1^0$. Those are
\begin{align*}
	&T_1:\, A\to B\to C\to H\qquad\textnormal{with}\qquad \sigma_{T_1}=(1,2,3),\\
	&T_2:\, A\to F\to G\to H\qquad\textnormal{with}\qquad \sigma_{T_2}=(3,1,2),\\
	&T_3:\, A\to B\to G\to H\qquad\textnormal{with}\qquad \sigma_{T_3}=(1,3,2),\\
	&T_4:\, A\to D\to C\to H\qquad\textnormal{with}\qquad \sigma_{T_4}=(2,1,3),\\
	&T_5:\, A\to F\to E\to H\qquad\textnormal{with}\qquad \sigma_{T_5}=(3,2,1),\\
	&T_6:\, A\to D\to E\to H\qquad\textnormal{with}\qquad \sigma_{T_6}=(2,3,1).
\end{align*}
\cref{fig:tetra1,fig:tetra2,fig:tetra3,fig:tetra4,fig:tetra5,fig:tetra6} capture the set $C_T$ for each $T\in\mathscr T_1^0$. 
The method used in \cref{fig:tetra} visualizes the respective tetrahedron by colouring three of the four triangles which form its surface.

In this article, however, it is advantageous for technical reasons to have a partition of the unit cube $[0,1)^d$.
As the polyhedrons $C_T$, $T\in\mathscr T_1^0$, do intersect on their boundary, we now slightly modify the sets to obtain
a $d!$-sized family of sets $D_T$, indexed by $T\in\mathscr T_1^0$, with $D_T\subset C_T$ and 
\begin{equation}\label{eqn:dartition}
	[0,1)^d\,=\,\dot{\bigcup_{T\in{\mathscr T}_1^0}}D_T.
\end{equation}
The symbol $<_{\sigma,i}$ for $i\in\{2,\dots,d\}$ and $\sigma\in \mathscr S_d$ denotes the relation on $\R$ which coincides with `$<$' in case $\sigma(i-1)<\sigma(i)$ and 
with `$\leq$' in case $\sigma(i-1)>\sigma(i)$. Now we define
\begin{equation*}
	D_T\,:=\,\Big\{\,x\in\R^d\,\Big|\,0\leq x_{\sigma_T(d)}<_{\sigma_T,d}x_{\sigma_T(d-1)}<_{\sigma_T,d-1}\dots<_{\sigma_T,2}x_{\sigma_T(1)}<1\,\Big\}.
\end{equation*}
Obviously, the topological closure $\overline {D_T}$ coincides with the convex hull $C_T$ of $\{T(i)|i=0,\dots,d\}$.
Before we give a quick argumentation why \cref{eqn:dartition} is satisfied with this definition, we remark
that for reasons of symmetry the volume $|D_T|$ equals $1/(d!)$ under \cref{eqn:dartition}.
Let now $x\in [0,1)^d$. To find the unique element $T\in\mathscr T_1^0$ such that $x\in D_T$
one simply takes the lexicographically ordered sequence, say $P_1>\dots>P_d$, of the tuples $\{(x_i,i)|i=1,\dots,d\}$. 
Then one defines $\sigma(1)$ to be the second component of $P_1$, $\sigma(2)$ to be the second component of $P_2$, and so forth.
It holds $x\in D_T$ if and only if $\sigma_T=\sigma$. 

With the help of the sets $D_T$ we now construct certain piecewise linear, continuous functions on $\R^d$, which turn out useful in the context of Mosco convergence.
For the moment $T\in\mathscr T_1^0$ is fixed.
Let $i\in\{0,\dots,d\}$. \newpage\noindent The hyperplane in $\R^{d+1}$ which interpolates the sample $(T(j),\eins_{\{i\}}(j))$, $j=0,\dots,d$, can be represented as the graph of
the function
\begin{equation*}
	H_T^i:\,\R^d\,\ni\, x\,\longmapsto\,
	\begin{cases}
		\,1-x_{\sigma_T(1)}&\textnormal{if }i=0,\\
		\,x_{\sigma_T(i)}-x_{\sigma_T(i+1)}&\textnormal{if }i\in\{1,\dots,d-1\}\textnormal{ and}\\
		\,x_{\sigma_T(d)}&\textnormal{if }i=d.
	\end{cases}
\end{equation*}
Indeed, given $j\in\{0,\dots,d\}$ and $k\in\{1,\dots,d\}$, the $\sigma_T(k)$-th component of $T(j)$ equals $1$ if $j\geq k$,
while the $\sigma_T(k)$-th component of $T(j)$ equals $0$ if $j<k$, due to \cref{eqn:suniv}. Hence, we verify
\begin{equation*}
	H_T^i(T(j))\,=\,
	\begin{cases}
		\,1-\eins_{\{1,\dots,d\}}(j)\,&=\,\eins_{\{i\}}(j)\quad\textnormal{if }i=0,\\
		\,\eins_{\{i,\dots,d\}}(j)-\eins_{\{i+1,\dots,d\}}(j)&=\,\eins_{\{i\}}(j)\quad\textnormal{if }i\in\{1,\dots,d-1\},\\
		\,\eins_{\{d\}}(j)&=\,\eins_{\{i\}}(j)\quad\textnormal{if }i=d.
	\end{cases}
\end{equation*}
In particular, it holds
\begin{equation}\label{eqn:partun}
	\sum_{i=0}^d H_T^i\,=\,\eins_{\R^d}.
\end{equation}
For $i\in\{0,\dots,d\}$ the gradient of $H_T^i$ is the constant vector
\begin{equation}\label{eqn:gradha}
	\sum_{k=1}^d\partial_kH^i_T\,\univ_k\,=\,\eins_{\{1,\dots,d\}}(i)\,\univ_{\sigma_T(i)}-\eins_{\{0,\dots,d-1\}}(i)\,\univ_{\sigma_T(i+1)}.
\end{equation}
Using \cref{eqn:gradha} to calculate
the euclidean scalar product of the gradients of $H_T^i$ and $H_T^j$ for $i,j\in\{0,\dots,d\}$ at a point $x\in\R^d$ we obtain
\begin{equation}\label{eqn:gaha}
	\sum_{k=1}^d\partial_kH^i_T(x)\,\partial_kH^j_T(x)\,=\,
	\begin{cases}
		\,\phantom{-}2&\textnormal{if }i=j\in\{1,\dots,d-1\},\\
		\,\phantom{-}1&\textnormal{if }i=j\in\{0,d\},\\
		\,-1&\textnormal{if }|i-j|=1\textnormal{ and}\\
		\,\phantom{-}0&\textnormal{else.}
	\end{cases}
\end{equation}

\begin{figure}[t]
	\begin{subfigure}[t]{0.35\textwidth}
				\begin{tikzpicture}[scale=1.4]
			\draw[thick, gray, fill=ForestGreen!60!white](-2,-1)--(-2,0)--(-1,0)--(-2,-1);
			\draw[thick, gray, fill=ForestGreen!60!white](-1,-1)--(-1,0)--(0,0)--(-1,-1);
			\draw[thick, gray, fill=ForestGreen!60!white](-1,0)--(-1,1)--(0,1)--(-1,0);
			
			\draw[thick, gray, fill=LimeGreen!80!white](-1,0)--(0,0)--(0,1)--(-1,0);
			\draw[thick, gray, fill=LimeGreen!80!white](-2,0)--(-1,0)--(-1,1)--(-2,0);
			\draw[thick, gray, fill=LimeGreen!80!white](-2,-1)--(-1,-1)--(-1,0)--(-2,-1);
			
			\draw[thick,->] (-1,0) -- (0.7,0) node[anchor=north east]{$x$};
			\draw[thick,->] (-1,0) -- (-1,1.7) node[anchor=north east]{$y$};
			\draw (0,0.1) -- (0,-0.1) node[below] at (0.1,0) {$1$};
			\draw (-0.9,1) -- (-1.1,1) node[left] at (-1,1.1) {$1$};			
		\end{tikzpicture}
		\caption{The six triangular domains comprising the support of $\chi_1^0$.}
		\label{fig:trian}
	\end{subfigure}
	\hfill
	\begin{subfigure}[t]{0.6\textwidth}
		\tdplotsetmaincoords{65}{10}

\begin{tikzpicture}[tdplot_main_coords,scale=0.9]
	\draw[thick,->] (-4,2,0) -- (-3,2,0) node[anchor=north east]{$x$};
	\draw[thick,->] (-4,2,0) -- (-4,3,0) node[anchor=north west]{$y$};
	\draw[thick,->] (-4,2,0) -- (-4,2,1) node[anchor=south]{$z$};

	\draw[fill=green!60!white] (-2,0,0)--(-2,-2,0)--(0,-2,0)--(2,0,0)--(2,2,0)--(0,2,0)--(-2,0,0);

	\draw[thick] (-2,0,0)--(-2,-2,0);
	\draw[thick] (-2,-2,0)--(0,-2,0);
	\draw[thick] (0,-2,0)--(2,0,0);
	\draw[dashed] (-2,0,0)--(0,2,0);
	\draw[dashed] (0,2,0)--(2,2,0);
	\draw[thick] (2,2,0)--(2,0,0);
	
	\draw[thick,gray] (-2,2,0)--(-2,-4,0);
	\draw[thick,gray] (-4,-2,0)--(2,-2,0);
	\draw[thick,gray] (-2,-4,0)--(4,2,0);
	\draw[thick,gray] (-3,-3,0)--(3,3,0);
	\draw[thick,gray] (-4,-2,0)--(2,4,0);
	\draw[thick,gray] (-2,2,0)--(4,2,0);
	\draw[thick,gray] (2,4,0)--(2,-2,0);
	\draw[thick,gray] (0,4,0)--(0,-4,0);
	\draw[thick,gray] (-4,0,0)--(4,0,0);

	\draw[dashed,red] (0,0,0)--(0,0,3);

	\draw[thick,blue] (-2,0,0)--(0,0,2);
	\draw[thick,blue] (-2,0,0)--(-2,-2,0);
	\draw[thick,blue]  (-2,-2,0)--(0,0,2);
	\draw[thick,blue]  (-2,-2,0)--(0,-2,0);
	\draw[thick,blue]  (0,-2,0)--(0,0,2);
	\draw[thick,blue!50!white]   (0,2,0);
	\draw[thick,blue]  (2,2,0)--(0,0,2);
	\draw[thick,blue]  (2,2,0)--(2,0,0);
	\draw[thick,blue]  (2,0,0)--(0,0,2);
	\draw[thick,blue]  (2,0,0)--(2,2,0);
	\draw[thick,blue]  (2,0,0)--(0,-2,0);
	
	\draw[thick,blue]  (-2,0,0)--(0,2,0);
	\draw[thick,blue]  (0,2,0)--(2,2,0);

	\draw[red] (0,0,0) node {$\bullet$};

	\path[fill=blue!70!white,opacity=.4]
	(-2,0,0)--(0,0,2)--(-2,-2,0);
	
	\path[fill=blue!30!white,opacity=.4]
	(-2,-2,0)--(0,0,2)--(0,-2,0);
	
	\path[fill=blue!25!white,opacity=.4]
	(0,-2,0)--(0,0,2)--(2,0,0);
	
	\path[fill=blue!40!white,opacity=.4]
	(2,0,0)--(0,0,2)--(2,2,0);
	
	\path[fill=blue!30!white,opacity=.4]
	(2,2,0)--(0,0,2)--(0,2,0);
	
	\path[fill=blue!30!white,opacity=.4]
	(-2,0,0)--(0,0,2)--(0,2,0);
	
	\node[right] at (0.1,-0.5,0) {$0$};
	
\end{tikzpicture}
		\caption{The graph of $\chi_1^0$ in the $(2+1)$-dimensional space.}\label{fig:tent}
	\end{subfigure}
	\caption{The primal tent function $\chi_1^0$ for $d=2$ with its hexagonal support.}
\end{figure}
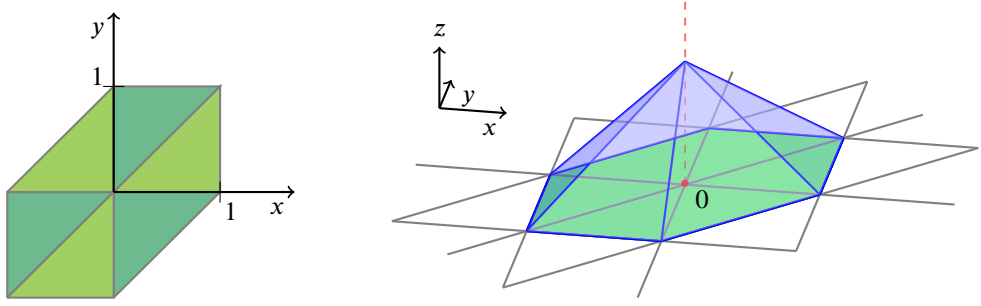

If we compose the function $\eins_{D_T}\cdot H_T^i$ with the shift $\R^d\ni x\mapsto x+T(i)$ and sum up 
over all $T\in\mathscr T_0^1$ and $i=0,\dots,d$ we arrive at the definition of the the \textbf{primal tent function}
\begin{equation*}
	\chi_1^0:\R^d\ni x\,\mapsto\,\sum_{T\in\mathscr T_1^0}\sum_{i=0}^{d}\eins_{D_T}(x+T(i))\,H_T^i(x+T(i)).
\end{equation*}
This piecewise definition glues together $(d+1)\,|\mathscr T_1^0|=(d+1)!$ many components, all of which are cut-off linear functions. 
For $T\in\mathscr T_1^0$ and $i=0,\dots,d$ the
indicator function $\eins_{D_T}(\,\cdot\,+T(i))$ $=\eins_{D_T-T(i)}(\cdot)$ - up to a set of Lebesgue measure zero - weights the convex hull of the points
\begin{equation*}
	P_0:=T(0)-T(i),\,P_1:=T(1)-T(i),\,\dots,\, P_d:=T(d)-T(i).
\end{equation*}
$P_0,\dots,P_d$ are the vertexes of a path of length $(d+1)$, which only travels in directions parallel to $\univ_1,\univ_2\dots,$ or $\univ_d$, uses each direction once, 
and visits the origin as its $i$-th vertex. The function $H_T^i(\,\cdot\,+T(i))$ is the linear interpolation of the sample ${(P_j,\delta_{\{i\}}(j))}$, $j=0,\dots,d$,
being $1$ in the origin and $0$ at all other nodes $P_j$, $j\neq i$.
\cref{fig:tent} shows the graph of the primal tent function for the case $d=2$. Its support, highlighted in \cref{fig:trian}, comprises six triangular domains.
\newpage\noindent These are
\begin{align*}
	\overline D_{T_1},\qquad\qquad &\phantom{-\univ_1-\univ_2+}\overline D_{T_2},\\
	-\univ_1+\overline D_{T_1},\qquad\qquad &\phantom{-\univ_2}-\univ_2+\overline D_{T_2},\\
	-\univ_1-\univ_2+\overline D_{T_1},\qquad\qquad &-\univ_1-\univ_2+\overline D_{T_2},
\end{align*}
where $T_1,T_2$ are the two elements of $\mathscr T_1^0$ in the two-dimensional case, defined via \cref{eqn:suniv} with $\sigma_{T_1}=(1,2)$ and $\sigma_{T_2}=(2,1)$.

The fact that $\chi_1^0$ is indeed an element of $\mathscr C$, i.e.~a primal function, becomes evident in the proceeding part.
We define the \textbf{tent function} with scaling parameter $r\in(0,\infty)$ and node $\alpha\in r\,\Z^d$ as $\chi_r^\alpha(x):=\chi_1^0((x-\alpha)/r)$, $x\in\R^d$.    
For $x\in\R^d$ it holds
\begin{align}\label{eqn:partuna}
	\sum_{\alpha\in r\Z^d}\chi_r^\alpha(x)\,
	=\,\sum_{\beta\in\Z^d}\,\sum_{T\in\mathscr T_1^0}\,\sum_{i=0}^{d}\,\eins_{D_T}\big(\,\frac{x}{r}+T(i)-\beta\,\big)\,H_T^i\big(\,\frac{x}{r}+T(i)-\beta\,\big).
\end{align}
Since $D_T\subset[0,1)^d$, a summand of the right hand side can only yield a non-zero value if
\begin{equation*}
	\beta\,=\,\big[\,\frac{x}{r}+T(i)\,\big]\,=\,\big[\,\frac{x}{r}\,\big]+T(i).
\end{equation*}
Hence, with \cref{eqn:partun} and \cref{eqn:dartition} we continue the calculation for the value of \cref{eqn:partuna} by
\begin{align}\label{eqn:partunb}
	\sum_{\alpha\in r\Z^d}\chi_r^\alpha(x)\,&=\,
	\sum_{T\in\mathscr T_1^0}\,\sum_{i=0}^{d}\,\eins_{D_T}\big(\,\frac{x}{r}-\big[\,\frac{x}{r}\,\big]\,\big)\,H_T^i\big(\,\frac{x}{r}-\big[\,\frac{x}{r}\,\big]\,\big)\nonumber\\
	&=\,\sum_{T\in\mathscr T_1^0}\,\eins_{D_T}\big(\,\frac{x}{r}-\big[\,\frac{x}{r}\,\big]\,\big)
	\,=\,1
\end{align}
for $x\in\R^d$.   
So, the family $\chi_r^\alpha$, $\alpha\in r\,\Z^d$ form a partition of unity for any $r\in(0,\infty)$. 
We sum up the properties which make it valuable for the proceeding discussion in \cref{thm:partun}.
Before we state the theorem we introduce some further terminology. The rescaled family $r\,D_T$, $T\in\mathscr T_1^0$,
partition the semi-open cube $[0,r)^d$ with side length $r$ for fixed $r\in(0,\infty)$. Hence the family $\alpha+r\,D_T$, over $\alpha\in r\,\Z^d$, $T\in\mathscr T_1^0$,
form a partition of $\R^d$. For $r\in(0,\infty)$ we set
\begin{equation*}
	\mathscr T_r\,:=\,\bigcup_{\alpha\in r\Z^d}\alpha+r\,\mathscr T^0_1
\end{equation*}
and if $T=\alpha+r\, T'\in\mathscr T_r$ with $T'\in\mathscr T_1^0$, $\alpha\in r\,\Z^d$, then
\begin{equation*}
	D_T\,:=\,\alpha+r\,D_{T'},\qquad
	\sigma_T\,:=\,\sigma_{T'}.
\end{equation*}
\begin{thm}\label{thm:partun}
	Let $r\in(0,\infty)$.
	The space $\R^d$ admits a partition $\{D_T|T\in \mathscr T_r\}$ and a partition of unity ${(\chi_r^{\alpha})}_{\alpha\in r\Z^d}$ with the following properties.
	\begin{thmlist}[wide, labelindent=0pt]
		\item \label{eqn:partuni}$\int_{\R^d}\chi_r^{\alpha} \de x =r^d$ for $\alpha\in r\Z^{d}$.
		\item\label{eqn:partunii} For $\alpha\in r\Z^d$ the function $\chi_r^{\alpha}$ is continuous on $\R^d$ with values in $[0,1]$ and support on $\alpha+[-r,r]^d$. 
		\item\label{eqn:partuniii} Let
		$w\in\R^{(r \Z^d)}$. The $i$-th weak partial derivative of the weighted sum $\sum_{\alpha}w_\alpha\,\chi_r^\alpha$ reads\label{thm:partuniii}
		\begin{equation*}
			\sum_{\alpha\in r\Z^d}w_\alpha\partial_i\chi^\alpha_r\,=\,\frac{1}{r}\sum_{\alpha\in r\Z^d}
			\sum_{\substack{T\in\mathscr T_r:\\T(\sigma^{-1}_T(i)-1)=\alpha}}(w_{\alpha+r\univ_i}-w_{\alpha})\,\eins_{D_T}
		\end{equation*}
		for $i=1,\dots,d$. The equality in the line above holds in $a.e.$-sense w.r.t.~the Lebesgue measure on $\R^d$. Thus the
		function on the right hand side is a version of the weak partial derivative of the left hand side.
		The squared norm of the weak gradient calculates as
		\begin{equation*}
			\sum_{i=1}^{d}\Big(\sum_{\alpha\in r\Z^d}w_\alpha\,\partial_i\chi^\alpha_r\Big)^2
			\,=\,\frac{1}{r^2}\sum_{T\in\mathscr T_r}\eins_{D_T}\sum_{i=1}^d\big(\,w_{T(i)}-w_{T(i-1)}\,\big)^2.
		\end{equation*}
		Again, the equation holds in $a.e.$-sense and the function on the right hand side is a version of the left hand side.
	\end{thmlist}
\end{thm}
\newcommand{\partuniii}{\cref{eqn:partuniii}}
\begin{proof}
	To proof (i) we use the shift invariance of the Lebesgue measure, the transformation formula of the integral, \cref{eqn:dartition} and \cref{eqn:partun} to calculate
	\begin{align*}
		\int\limits_{\R^d}\chi_r^\alpha\de x\,&=\, \sum_{T\in\mathscr T_1^0}\sum_{i=0}^d\,\int\limits_{rD_T} H_T^i\big(\frac{x}{r}\big)\de x\\
		&=\, r^d\,\sum_{T\in\mathscr T_1^0}\sum_{i=0}^d\,\int\limits_{D_T} H_T^i(x)\de x
		\,=\,\sum_{T\in\mathscr T_1^0}r^d\int\limits_{\R^d} \eins_{D_T}\de x\,=\,r^d.
	\end{align*}	
	for $\alpha\in r\,\Z^d$.
	
	For Item (ii), we use and prove an auxiliary result characterizing the support of the primal tent function $\chi_1^0$. Let 
	\begin{align*}
		A:=\Big\{\,x\in\R^d\,\Big|\,&-1\leq x_k<1\quad\textnormal{ for }k\in\{1,\dots,d\},\\
		&-1\leq x_l-x_k < 1\quad\textnormal{for }k,l\in\{1,\dots,d\}\textnormal{ with } k<l\,\Big\}.
	\end{align*}
	For $x\in\R^d$ we write $x_\text{min}:=\min(\{x_k$ $|$ $k=1,\dots,d\})$ and $x_\text{max}:=\max(\{x_k$ $|$ $k=1,\dots,d\})$.
	The statement, we need as a preliminary, fixes a point $x\in\R^d$ and reads:  \textit{
		\begin{description}[font=\normalfont]
			\item[(P)]\quad There exists $i\in\{0,\dots,d\}$ and $T\in\mathscr T_1^0$ such that $x+T(i)\in D_T$
			if and only if $x\in A$.\\
			Moreover, if $x+T(i)\in D_T$, the choices of $T\in\mathscr T_1^0$ and $i\in\{0,\dots,d\}$ are unique and it holds:
			\begin{equation*}
				H_T^i(x+T(i))=\min\Big(\Big\{\,1+x_\textnormal{min}\,,\,1-x_\textnormal{max}\,,\, 1+x_\textnormal{min}-x_\textnormal{max}\,\Big\}\Big).
			\end{equation*}
	\end{description}}

	Before we set out to prove this statement, we argue why it implies Item (ii) right away. If (P) is true, then
	$\chi_1^0(x)=0$ for $x\in\R^d\setminus A$ holds by definition of $\chi_1^0$ and hence $\supp[\chi_1 ^0]$ $\subset[-1,1]$.  
	Moreover, for $x\in A$ the estimate
	\begin{equation*}
		0\leq \min\Big(\Big\{\,1+x_\textnormal{min}\,,\,1-x_\textnormal{max}\,,\, 1+x_\textnormal{min}-x_\textnormal{max}\,\Big\}\Big)\leq 1
	\end{equation*}
	is valid. Therefore, the function $\chi_1^0$ on $\R^d$ takes values in $[0,1]$ under the assumption of (P).
	The continuity of $\chi_1^0$ under (P) is seen as follows. For $x\in\R^d\setminus A$ it holds
	\begin{equation*}
		\min\Big(\Big\{\,1+x_\textnormal{min}\,,\,1-x_\textnormal{max}\,,\, 1+x_\textnormal{min}-x_\textnormal{max}\,\Big\}\Big)\leq 0.
	\end{equation*}
	Consequently, 
	\begin{align*}
		\chi_1^0(x)&=\eins_A(x)\,\min\Big(\Big\{\,1+x_\textnormal{min}\,,\,1-x_\textnormal{max}\,,\, 1+x_\textnormal{min}-x_\textnormal{max}\,\Big\}\Big)\\
		&=\max\Big(\Big\{\,0\,,\, \min\Big(\Big\{\,1+x_\textnormal{min}\,,\,1-x_\textnormal{max}\,,\, 1+x_\textnormal{min}-x_\textnormal{max}\,\Big\}\Big)\,\Big\}\Big), \quad x\in\R^d.
	\end{align*}
	The expression on the right-hand-side, however, yields a continuous function on $\R^d$ in the variable $x$.
	Under the assumption of (P), we have shown the claim of Item (ii) for the case $\alpha=0$ and $r=1$. Then, Item (ii) in the general case is clear by the definition of $\chi_r^\alpha$ for $\alpha\in r\,\Z^d$ and $r\in(0,\infty)$. 
	
	Only the preliminary statement (P) is left to show concerning Item (ii). Let $x\in\R^d$. We start by assuming $x+T(i)\in D_T$ for some $i\in\{0,\dots,d\}$ and $T\in\mathscr T_1^0$.
	Since the sets $D_S$, $S\in\mathscr T_1^0$, form a partition of $[0,1)^d$, the condition $x+\alpha\in D_S$ for some $S\in\mathscr T_1^0$ and $\alpha\in\Z^d$ define $S$ and $\alpha$ uniquely.
	Hence, in our assumption, $T\in\mathscr T_1^0$ and $i\in\{0,\dots,d\}$ are uniquely determined. We now show that $x\in A$
	and prove that $H_T^i(x+T(i))$ admits the expression claimed in (P). The inequality $-1\leq x_k<1$ for $k\in\{1,\dots,d\}$ holds true, because of $(T(i))_k\in\{0,1\}$ and $x_k+(T(i))_k\in[0,1)$.
	To prove the inequalities 
	\begin{equation}\label{eqn:rev1}
		-1\leq x_l-x_k < 1\quad\textnormal{for }k,l\in\{1,\dots,d\}\textnormal{ with } k<l,
	\end{equation}
	we do a case distinction w.r.t.~the value of $i$. 
	
	First, if $i=0$, then $T(i)=0$ and therefore $x\in D_T$. 
	In particular, $x\in[0,1)^d$ and \cref{eqn:rev1} holds even with 
	strict inequalities on both sides. Furthermore, $x\in D_T$ implies $x_\textnormal{max}=x_{\sigma_T(1)}$. Since all components of $x$ are non-negative, we have 
	\begin{align*}
		\min\Big(\Big\{\,1+x_\textnormal{min}\,,\,1-x_\textnormal{max}\,,\, 1+x_\textnormal{min}-x_\textnormal{max}\,\Big\}\Big)&=1-x_\textnormal{max}=1-x_{\sigma_T(1)}
		\\&=H^{0}_T(x)=H^{0}_T(x+T(0)),
	\end{align*}
	as desired.
	
	In the next case, we assume $i=d$. Then $x+\univ \in D_T$ as $T(d)=\univ$. In particular, $x\in[-1,0)^d$ and \cref{eqn:rev1} again holds even with 
	strict inequalities on both sides. 
	From $x+\univ\in D_T$ it follows moreover $x_\textnormal{min}=x_{\sigma_T(d)}$. Consequently, as all components of $x$ are non-positive, we have
	\begin{align*}
		\min\Big(\Big\{\,1+x_\textnormal{min}\,,\,1-x_\textnormal{max}\,,\, 1+x_\textnormal{min}-x_\textnormal{max}\,\Big\}\Big)&=1+x_\textnormal{min}=1+x_{\sigma_T(d)}\\
		&=H^{d}_T(x+\univ)=H^{d}_T(x+T(d)),
	\end{align*}
	as desired for this case as well.
	
	Finally, we consider the case where $i\in\{1,\dots,d-1\}$. By virtue of \cref{eqn:suniv} it holds
	\begin{equation*}
		{(T(i))}_{\sigma_T(j)}=
		\begin{cases}0&\textnormal{if }i<j,\\
			1&\textnormal{if }i\geq j
		\end{cases}
	\end{equation*}
	for $j\in\{1,\dots,d\}$.
	So, the assumption $x+T(i)\in D_T$ is equivalent to the following system of inequalities:
	\begin{align}\label{eqn:rev2}
		0&\leq x_{\sigma_T(d)},\nonumber\\
		x_{\sigma_T(i+1)}&<_{\sigma_T,i+1}x_{\sigma_T(i)}+1,\nonumber\\
		x_{\sigma_T(j+1)}&<_{\sigma_T,j+1}x_{\sigma_T(j)}\quad\textnormal{for }j\in\{1,\dots,d-1\}\setminus\{i\},\nonumber\\
		x_{\sigma_T(1)}&<0.
	\end{align}
	The system of \cref{eqn:rev2} implies  
	\begin{equation*}
		x_\textnormal{max}=x_{\sigma_T(i+1)}\geq 0\quad\textnormal{and}\quad x_\textnormal{min}=x_{\sigma_T(i)}< 0.
	\end{equation*}
	Consequently,
	\begin{multline*}
		\min\Big(\Big\{\,1+x_\textnormal{min}\,,\,1-x_\textnormal{max}\,,\, 1+x_\textnormal{min}-x_\textnormal{max}\,\Big\}\Big)=1+x_\textnormal{min}-x_\textnormal{max}\\
		=1+x_{\sigma_T(i)}-x_{\sigma_T(i+1)}
		=x_{\sigma_T(i)}+1-\big(x_{\sigma_T(i+1)}+0\big)=H^{i}_T(x+T(i)),
	\end{multline*}
	as desired. We continue to prove \cref{eqn:rev1}.
	Since $x_\textnormal{max}-x_\textnormal{min}\leq 1$ is obvious from \cref{eqn:rev2}, it suffices to show that for any $k,l\in\{1,\dots,d\}$ the 
	equality $x_k-x_l=1$ necessitates $k<l$. To prove this, we first observe that 
	\begin{align*}
		x_{\sigma_T(j)}&\geq 0\quad\textnormal{for }j\in\{i+1,\dots,d\},\\
		x_{\sigma_T(j)}&< 0\quad\textnormal{for }j\in\{1,\dots,i\},
	\end{align*}
	holds as a consequence of the system \cref{eqn:rev2}.
	So, the equality $x_k-x_l=1$ implies $x_k=x_\textnormal{max}=x_{\sigma_T(i+1)}$ and $k=\sigma_T(j)$ for some $j\in\{i+1,\dots,d\}$ on the one hand.
	On the other hand, it implies $x_l=x_\textnormal{min}=x_{\sigma_T(i)}$ and $l=\sigma_T(j)$ for some $j\in\{1,\dots,i\}$.
	Then, due to \cref{eqn:rev2}, the equality $x_k-x_l=1$ leads to
	\begin{align*}
		x_{\sigma_T(i+1)}&=x_{\sigma_T(j)}\quad\textnormal{for }j=i+1,\dots,\sigma_T^{-1}(k),\\
		x_{\sigma_T(i+1)}&=x_{\sigma_T(i)}+1,\\
		x_{\sigma_T(i)}&=x_{\sigma_T(j)}\quad\textnormal{for }j=\sigma_T^{-1}(l),\dots,i.
	\end{align*}
	Furthermore, by definition of the relation $<_{\sigma_{T,j}}$ for $j\in\{2,\dots,d\}$, we must have
	\begin{equation*}
		\sigma_T(j)>\sigma_T(j+1)\quad\textnormal{for }j=\sigma_T^{-1}(l),\dots,\sigma_T^{-1}(k)-1.
	\end{equation*}
	This yields $k<l$. Summing up, we have shown that $x_\textnormal{max}-x_\textnormal{min}\leq 1$ and that the equality $x_k-x_l= 1$ for some $k,l\in\{1,\dots,d\}$ implies $k<l$.
	The proof of \cref{eqn:rev1} is now complete. 
	
	So far, one direction of the equivalence claimed in (P) has been proven, namely that
	$x+T(i)\in D_T$ for some $T\in\mathscr T_1^0$ and $i\in\{0,\dots,d\}$ implies  $x\in A$. Moreover, we verified for given $x\in\R^d$, that the claimed expression for the
	value of $H_T^i(x+T(i))$ holds true under the condition $x+T(i)\in D_T$ and that $T\in\mathscr T_1^0$ and $i\in\{0,\dots,d\}$ are uniquely determined by that condition.

	We now address the other direction of the equivalence claimed in (P). Let $x\in A$. Since $A\subset[-1,1)^d$
	there is a unique $\alpha\in\{0,1\}^d$ such that $x+\alpha\in [0,1)^d$, which is given by
	\begin{equation*}
		\alpha_k=
		\begin{cases}0&\textnormal{if }x_k\geq 0,\\
			1&\textnormal{if }x_k<0,
		\end{cases}
	\end{equation*} 
	for $k\in\{1,\dots,d\}$. Moreover, since the sets $D_T$, $T\in\mathscr T_1^0$, form a partition of $[0,1)^d$, there is a unique $T\in\mathscr T_1^0$ such that
	$x+\alpha\in D_T$.
	It suffices to show $\alpha=T(i)$ for some $i\in\{0,\dots,d\}$.
	We set $I^-:=\{k$ $|$ $k\in\{1,\dots,d\}, x_k<0\}$ and $I^+:=\{k$ $|$ $k\in\{1,\dots,d\}, x_k\geq 0\}$. 
	At first, we claim, if $l\in\sigma_T^{-1}(I^-)$ and $l\geq 2$ then $l-1$ is also an element of $\sigma_T^{-1}(I^-)$, as the converse leads to a contradiction.
	Let us assume $l\in\sigma_T^{-1}(I^-)$, $l\geq 2$, and $l-1\in \sigma_T^{-1}(I^+)$.
	Then, the condition $x+\alpha\in D_T$ yields
	\begin{equation*}
		x_{\sigma_T(l)}+1=x_{\sigma_T(l)}+\alpha_{\sigma_T(l)}<_{\sigma_T,l} x_{\sigma_T(l-1)}+\alpha_{\sigma_T(l-1)}=x_{\sigma_T(l-1)}.
	\end{equation*}
	Since $x\in A$ and hence $x_\textnormal{max}-x_\textnormal{min}\leq 1$, a strict inequality in the line above is impossible.
	Equality is only possible if $\sigma_T(l+1)<\sigma_T(l)$ by the definition of $<_{\sigma_T,l}$.
	However, this is a contradiction, as $x\in A$ would imply $x_{\sigma_T(l)}-x_{\sigma_T(l+1)}<1$.
	
	Now, let $i:=\max(\sigma_T^{-1}(I^-))$ in case $\sigma_T^{-1}(I^-)\neq \emptyset$ and $i:=0$ in case $\sigma_T^{-1}(I^-)=\emptyset$. 
	From the arguments above, it follows $\sigma_T^{-1}(I^-)=\{1,\dots,i\}$ if $i>0$. Consequently, for $j\in\{1,\dots,d\}$ we have
	\begin{equation*}
		{\alpha}_{\sigma_T(j)}=
		\begin{cases}0&\textnormal{if }i<j,\\
			1&\textnormal{if }i\geq j.
		\end{cases}
	\end{equation*}
	Hence, $\alpha=T(i)$ in view of \cref{eqn:suniv}. This concludes the proof of Item (ii).

	To show (iii) we fix a point $x\in\R^d$ and choose the unique element $T\in\mathscr T_r$, say $T=\gamma+rT'$ for some 
	$T'\in\mathscr T_1^0$, such that $x\in D_T=\gamma+r D_{T'}$. 
	We assume that $x$ is contained in the interior of $D_T$. This poses no restriction since the statements we want to show refer to an `a.e.'~way of reading and
	the boundary of $D_T$ has Lebesgue measure zero.
	Let $\alpha\in r\,\Z^d$.
	A necessary condition for $x\in\supp[\chi_r^\alpha]$ is
	\begin{equation*}
		x  \in \alpha-r\,\tilde T(j)+r\, D_{\tilde T}\qquad\textnormal{for some }\tilde T\in\mathscr T_1 ^0\textnormal{ and }j=0,\dots,d.
	\end{equation*}
	Since $\mathscr T_r$ yields a partition of $\R^d$, this condition is equivalent to 
	\begin{equation*}
		\tilde T=T'\qquad\textnormal{with }\gamma=\alpha-r\tilde T(j)\qquad\textnormal{for some }j=0,\dots,d,
	\end{equation*}
	which is in turn equivalent to
	$\tilde T=T'$ with $\alpha=\gamma+rT'(j)$ $(=T(j))$ for some $j=0,\dots,d$.
	In particular, $\{\alpha\in r\,\Z^d$ $|$ $x\in\textnormal{supp}[\chi_r^\alpha]\}$ $=$ $\{T(j)$ $|$ $j=0,1,\dots ,d\}$.
	Using these equivalences, we conclude that for each point $y$ from the interior of $D_T$ it holds
	\begin{align*}
		\chi^{T(j)}_r(y)\,=\,H^j_{T'}\big(\,\frac{y-T(j)}{r}+T'(j)\,\big)\,=\,H^j_{T'}\big(\,\frac{y-\gamma}{r}\,\big).
	\end{align*}
	Therefore, for $i\in\{1,\dots,d\}$, we compute, using \cref{eqn:gradha} 
	in the third equality, 
	\begin{multline}\label{eqn:weakpap}
		\sum_{\alpha\in r\Z^d}w_\alpha\,\partial_i\chi_r^\alpha(x)
		\,=\,\sum_{j=0}^dw_{T(j)}\,\frac{1}{r}\,\partial_i H^j_{T'}\big(\,\frac{x-\gamma}{r}\,\big)\\
		=\,\frac{1}{r}\sum_{j=0}^dw_{T(j)}\,{\Big(\,\eins_{\{1,\dots,d\}}(j)\,\univ_{\sigma_{T}(j)}
			-\eins_{\{0,\dots,d-1\}}(j)\,\univ_{\sigma_{T}(j+1)}\,\Big)}^\textnormal{T}\univ_i\\
		=\,\frac{1}{r}\sum_{j=1}^dw_{T(j)}\,\eins_{\{i\}}(\sigma_T(j))-\frac{1}{r}\sum_{j=0}^{d-1}w_{T(j)}\,\eins_{\{i\}}(\sigma_T(j+1))\\
		=\,\frac{1}{r}\,\big(\,w_{T(\sigma_T^{-1}(i))}-w_{T(\sigma_T^{-1}(i)-1)}\,\big)\\
		=\,\frac{1}{r}\,\big(\,w_{T(\sigma_T^{-1}(i)-1)+r\univ_i}-w_{T(\sigma_T^{-1}(i)-1)}\,\big).
	\end{multline}
	The last equality holds due to \cref{eqn:tirect}. Consequently,
	\begin{align*}
		\sum_{\alpha\in r\Z^d}w_\alpha\,\partial_i\chi_r^\alpha(x)\,
		&=\,\frac{1}{r}\,\sum_{T\in\mathscr T_r}\big(\,w_{T(\sigma_T^{-1}(i)-1)+r\univ_i}-w_{T(\sigma_T^{-1}(i)-1)}\,\big)\,\eins_{D_{T}}(x)\\
		&=\,\frac{1}{r}\,\sum_{\alpha\in r\Z^d}\,\sum_{\substack{T\in\mathscr T_r:\\T(\sigma_T^{-1}(i)-1)=\alpha}}
		\big(\,w_{T(\sigma_T^{-1}(i)-1)+r\univ_i}-w_{T(\sigma_T^{-1}(i)-1)}\,\big)\,\eins_{D_{T}}(x)
	\end{align*}
	as desired. To calculate the squared norm of the weak gradient, we sum up the squared values of \cref{eqn:weakpap} over the index $i=1,\dots d$
	and obtain
	\begin{multline*}
		\sum_{i=1}^{d}\Big(\sum_{\alpha\in r\Z^d}w_\alpha\,\partial_i\chi^\alpha_r\Big)^2(x)\,
		=\,\frac{1}{r^2}\,\sum_{T\in\mathscr T_r}\,\sum_{i=1}^d{\big(\,w_{T(\sigma_T^{-1}(i)-1)+r\univ_i}-w_{T(\sigma_T^{-1}(i)-1)}\,\big)}^2\,\eins_{D_{T}}(x)\\
		=\,\frac{1}{r^2}\,\sum_{T\in\mathscr T_r}\,\sum_{\substack{j\in\{1,\dots,d\}:\\\sigma_T(j)=i}}{\big(\,w_{T(j-1)+r\univ_{\sigma_T(j)}}-w_{T(j-1)}\,\big)}^2\,\eins_{D_{T}}(x)
		\\=\,\frac{1}{r^2}\,\sum_{T\in\mathscr T_r}\,\sum_{j=1}^d\big(\,w_{T(j)}-w_{T(j-1)}\,\big)^2\,\eins_{D_{T}}(x).\\
	\end{multline*}
	The last equality holds due to \cref{eqn:tirect} and concludes the proof.
\end{proof}

\subsection{\texorpdfstring{$L^2$}{L2} and energy estimates} \label{sec:est}
In this section we investigate the approximative quality of such weighted sums as have been considered in \partuniii{} regarding certain symmetric bilinear forms
of $L^2$ and energy type. Let $\varrho:\R^d\to[0,\infty)$ be a probability density w.r.t.~the Lebesgue measure on $\R^d$.
The energy
\begin{equation*}
	E^\varrho(f,g)\,:=\,\int\limits_{\R^d}\Gamma(f,g)\,\varrho\de x
	\qquad\textnormal{with}\qquad\Gamma(f,g)\,:=\,\sum_{i=1}^{d}\partial_if\,\partial_ig
\end{equation*}
shall be defined for $f,g\in\mathscr D$. We denote by
\begin{equation*}
	\mathscr D:= \Big\{\,f\in\cb(\R^d)\,\Big|\,\partial_i f\textnormal{ exists in weak sense and }\partial_if\in L^\infty(\R^d,\de x)\textnormal{ for }i=1,\dots,d\,\Big\}.
\end{equation*}
that linear subspace of the bounded, continuous functions $\cb(\R^d)$ whose elements are representatives for elements of the Sobolev space $H^{1,\infty}(\R^d)$.
The resulting Lemma is the key ingredient to the principle theorem of this paper, \cref{thm:m1} of the next section.
The approximative quality translates into conditions on the density $\varrho$ in terms of the $\varphi,\eta$-\textbf{residual} $R_{r}^{\varphi,\eta}(\varrho)$ and the 
$\varphi,\eta$-\textbf{perturbation} $I_{r}^{\varphi,\eta}(\varrho)$.
These quantities are defined depending on the primal functions $\varphi,\eta:\R^d\to[0,1]$ and the parameter $r\in(0,\infty)$, 
as mappings from the set of non-negative, integrable functions into itself.
\begin{align}
	R_{r}^{\varphi,\eta}(g)&\,:=\,
	\sum_{\alpha\in r\Z^d}r^{-d}\int\limits_{\R^d}|g(\,\cdot\,)-g(x)|\,\eta_r^\alpha(x)\de x\,\varphi_r^\alpha(\,\cdot\,)\label{eqn:phetares}\\
	\textnormal{and}\quad I_{r}^{\varphi,\eta}(g)&\,:=\,
	\sum_{\alpha\in r\Z^d}r^{-d}\int\limits_{\R^d}\eta_r^\alpha(x)\,g(x)\de x\,\varphi_r^\alpha(\,\cdot\,)\label{eqn:phetaper}
\end{align}
for a non-negative, measurable function $g\in\mplus(\R^d)$ with $\int_{\R^d}g\de x<\infty$. Since 
\begin{equation*}
	\int\limits_{\R^d}\eta_r^\alpha(x)\de x\,=\,r^d\qquad\textnormal{and}\qquad\int\limits_{\R^d}\varphi_r^\alpha(x)\de x\,=\,r^d,
\end{equation*}
it holds 
\begin{equation*}
	\int\limits_{\R^d}I_{r}^{\varphi,\eta}(g)\,(x)\de x\,=\,
	\sum_{\alpha\in r\Z^d}\,\int\limits_{\R^d}g(x)\,\eta_r^\alpha(x)\de x\,=\,
	\int\limits_{\R^d} g\de x
\end{equation*}
as well as
\begin{multline*}
	\int\limits_{\R^d}R_{r}^{\varphi,\eta}(g)\,(x)\de x\,
	\leq\,\int\limits_{\R^d}\,\sum_{\alpha\in r\Z^d}r^{-d}\int\limits_{\R^d}\big(g(y)+g(x)\big)\,\eta_r^\alpha(x)\de x\,\varphi_r^\alpha(y)\de y\\
	=\,\sum_{\alpha\in r\Z^d}\,\int\limits_{\R^d}g(y)\,\varphi_r^\alpha(y)\de y+\sum_{\alpha\in r\Z^d}\,\int\limits_{\R^d}g(x)\,\eta_r^\alpha(x)\de x
	\,=\,2\int\limits_{\R^d}g(x)\de x
\end{multline*}
for $g\in\mplus (\R^d)$. Let $m:=\varrho\,\de x$. 
The space of bounded measurable functions $\mb(\R^d)$ is dense in $L^2(\R^d,m)$. We use the convention $1/0=\infty$ and $\infty\cdot 0=0$.
For a measurable function $\kappa:\R^d\to[0,1]$ we set
\begin{align}
	\delta^\kappa_r(m)&\,:=\,\sup_{\varphi,\eta\in \mathscr C}\,\sup_{f\in\mb(\R^d)}
	\,\frac{|\int_{\R^d} f(x)\,R_{r}^{\varphi,\eta}(\kappa\,\varrho)\,(x)\de x\,|}{\|f\|_{L^2(m)}}\,\in [0,\infty]\label{eqn:delta}\\
	\textnormal{and}\quad C^\kappa_r(m)&\,:=\,\sup_{\varphi,\eta\in \mathscr C}\,\sup_{f\in\mb(\R^d)}
	\,\frac{|\int_{\R^d} f(x)\,I_{r}^{\varphi,\eta}(\kappa\,\varrho)\,(x)\de x\,|}{\|f\|_{L^1(m)}}\in [0,\infty].
	\label{eqn:bic}
\end{align}
A finite value for $C^\kappa_r(m)$ allows to extend the linear functional 
\begin{equation*}
	f\mapsto \int_{\R^d} \,f(x)\,I_{r}^{\varphi,\eta}(\kappa\,\varrho)(x)\de x
\end{equation*}
 to an element from the dual of $L^1(\R^d,m)$ with norm smaller equal $C^\kappa_r(m)$ via the BLT theorem.
In the same way, a finite value for $\delta^\kappa_r(m)$ allows to extend the linear functional $f\mapsto$  $\int_{\R^d} \,f(x)\,R_{r}^{\varphi,\eta}(\kappa\,\varrho)(x)\de x$ 
to an element from the dual
of $L^2(\R^d,m)$ with norm smaller equal $\delta^\kappa_r(m)$. We hint at a consequence of the Riesz isomorphism.
\begin{rem}\label{rem:blt}
	Let	$R\in\mplus(\R^d)$ with $\int_{\R^d}R\de x<\infty$.
	\begin{equation*}
		\sup_{f\in\mb(\R^d)}\,\frac{|\int_{\R^d} f(x)\,R(x)\de x\,|}{\|f\|_{L^2(m)}}\,=:\,D<\infty
	\end{equation*}
	if and only if there exists $\hat R\in L^2(\R^d,m)$ with $\|\hat R\|_{L^2(m)}=D$ and
	\begin{equation*}
		\int\limits_{\R^d} f(x)\,R(x)\de x\,=\,\int\limits_{\R^d} f(x)\,\hat R(x)\,\varrho(x)\,\de x.
	\end{equation*}
	for $f\in\mb(\R^d)$. The latter is equivalent to $R/\varrho\in L^2(m)$ with
	\begin{equation*}
		D\,:=\,\Big(\int\limits_{\R^d}|R(x)|^2\,\frac{1}{\varrho(x)}\de x\Big)^{\frac{1}{2}}.
	\end{equation*}
\end{rem}

To make the analysis useful for the purpose of the next section we consider direct integrals of such energies.
For a Polish space $S$  we denote the Borel  $\sigma$-algebra by $\mathcal B(S)$ and the bounded, measurable functions from $S$ to $\R$ by $\mb(S)$.
Let $(S,\mathcal B(S),\nu)$ be a probability space over a Polish space $S$. We assume
\begin{equation*}
	\de m_s(x)\,:=\,\varrho(s,x)\de x
\end{equation*}
is a probability measure on $\R^d$ with density $\varrho(s,\,\cdot\,)$ for $s\in S$.
Moreover, let $S\ni s\mapsto m_s(A)\in[0,1]$ be measurable for $A\in\mathcal B(\R^d)$.
The direct integral then reads
\begin{equation*}\label{eqn:dirint}
	\big(\textstyle\int^\oplus\displaystyle E^{\varrho(s,\cdot)}\de\nu\big)(u,v)\,:=\,\int\limits_S E^{\varrho(s,\cdot)}\big(u(s,\cdot),v(s,\cdot)\big)\de\nu(s)
\end{equation*}
for $u,v$ from its domain
\begin{align*}\label{eqn:domdiri}
	\dom\big(\textstyle{\int^\oplus} E^{\varrho(s,\cdot)}\de\nu\big)\,:=\Big\{\,u\in\mb(S\times\R^d)\,\Big|\,&u(s,\cdot)\in \mathscr D\textnormal{ for }s\in S\textnormal{ and}\nonumber\\
	&S\ni s\,\mapsto\,E^{\varrho(s,\cdot)}(u(s,\cdot),u(s,\cdot))\nonumber\\&\textnormal{ is measurable and integrable w.r.t.~}\nu\,\Big\}.
\end{align*} 

\noindent For given $r\in(0,\infty)$ we recall the family $\chi_r^\alpha$, $\alpha\in r\,\Z^d$, from \cref{sec:kuhn} and define a linear subspace of 
$\dom\big(\textstyle{\int^\oplus} E^{\varrho(s,\cdot)}\de\nu\big)$ by
\begin{align*}
	\mathscr L_r\,:=\,\Big\{\lambda:S\times \R^d\to\R\,\Big|\, 
	&\textnormal{there exists }M\in(0,\infty)\textnormal{ such that}\\
	&\begin{aligned}
		\lambda(s,x)=\sum_{\alpha\in r\Z^d}\lambda^\alpha(s)\,\chi_r^\alpha(x),\, s\in S,x\in\R^d,\,\textnormal{with}&\\
		\textnormal{measurable }\lambda^\alpha:S\to[-M,M]\textnormal{ for }\alpha\in r\,\Z^d&\,\Big\}.
	\end{aligned}
\end{align*}
Since the support of $\chi_r^\alpha$ is contained in the cube $\alpha+[-r,r]^d$
for every $\alpha\in r\,Z^d$ and $r\in(0,\infty)$, the number of indexes $\alpha$ for which, at a given point $x\in\R^d$, the function 
$\chi_r^\alpha$ does not vanish at $x$ is bounded by $2^d$. Therefore, $\mathscr L_r$ is a subspace of $\mb(S\times\R^d)$ for every $r\in(0,\infty)$.
The approximative qualities of the subspace $\mathscr L_r$ are in the focus of the next Lemma.
For $g\in\cc(\R^d)$ and $\varepsilon\in(0,\infty)$, the continuous functions on $\R^d$ with compact support, we write
\begin{equation}\label{eqn:modoco}
	\omega_\varepsilon^g\,:=\,\sup_{\substack{x,y\in\R^d\\\textnormal{max}_j|x_j-y_j|\leq 4\varepsilon}}|g(x)-g(y)|\,.
\end{equation}
Furthermore, if $g\in\cic(\R^d)$, the continuously differentiable functions on $\R^d$ with compact support, then
\begin{equation*}
	\omega_\varepsilon^{\nabla g}\,:=\,\max_{i=1,\dots,d}\,\omega_\varepsilon^{\partial_i g}
	\qquad\textnormal{and}\qquad
	D^g\,:=\max_{i=1,\dots,d}\|\partial_i g\|_\infty\,.
\end{equation*}

\begin{lem}\label{lem:polygon}
	Let $f\in\mb(S)$, $g\in\cic(\R^d)$ and $\kappa:S\times \R^d\to [0,1]$ be measurable.
	For $u\in\dom (\int^\oplus E^{\varrho(s,\cdot)}\de\nu)$ with $-1\leq u(\cdot)\leq 1$ and $r\in(0,\infty)$ there exists $\lambda\in \mathscr L_r$ 
	such that each of the inequalities holds true.\vspace{1ex}
	\begin{thmlist}[wide, labelindent=0pt]
		\item\label{eqn:polygon0} $|\lambda^\alpha(\,\cdot\,)|\,\leq\,1$
		where $\lambda^\alpha$ is the coefficient of $\lambda$ with index $\alpha\in r\,\Z^d$.\vspace{1ex}
		\item\label{eqn:polygon1} $\displaystyle\Big|\int\limits_{S\times\R^d}f(s)\,g(x)\,\big(u(s,x)-\lambda(s,x)\big)\,\kappa(s,x)\de m_s(x)\de\nu (s)\,\Big|$\vspace{1ex}\\
		$\phantom{(ii)}\,\displaystyle\leq\,\|f\|_{L^\infty(\nu)}\,\big(\omega_r^g+\|g\|_\infty\,\|\delta^{\kappa(s,\cdot)}_r(m_s)\|_{L^2(\nu)}\big)\,.$\vspace{1ex}
		\item\label{eqn:polygon2} 
		$\displaystyle \Big|\int\limits_Sf(s)\,E^{\kappa\varrho(s,\cdot)}\big(g,u(s,\cdot)-\lambda(s,\cdot)\big)\de\nu(s)\,\Big|$\vspace{1ex}\\
		$\phantom{(iii)}\,\displaystyle\leq\,\sqrt{d}\,\|f\|_{L^\infty(\nu)}\,\big(\omega_r^{\nabla g}+D^g\,\|\delta^{\kappa(s,\cdot)}_r(m_s)\|_{L^2(\nu)}\big)\,
		\Big(\int\limits_SE^{\varrho(s,\cdot)}\big(u(s,\cdot),u(s,\cdot)\big)\de\nu(s)\Big)^{\frac{1}{2}}\,.$\vspace{1ex}
		\item\label{eqn:polygon3} 
		$\displaystyle\int\limits_SE^{\kappa\varrho(s,\cdot)}\big(\lambda(s,\cdot),\lambda(s,\cdot)\big)\de\nu(s)$
		$\,\displaystyle\leq\,\|C^{\kappa(s,\cdot)}_r(m_s)\|_{L^\infty(\nu)}\,\int\limits_SE^{\varrho(s,\cdot)}\big(u(s,\cdot),u(s,\cdot)\big)\de\nu(s).$\vspace{1ex}
	\end{thmlist}
\end{lem}
\newcommand{\polygonz}{\cref{eqn:polygon0}}
\newcommand{\polygoni}{\cref{eqn:polygon1}}
\newcommand{\polygonii}{\cref{eqn:polygon3}}
\begin{proof}
	Let $r\in(0,\infty)$ be fixed throughout this proof.
	We start with an abstract estimate, which will be used in two separate instances in the subsequent course of this proof.
	By $\rho:\R^d\to[0,\infty)$ we denote a generic non-negative, measurable function with $\int\rho\de x\leq 1$.
	Let $\varphi$, $\eta$ denote two generic elements from $\mathscr C$ and $h\in\cc(\R^d)$. 
	Since the family ${(\varphi_r^\alpha)}_\alpha$ sum up to one and $r^{-d}\eta_r^\alpha$ is a probability density
	for $\alpha\in r\,\Z^d$, we estimate 
	\begin{multline*}\label{eqn:Absineq}
		\Big|\,h(x)\,\rho(x)-\sum_{\alpha\in r\Z^d}\varphi^\alpha_r(x)
		\,r^{-d}\int\limits_{\R^d}\eta_r^\alpha(y)\, h(y)\,\rho(y)\de y\,\Big|\\
		\begin{aligned}
			&=\,\Big|\sum_{\alpha\in r\Z^d}\varphi^\alpha_r(x)
			\,r^{-d}\int\limits_{\R^d}(h(x)\,\rho(x)- h(y)\,\rho(y))\,\eta_r^\alpha(y)\de y\Big|\\
			&\leq\,\sum_{\alpha\in r\Z^d}\varphi^\alpha_r(x)\,\Big(\omega_r^h\,\rho(x)+\|h\|_\infty\,r^{-d}\int\limits_{\R^d}|\rho(x)-\rho(y)|\,\eta_r^\alpha(y) \de y\Big)
			\end{aligned}\\
			 =\, \omega_r^h\,\rho(x)+\|h\|_\infty\,R^{\varphi,\eta}_r(\rho)(x)\numberthis
	\end{multline*}
	for each point $x\in\R^d$.
	\newcommand{\eqnAbsineq}{\cref{eqn:Absineq}}
	
	As in the claim of this lemma, we fix $u\in\dom (\int^\oplus E^{\varrho(s,\cdot)}\de\nu)\cap\mathcal M(S\times\R^d,[-1,1])$, $g\in\cic(\R^d)$ and $f\in\mb(S)$.
	The definition
	\begin{equation}\label{eqn:defla}
		\lambda(s,x)\,:=\,\sum_{\alpha\in r\Z^d}\frac{1}{r^d}\int\limits_{\R^d} u(s,y)\,\eins_{(-r,0]^d}(\alpha-y)\de y\,\chi_r^\alpha(x)
	\end{equation}
	for $x\in\R^d$ and $s\in S$ is in accordance with (i).
	We now dedicate ourselves to the verification of (ii) to (iv).
	
	We start with Item (ii). The integrals, which needs to be computed, involves only integrands of bounded functions. By linearity, we obtain
	\begin{multline*}\label{eqn:absineqii}
		\Big|\int\limits_{S\times\R^d}f(s)\,g(x)\,\big(u(s,x)-\lambda(s,x)\big)\,\kappa(s,x)\,\de m_s(x)\de\nu(s)\,\Big|\\
		\begin{aligned}
			&=\Big|\int\limits_{S\times\R^d}f(s)\,u(s,x)\,g(x)\,\kappa(s,x)\,\de m_s(x)\de\nu(s)\\
			&\phantom{=}-\int\limits_{S\times\R^d}f(s)\,u(s,x)
			\sum_{\alpha\in r\Z^d}\Big(\eins_{(-r,0]^d}(\alpha-x)\,
			r^{-d}\int\limits_{\R^d}g(y)\,\kappa(s,y)\,\chi_r^\alpha(y)\de m_s(y)\Big)\de x\de\nu(s)\Big|.
			\end{aligned}\\\numberthis
	\end{multline*}
	\newcommand{\eqnAbsineqii}{\cref{eqn:absineqii}}
	To get the subtracting term into the form as it is written in the line above, we plug in \cref{eqn:defla}, use Fubini's theorem and 
	after changing the order of integration we also exchange the names of the variables $x$ and $y$.  	
	To go on, we put the subtraction inside the integral again and for each $s\in S$ and $x\in\R^d$ use \eqnAbsineq~with the choices 
	$\kappa(s,\cdot)\,\varrho(s,\cdot)$ as $\rho$, $\eins_{[0,1)^d}$ as $\varphi$, the tent function $\chi_1^0$ as $\eta$ and $h=g$, to estimate 
	\eqnAbsineqii~from above with
	\begin{multline*}
		\int\limits_S\int\limits_{\R^d}|f(s)\,u(s,x)|\,\Big(\,\omega_r^g\,\kappa (s,x)\,\varrho(s,x)
		+\,\|g\|_\infty\,R^{\eins_{[0,1)},\chi_1^0}_r((\kappa\,\varrho)(s,\cdot))\,(x)\,\Big)\de  x\de\nu(s)\\
		\begin{aligned}
			&\leq\,\omega_r^g\, \|f\|_{L^\infty(\nu)}\int\limits_S\int\limits_{\R^d}|u(s,x)|\de m_s(x)\de\nu(s)\\
			&\phantom{\leq}\,+\|g\|_\infty\,\|f\|_{L^\infty(\nu)}\,\int\limits_S\delta_r^{\kappa(s,\cdot)}(m_s)\,\Big(\int\limits_{\R^d}|u(s,x)|^2\de m_s(x)\Big)^{\frac{1}{2}}\de\nu(s).
		\end{aligned}
	\end{multline*}
	The claim of (ii) now follows with the Cauchy-Schwarz inequality.
	
	We now turn to the proof of (iii). We fix $i\in\{1,\dots,d\}$. For $\alpha\in r\,\Z^d$ and $q\in(0,\infty)$ we set
	\begin{equation*}
		\mathscr T_q^{\alpha,i}:=\Big\{T\in\mathscr T_q \,\Big|\,T(\sigma_T^{-1}(i)-1)=\alpha\Big\}
	\end{equation*}
	with the notation of \cref{sec:kuhn}.
	In view of \cref{eqn:tirect} this condition means, that after hitting $\alpha$ the path of $T$ takes the direction of $\univ_i$ 
	until it reaches the next point on the lattice $\alpha+q\,\univ_i$.
	We note, that by definition of $ \mathscr T_q$ and \cref{eqn:suniv} 
	an element from $T\in \mathscr T_q$ is already uniquely defined if we are given its corresponding permutation $\sigma_T\in \mathscr S_d$ 
	and one of its vertexes $T(j)\in q\,\Z^d$ for an arbitrary index $j\in\{0,\dots,d\}$.
	So, the size of $\mathscr T_q^{0,i}$ calculates as
	\begin{align*}
		|\mathscr T_q^{0,i}|\,&=\,\sum_{j=0}^{d-1}\Big|\Big\{T\in\mathscr T_q \,\Big|\,T(j)=0\textnormal{ and }\,T(\sigma_T^{-1}(i)-1)=T(j)\Big\}\Big|\\
		&=\,\sum_{j=0}^{d-1}\Big|\Big\{T\in\mathscr T_q \,\Big|\,T(j)=0\textnormal{ and }\,\sigma_T(j+1)=i\Big\}\Big|\\
		&=(d-1)!\,d=d!.
	\end{align*}
	Hence, 
	\begin{equation*}
		\tilde\eta:=\sum_{T\in\mathscr T_1^{0,i}}\eins_{D_T}
	\end{equation*}
	defines an element from $\mathscr C$, because $|D_T|=1/(d!)$ for $T\in\mathscr T_1$ and moreover
	\begin{equation*}
		\sum_{\alpha\in \Z^d}\tilde\eta_1^\alpha\,=\,\sum_{\alpha\in \Z^d}\sum_{T\in\mathscr T_1^{\alpha,i}}\eins_{D_T}\,=\sum_{T\in\mathscr T_1}\eins_{D_T}\,=\,1.
	\end{equation*}  	
	We define another element from $\mathscr C$, 
	\begin{equation*}
		\tilde\varphi:\R^d\ni y\mapsto \int_0^1\eins_{[0,1)^d}(y-t\,\univ_i)\de t .
	\end{equation*}
	Let $s\in S$ and $i\in\{1,\dots,d\}$. 
	The main effort in the proof of the estimate of (iii) is done by a preliminary transformation of a relevant integral.
	To shorten the notation in the next lines we set 
	\begin{equation*}
		I_T\,:=\,\int\limits_{D_T}\partial_i g(x)\,\kappa(s,x)\,\varrho(s,x)\de x
	\end{equation*}
	for $T\in\mathscr T_r$.
	To obtain the equivalences as follows, we first use Fubini's theorem, then apply \partuniii{} before we use 
	the translation invariance of the Lebesgue measure and the fundamental theorem of calculus.
	\begin{align}\label{eqn:funda}
		&\int\limits_{\R^d}\partial_{x_i}\lambda(s,x)\,\partial_i g(x)\,\kappa(s,x)\de m_s(x)\nonumber\\
		&=\sum_{\alpha\in r\Z^d}r^{-d}\int\limits_{\R^d} u(s,y)\,\eins_{(-r,0]^d}(\alpha-y)\de y\int\limits_{\R^d}\partial_i\chi^\alpha_r(x)\,\partial_i g(x)\,(\kappa\,\varrho)(s,x)\de x\nonumber\\
		&=\,\sum_{\alpha\in r\Z^d}\sum_{T\in\mathscr T_r^{\alpha,i}}r^{-d-1}\int\limits_{\R^d} u(s,y)\,\big(\eins_{(-r,0]^d}(\alpha+r\,\univ_i-y)-\eins_{(-r,0]^d}(\alpha-y)\big)\de y
		\,I_T\nonumber\\
		&=\,\sum_{\alpha\in r\Z^d}\sum_{T\in\mathscr T_r^{\alpha,i}}r^{-d-1}\int\limits_{\R^d} \big(u(s,y+r\,\univ_i)-u(s,y)\big)\,\eins_{(-r,0]^d}(\alpha-y)\de y
		\,I_T\nonumber\\
		&=\,\sum_{\alpha\in r\Z^d}\sum_{T\in\mathscr T_r^{\alpha,i}}r^{-d-1}\int\limits_{\R^d} \int_{0}^r\frac{\partial u(s,x+t\,\univ_i)}{\partial x_i}\de t\,\eins_{(-r,0]^d}(\alpha-x)\de x
		\,I_T\nonumber\\
		&=\,\sum_{\alpha\in r\Z^d}\sum_{T\in\mathscr T_r^{\alpha,i}}r^{-d-1}\int\limits_{\R^d}\partial_{x_i}u(s,x)\int_{0}^r\eins_{(-r,0]^d}(\alpha-x+t\,\univ_i)\de t\de x
		\,I_T\nonumber\\
		&=\,\int\limits_{\R^d}\partial_{x_i}u(s,x)\,\sum_{\alpha\in r\Z^d}\tilde\varphi_r^\alpha(x)\,
		r^{-d}\,\int\limits_{\R^d}\tilde \eta_r^\alpha(y)\, \partial_i g(y)\,(\kappa\,\varrho)(s,y) \de y\de x.
	\end{align}
	In the second to last step we change the order in which we integrate w.r.t.~$\de t$ and $\de x_i$, then use the translation invariance of $\de x_i$
	before we change back.
	If we subtract the integral calculated in \cref{eqn:funda} from the term $\int_{\R^d}\partial_i u(s,x)$ $\partial_i g(x)$ $\kappa(s,x)\de m_s(x)$,
	then we can make use of \eqnAbsineq{} by choosing 
	$\kappa(s,\cdot)\,\varrho(s,\cdot)$ as $\rho$, the primal functions $\tilde\varphi$,  $\tilde \eta$, as $\varphi$, respectively $\eta$, and $\partial_i g$ as $h$.
	Summing up over $i=1,\dots,d$ and integrating w.r.t. $f\de\nu$ over $S$, we obtain the inequality
	\begin{multline*}
		\Big|\int\limits_Sf(s)\,E^{\kappa\varrho(s,\cdot)}\big(g,u(s,\cdot)-\lambda(s,\cdot)\big)\de\nu(s)\,\Big|\\
		\begin{aligned}
			&\leq\, \|f\|_{L^\infty(\nu)}\,\omega_{r}^{\nabla g}\,\sum_{i=1}^d\int\limits_S\int\limits_{\R^d}|\partial_{x_i}u(s,x)|\de m_s(x)\de\nu(s)\\
			&\phantom{\leq}\,+\|f\|_{L^\infty(\nu)}\,D_g\,\sum_{i=1}^d\int\limits_S\delta^{\kappa(s,\cdot)}_{r}(m_s)\,\Big(\int\limits_{\R^d}|\partial_{x_i} u(s,x)|^2\de m_s(x)\Big)^{\frac{1}{2}}\de\nu(s)\,.
		\end{aligned}
	\end{multline*} 
	Now, (iii) follows by applying the Cauchy-Schwarz inequality twice on each summand.
	
	We approach the missing proof of Item (iv). In a similar calculation as done in \cref{eqn:funda},
	we first use \partuniii{} and the shift invariance of the Lebesgue measure,
	before we estimate twice with the Cauchy-Schwarz inequality to obtain 
	\begin{align}\label{eqn:nablest}
		&\sum_{i=1}^d|\,\partial_{i}\lambda(s,\cdot)\,|^2
		\,=\,\sum_{i=1}^d\Big|\sum_{\alpha\in r\Z^d}r^{-d}\int\limits_{\R^d} u(s,y)\,\eins_{(-r,0]^d}(\alpha-y)\de y\,\partial_i\chi_r^\alpha\,\Big|^2\nonumber\\
		&=\,\frac{1}{r^2}\sum_{T\in\mathscr T_r}\eins_{D_T}\sum_{i=1}^d\Big|\,r^{-d}\int\limits_{\R^d} u(s,y)\,\big(\eins_{(-r,0]^d}(T(i)-y)-\eins_{(-r,0]^d}(T(i-1)-y)\big)\de y\,\Big|^2\nonumber\\
		&=\,\frac{1}{r^2}\sum_{i=1}^d\sum_{\alpha\in r\Z^d}\sum_{T\in\mathscr T_r^{\alpha,i}}\eins_{D_T}\,\Big|\,r^{-d}
		\int\limits_{\R^d} u(s,y)\,\big(\eins_{(-r,0]^d}(\alpha+r\,\univ_{i}-y)-\eins_{(-r,0]^d}(\alpha-y)\big)\de y\,\Big|^2\nonumber\\
		&=\,\frac{1}{r^2}\sum_{i=1}^d\sum_{\alpha\in r\Z^d}\sum_{T\in\mathscr T_r^{\alpha,i}}\eins_{D_T}\,\Big|\,r^{-d}
		\int\limits_{\R^d} \big(u(s,y+r\,\univ_i)-u(s,y)\big)\,\eins_{(-r,0]^d}(\alpha-y)\de y\,\Big|^2\nonumber\\
		&\leq\,\frac{1}{r^2}\sum_{i=1}^d\sum_{\alpha\in r\Z^d}\sum_{T\in\mathscr T_r^{\alpha,i}}\eins_{D_T}\,r^{-d}
		\int\limits_{\R^d} |u(s,y+r\,\univ_i)-u(s,y)|^2\,\eins_{(-r,0]^d}(\alpha-y)\de y\nonumber\\
		&\leq\,\sum_{i=1}^d\sum_{\alpha\in r\Z^d}\sum_{T\in\mathscr T_r^{\alpha,i}}r^{-d-1}\,\eins_{D_T}
		\int\limits_{\R^d}\int\limits_{0}^r \Big|\frac{\partial u(s,y+t\,\univ_i)}{\partial y_i}\Big|^2\de t\,\eins_{(-r,0]^d}(\alpha-y)\de y\nonumber\\
		&=\,\sum_{i=1}^d\sum_{\alpha\in r\Z^d}r^{-d}\,\tilde\eta^\alpha_r\int\limits_{\R^d}|\partial_{y_i}u(s,y)|^2\,\tilde\varphi^\alpha_r(y)\de y.
	\end{align}
	The fundamental theorem of calculus is used in the second to last step.
	In the last equality, for each $i=1,\dots,d$, we change the order in which we integrate w.r.t.~$\de t$ and $\de y_i$, then use the translation invariance of $\de y_i$
	before we change back.
	
	We first integrate the function of \cref{eqn:nablest} w.r.t.~$\kappa(s,\cdot)$ $\de m_s$ over $\R^d$ and then we integrate w.r.t.~$\de\nu$ over the variable $s\in S$.
	Looking at the integrated version of \cref{eqn:nablest}, we observe that
	the left hand side coincides with the left hand side of \cref{eqn:polygon3}.
	The right hand side yields the desired upper bound, because with Fubini's theorem we write
	\begin{align*}
		&\int\limits_{S}\int\limits_{\R^d}
		\sum_{\alpha\in r\Z^d}r^{-d}\,\tilde\eta^\alpha_r(x)\int\limits_{\R^d}\sum_{i=1}^d|\partial_{y_i}u(s,y)|^2\,\tilde\varphi^\alpha_r(y)\de y\,\kappa(s,x)\de m_s(x)\de\nu(s)\\
		&=\,\int\limits_{S}\int\limits_{\R^d}\sum_{i=1}^d|\partial_{i}u(s,\cdot)|^2\, I_r^{\tilde\varphi,\tilde\eta}(\,\kappa(s,\cdot)\,\,\varrho(s,\cdot)\,)\de y\de\nu(s)\\
		&\leq\,\int\limits_{S}C_r^{\kappa(s,\cdot)}(m_s)\int\limits_{\R^d}\sum_{i=1}^d|\partial_{i}u(s,\cdot)|^2\,\de m_s\de\nu(s)
	\end{align*}
	This concludes the proof.
\end{proof}

\section{Preliminaries on Mosco convergence and main results}\label{sec:lowse}
\subsection{Basic terminology and the theorem of Mosco-Kuwae-Shioya}\label{sec:mokushi}
\newcounter{c}
\setcounter{c}{\value{section}}
\newcommand{\mimii}{\cref{eqn:m1m2m}}
\newcommand{\kushi}{\cref{rem:extop}}
For the convenience of the reader  we give a self contained introduction to the most elementary concepts developed in \cite{umberto, kuwae}. 
The section comprises all the aspects from this theory which are relevant to this article.
The theorem of Mosco-Kuwae-Shioya defines a notion of convergence for spectral structures over varying Hilbert spaces, indexed by $N$, and
finds equivalent formulations in terms of semigroup ${(T_t^N)}_{t>0}$, resolvent ${(G_\alpha^N)}_{\alpha>0}$, or symmetric closed form $\E^N$.
The manifestation of the central theorem, as it is arranged in this text, contains a simplification for the condition of {(M1)}.
In both of the original papers \cite{umberto, kuwae} the validity of this modification is evident from their proofs, however has not been stated explicitly.
In its traditional formulation {(M1)} reads exactly as Property (a) of \mimii{}. It demands a verification for the sequential lower-semi-continuity
of ${(\E^N)}_N$ considering an abstract sequence ${(u_N)}_{N\in\N}$ and its weak limit.
Now, \amind{} says that we may restrict to the case where $u_N$ is in the image set defined by the action of $G_\alpha^N$ on a certain well-known class of pre-images for $N\in\N$
and some fixed value $\alpha>0$.
This observation is particularly useful in the context of Dirichlet forms with $\alpha\,G_\alpha ^N$ being sub-Markovian.
In the proof of \cref{thm:m1} we can benefit from it.

All abstract Hilbert spaces are assumed to be real and separable. A sequence of converging Hilbert spaces comprises linear maps
\[\Psi_N:\mathcal C\to H_N\] indexed by the parameter $N\in\overline\N$, where $\mathcal C$ is a dense linear subspace of a 
Hilbert space $(H_\infty,{\langle\cdot,\cdot\rangle}_\infty)$ and the image space $(H_N,{\langle\cdot,\cdot\rangle}_N)$ is Hilbert as well.
Apart from that, the asymptotic equations
\begin{align}\label{eqn:initop}
	\Psi_\infty\varphi\,&=\,\varphi,\nonumber\\
	\lim_{N\to\infty}\,{\langle\Psi_N\varphi,\Psi_N\varphi\rangle}_N\,&=\,{\langle\varphi,\varphi\rangle}_\infty
\end{align}
are assumed to hold for $\varphi\in\mathcal C$. For $N\in\overline\N$ the norm on $H_N$ is denoted by $\|\cdot\|_N$.
An element of 
\begin{equation}\label{eqn:construct}
	\mathcal H:=\prod_{N\in\overline \N}H_N
\end{equation}
is referred to as a \textbf {section} in this article.
The reasoning behind this terminology becomes clear in \kushi{} after the next lemma.
Moreover, we say that a section ${(u_N)}_{N\in\overline\N}$ is \textbf{strongly convergent} if
\begin{align*}
	\lim_{N\to\infty}{\langle u_N,\Psi_N\varphi\rangle}_N\,&=\,{\langle u_\infty,\varphi\rangle}_\infty\qquad\text{for }\varphi\in\mathcal C\\
	\textnormal{and}\qquad\lim_N{\|u_N\|}_N\,&=\,{\|u_\infty\|}_\infty.
\end{align*}
Building a dual notion, the section is called \textbf{weakly convergent} if
\begin{equation*}
	\lim_{N\to\infty}{\langle u_N,v_N\rangle}_N\,=\,{\langle u_\infty,v_\infty\rangle}_\infty
\end{equation*}
holds true for every strongly convergent section ${(v_N)}_{N\in\overline\N}$. What is more, if ${(N_k)}_{k\in\overline\N}\subset\overline\N$ is strictly increasing in $k$, then
${(u_{N_k})}_{k\in\overline\N}$ is referred to as a \textbf{subsection} of ${(u_N)}_{N\in\overline\N}$.
The terminology of strong and weak convergence naturally applies to subsections as well. If ${(u_N)}_{N\in\overline\N}$ is a section which has a weakly (or strongly) convergent subsection, then $u_\infty$
is called a weak (respectively strong) accumulation point of ${(u_N)}_{N\in\N}$.
The next lemma cites some results from \cite{kuwae,koles,kolesni} to better understand the newly introduced terminology.
The proof given here takes the same route as the one in \cite{koles}.
\begin{lem}\label{lem:ks}
	Let $H_N$, $N\in\N$, be a sequence of converging Hilbert spaces with asymptotic space $H_\infty$.
	\begin{thmlist}[wide, labelindent=0pt]
		\item\label{lem:ksi} For every $u\in H_\infty$ there is a strongly convergent section which has $u$ as its asymptotic element.
		\item\label{lem:ksii} A section ${(u_N)}_{N\in\overline\N}$ is strongly convergent if and only if 
		\begin{equation*}
			\lim_{N\to\infty}{\langle u_N,v_N\rangle}_N\,=\,{\langle u_\infty,v_\infty\rangle}_\infty
		\end{equation*}
		holds true for every weakly convergent section ${(v_N)}_{N\in\overline\N}$.
		\item\label{lem:ksiii} The norm is weakly lower semi-continuous. By this we mean 
		\begin{equation*}
			{\|u_\infty\|}_\infty\,\leq\,\liminf_{N\in\N}\,{\|u_N\|}_{N} 
		\end{equation*}
		for every weakly convergent section ${(u_N)}_{N\in\overline\N}$. Moreover, the right hand side of this inequality takes a finite value.
		\item\label{lem:ksiv} If $u_N\in\{u\in H_N\,|\,\|u\|_N\leq 1\}$ for $N\in\N$, then there exists a weak accumulation point of ${(u_N)}_{N\in\N}$.
	\end{thmlist}
\end{lem}
\newcommand{\ksi}{\cref{lem:ksi}}
\newcommand{\ksii}{\cref{lem:ksii}}
\newcommand{\lose}{\cref{lem:ksiii}}
\newcommand{\wesecomp}{\cref{lem:ksiv}}
\begin{proof}
	Let $\varphi_1,\varphi_2,\dots$ be elements from $\mathcal C$ which form an orthonormal basis for $H_\infty$.
	All the statements are clear if $\Psi_N$ is isometric for $N\in\overline \N$, since in this case there is a one-to-one identification 
	\begin{equation*}
		\mathcal H\ni {(u_N)}_{N\in\overline\N}\,\mapsto\,[{(\tilde u_N)}_{N\in\N}\,,\,\tilde u_\infty]
		\in {(l^2)}^\N\times l^2
	\end{equation*}
	if we set
	\begin{equation*}
		\tilde u_N\,:=\,\big(\,{\langle u_N,\Psi_N\varphi_1\rangle}_{N}\,,\,{\langle u_N,\Psi_N\varphi_2\rangle}_{N}\,,\,\cdots\big)\in l^2
	\end{equation*}
	for $N\in\overline\N$.
	It correctly explains the notion of strongly and weakly convergent sections through the usual strong and weak topology on $l^2$. This means
	${(u_N)}_{N\in\overline\N}$ is a strongly (respectively weakly) convergent section if and only if
	$\lim_{N\to\infty}\tilde u_N=\tilde u_\infty$ holds in the strong (respectively weak) topology of $l^2$.
	Moreover, $\|\tilde u_N\|_{l^2}=\|u_N\|_N$ for $N\in\overline \N$. So, if $\Psi_N$ is isometric, then the claim of the lemma follows from the analogous 
	facts of (i) to (iv) for $l^2$.
	
	We now construct isometric isomorphisms $\hat\Psi_N:H_\infty\to H_N$ 
	for $N\in\N$ which yield the same notion of strong and weak convergence for elements of $\mathcal H$ as the given ones.
	For $N,m\in\N$ define
	\begin{equation*}
		A^{N,m}\,:=\,{\big[a^{N,m}_{i,j}\big]}_{i,j=1}^{m}\,:=\,{\big[{\langle\Psi_N\varphi_i,\Psi_N\varphi_j\rangle}_N\big]}_{i,j=1}^{m}.
	\end{equation*}
	For fixed $m\in\N$ we have \[\lim_N A^{N,m}=\textnormal{id}\in\R^{m\times m}.\]
	Hence, for $m\in\N$ there is $N_m\in\N$ such that for $N\geq N_m$ the following is true: There exists $B^{N,m}:={[b^{N,m}_{i,j}]}_{i,j=1}^m\in\R^{m\times m}$ with
	\begin{align}
		\| \,B^{N,m}-\textnormal{id}\,\|_{\textnormal{op},\infty}\,&\leq\,\frac{1}{m},\label{eqn:simeq2}\\
		{(B^{N,m})}^{\textnormal{T}}\,A^{N,m}\,B^{N,m}\,&=\,\textnormal{id}\in\R^{m\times m},\label{eqn:aba}\\
		\textnormal{and}\qquad\|\Psi_N\varphi_i\|_N\,&\leq\,2,\quad i=1,\dots,m.\label{eqn:simeq1}
	\end{align}
	Indeed, $B^{N,m}$ can be set by $B^{N,m}:=(A^{N,m})^{-\frac{1}{2}}$ for sufficiently large $N$.
	In the line above $\|\cdot\|_{\textnormal{op},\infty}$ denotes the operator norm on $\R^{m\times m}$ w.r.t.~the supremum norm on $\R^m$.
	Now, for fixed $N\in\N$ we choose $m_N\in\N$ as the maximal $m\in\N$ for which $N_m\leq N$ and define
	\begin{equation*}
		\hat\Psi_N\varphi_j\,:=\,\sum_{i=1}^{m_N}b_{i,j}^{N,m_N}\,\Psi_N\varphi_i
	\end{equation*}
	for $j=1,\dots, m_N$. 
	\begin{align*}
		{\big[{\langle \hat\Psi_N\varphi_i,\hat\Psi_N\varphi_j\rangle}_N\big]}_{i,j=1}^{m_N}\,=\,\textnormal{id}\in\R^{m_N\times m_N}
	\end{align*}
	holds true due to \cref{eqn:aba}.
	So, $\hat\Psi_N$ can be extended to an isometric isomorphism $H_\infty\to H_N$ which we again denote by $\hat\Psi_N$. We further define $\hat\Psi_\infty:=\textnormal{id}_{H_\infty}$.
	
	To see that ${(\hat\Psi_N)}_{N\in\overline\N}$ indeed yield the same notion of strong convergence for elements of $\mathcal H$ as the given one, it suffices to check the asymptotic equality
	\begin{equation*}
		\lim_{N\to\infty}\Big|\,\langle u_N,\hat\Psi_N\varphi_j\rangle_N-\langle u_N,\Psi_N\varphi_j\rangle_N\,\Big|\,=\,0
	\end{equation*}
	for given $j\in\N$ and $u_N\in \{u\in H_N$ $|$ $\|u\|_N\leq 1\}$, $N\in\N$.
	For given $j\in\N$ and $N\in\N$ such that $m_N\geq j$ we calculate 
	\begin{equation*}
		|\,{\langle u_N,\hat\Psi_N\varphi_j-\Psi_N\varphi_j\rangle}_N\,|\,\leq\,\| \,B^{N,m_N}-\textnormal{id}\,\|_{\textnormal{op},\infty}\,
		\sup_{1\leq i\leq m_N}|\langle u_N,\Psi_N\varphi_i\rangle_N|\,\leq\,\frac{2}{m_N}
	\end{equation*}
	using \cref{eqn:simeq1,eqn:simeq2}. Since ${(\hat\Psi_N)}_{N\in\overline\N}$ induces the same notion of strong convergence for elements of $\mathcal H$ 
	as does ${(\Psi_N)}_{N\in\overline\N}$, the analogue statement concerning weak convergence is also true via duality.  This concludes the proof.
\end{proof}
We state some observations regarding \cref{lem:ks}.
\begin{rem}\label{rem:kurem}
	\begin{thmlist}[wide,labelindent=0pt]
		\item\label{rem:extop}Let ${(u_N)}_{N\in\overline\N}\in\mathcal H$. The map $N\mapsto u_N$ can be regarded as a section from $\overline\N$ into
		\begin{equation*}
			\mathcal KS \,:=\,\Big\{[N,u]\,\Big|\,N\in\overline\N,u\in H_N\Big\},
		\end{equation*}
		the disjoint union of the Hilbert spaces $H_N$, $N\in\overline\N$. We use the term `section' here in analogy to its meaning in the theory of vector bundles,
		where a section denotes a right inverse of the projection onto the base space.
		There are topologies $\tau_\textnormal{w}$ and $\tau_\textnormal{s}$ on $\mathcal KS$ such that for a section ${(u_N)}_{N\in\overline\N}\in\mathcal H$ 
		the strong (or weak) convergence as defined above is equivalent to $\lim_{N\to\infty}u_N=u_\infty$ w.r.t.~$\tau_\textnormal{s}$ (respectively $\tau_\textnormal{w}$).
		Indeed, these can be easily written down as initial topologies. Concerning $\tau_\textnormal{s}$, it is the initial topology generated by the family of maps
		\begin{align*}
			\Big\{\,\mathcal  KS\ni[N,u]\,&\mapsto\, N\in\overline\N,\\
			\mathcal  KS\ni[N,u]\,&\mapsto\, \|u\|_N\in\R,\\
			\mathcal KS\ni[N,u]\,&\mapsto\,{\langle\Psi_N\varphi,u\rangle}_{N}\in \R\,\big|\quad\varphi\in\mathcal C\,\Big\}.
		\end{align*}
		The open sets in $\overline\N$ are given as
		\begin{align*}
			\{\,U\,|\,U\subset \N\}\cup\Big\{\,U\cup\{\infty\}\,\Big|\,&U\subset\N\text{ and there exists }m\in\N
			\\&\text{such that }l\in U\text{ is true for }l\in\N,\,l\geq m\,\Big\}.
		\end{align*}
		Hence, for a convergent sequence ${(N_k)}_{k\in\N}$ in $\overline \N$ there are only two possibilities.
		Either there exists $m\in\N$ and $N^*\in\N$ such that $N_k=N^*$ for $k\in\N$, $k\geq m$, 
		or for arbitrarily chosen $M\in\N$ there exists $m\in\N$ such that $N_k\geq M$ holds for $k\in\N$, $k\geq m$.
		
		Then, concerning $\tau_\textnormal{w}$, it is the initial topology generated by the family of maps
		\begin{align*}
			\Big\{\,\mathcal KS\ni [N,u]\,&\mapsto\,N\in\overline\N,\\
			\mathcal KS\ni [N,u]\,&\mapsto\,{\langle u,v_N\rangle}_N\,\Big|\,{(v_N)}_{N\in\overline \N}\in\mathcal H
			\textnormal{ is a strongly convergent section}\,\Big\}.
		\end{align*}
		\item \label{rem:extopii}Due to \cref{rem:extop} a section ${(u_N)}_{N\in\overline\N}\in\mathcal H$ converges strongly (or weakly),
		if for every subsection there is a (sub-)subsection which does so. 
		\item\label{rem:wecrit}Let $V\subset H_\infty$ be a dense linear subspace, $u_\infty\in H_\infty$ and $u_N\in\{H_N$ $|$ $\|u_N\|_N\leq 1\}$ for $N\in\N$. 
		As a consequence of \cref{rem:extopii} and \wesecomp{} we obtain a sufficient criterion for 
		${(u_N)}_{N\in\overline\N}$ to form a weakly convergent section: For every $v_\infty\in V$
		there is a strongly convergent section with asymptotic element $v_\infty$, say ${(v_N)}_{N\in\overline\N}$, such that
		\begin{equation}\label{eqn:inwecrit}
			\lim_{N\to\infty}\,{\langle u_N,v_N\rangle}_N\,=\,{\langle u_\infty, v_\infty\rangle}_\infty.
		\end{equation}
		Indeed, \cref{eqn:inwecrit} allows to identify all weak accumulation points of ${(u_N)}_{N\in\overline\N}$ with the element $u_\infty$.
		\item \label{rem:cocrit}The proof of \cref{lem:ks} motivates to ask: Would the same notion of weakly and strongly convergent sections have emerged, had  
		the construction been initiated with a different choice $\tilde\Psi_N:\mathcal D\to\ H_N$ instead of the original map $\Psi_N$ for $N\in\overline\N$ ?
		Of course, the question only makes sense if ${(\tilde\Psi_N)}_{n\in\overline\N}$ meets the analogue of \cref{eqn:initop} 
		w.r.t.~the dense linear subspace $\mathcal D\subset H_\infty$. The answer is affirmative if and only if
		\begin{equation}\label{eqn:confirm}
			\lim_{N\to\infty}{\langle\tilde\Psi_N\varphi,\Psi_N\eta\rangle}_N\,=\,{\langle\varphi,\eta\rangle}_\infty
		\end{equation}
		for $\varphi\in\mathcal D$ and $\eta\in\mathcal C$. The necessity of \cref{eqn:confirm} is clear indeed. On the other hand, \cref{eqn:confirm} 
		implies the strong convergence of the section ${(\tilde\Psi_N \varphi)}_{N\in\overline\N}$ w.r.t~the notion induced by ${(\Psi_N)}_{N\in\overline\N}$ and vice versa.
		Hence, concerning the notion of strong convergence, the answer to the question is affirmative under \cref{eqn:confirm}.
		The equivalence of weak convergence follows from the equivalence of strong convergence.
		
		The proof for the equivalence between the strong convergence is easily done under \cref{eqn:confirm}. Indeed, for $\eta\in\mathcal C$ we have the weak convergence of ${(\Psi_N\eta)}_{N\in\N}$ to $\eta$ w.r.t.~the topology induced by ${(\tilde\Psi_N)}_{N\in\overline\N}$ by way of \cref{rem:wecrit}. 
		Then, for any strongly convergent section ${(v_N)}_{N\in\overline\N}$ w.r.t.~the topology induced by ${(\tilde\Psi_N)}_{N\in\overline\N}$, we have the convergence of  
		${({\langle v_N,\Psi_N\eta\rangle}_N)}_{N\in\N}$ to $\langle v_\infty,\eta\rangle_\infty$ for $\eta\in\mathcal C$ 
		and the convergence of ${({\|v_N\|}_N)}_{N\in\N}$ to $\|v_\infty\|_\infty$.
		Hence, the strong convergence of ${(v_N)}_{N\in\N}$ to $v_\infty$ also holds true w.r.t.~the topology induced by ${(\Psi_N)}_{N\in\overline\N}$.
		
		\item As we learn in the proof of \cref{lem:ks} there are isometric isomorphisms $\hat\Psi_N:H_\infty\to H_N$ for $N\in\N$ such that \cref{eqn:confirm} holds
		with $\mathcal D:= H_\infty$.
	\end{thmlist}
\end{rem}
\newcommand{\extopii}{\cref{rem:extopii}}
\newcommand{\wecrit}{\cref{rem:wecrit}}
\newcommand{\cocrit}{\cref{rem:cocrit}}
We are now preparing to state the theorem of Mosco-Kuwae-Shioya. Again $H_N$, $N\in\N$, is a sequence of converging Hilbert spaces with asymptotic space $H_\infty$
There is a natural way to introduce a notion of convergence for
an element ${(L_N)}_{N\in\overline\N}$ of \[\mathcal L(\mathcal H):=\prod_{N\in\overline \N}L(H_N).\]
$L(H_N)$ denotes the Banach space of bounded linear operators on $H_N$ with operator norm ${\|\cdot\|_{L(H_N)}}$ for $N\in\overline\N$.
Again it makes sense to refer to the elements of $\mathcal L(\mathcal H)$ as sections.
The section ${(L_N)}_{N\in\overline\N}$ is called strongly convergent if ${(L_Nu_N)}_{N\in\overline\N}\in\mathcal H$ converges strongly
for any strongly convergent section ${(u_N)}_{N\in\overline\N}$. 
\begin{rem}\label{rem:adjoint}
	Let ${(L_N)}_{N\in\overline\N}$, ${(L^*_N)}_{N\in\overline\N}$ be elements of the set
	\begin{equation*}
		\prod_{N\in\overline \N}\{L\in L(H_N)\,|\,\|L\|_{L(H_N)}\leq 1\},
	\end{equation*}
	where $L^*_N$ denotes the adjoint of $L_N$ for $N\in\overline\N$.
	Due to \ksii{} the strong convergence of ${(L_N)}_{N\in\overline\N}\in\mathcal L(\mathcal H)$ is equivalent to the following:
	${(L_N^*u_N)}_{N\in\overline\N}\in\mathcal H$ converges weakly for any weakly convergent section ${(u_N)}_{N\in\overline\N}$.
	By \wecrit{} however, the latter condition is in turn equivalent to the following property:
	For every $\varphi\in\mathcal C$ the section ${(L_N\Psi_N\varphi)}_{N\in\overline\N}$ is strongly convergent.
\end{rem}
The theorem of Mosco-Kuwae-Shioya throws some light on a family $\{{(G_\alpha^N)}_{N\in\overline \N}$ $|\alpha>0\}$ $\subset \mathcal L(\mathcal H)$ of strongly convergent sections,
where ${(G_\alpha^N)}_{\alpha>0}$ is assumed to form a strongly continuous contraction resolvent of symmetric operators on $H_N$ for fixed $N\in\overline\N$.
Its associated generator
\begin{equation*}
	\Delta_N\,:=\,\textnormal{id}-{(G_1^N)}^{-1},\qquad\dom(\Delta_N)\,:=\,\textnormal{Im}(G_1^N),
\end{equation*}
is densely defined and induces a non-negative, symmetric bilinear form 
\begin{equation*}
	\E^N(u,v)\,:=\,\langle u,-\Delta_N v\rangle_N,\qquad u,v\in\dom(\Delta_N).
\end{equation*}
This form is closable on $H_N$. Its closure is denoted by $(\E^N,\dom(\E^N))$ and satisfies
\begin{equation*}
	\E^N(G_\alpha^Nu,v)+\alpha\,\langle G_\alpha^Nu,v\rangle_N\,=\,\langle u,v\rangle_N
\end{equation*}
for $u\in H_N$, $v\in\dom(\E^N)$ and $\alpha>0$ due to the identity
\begin{equation*}
	G_\alpha^N\,=\,{(\alpha-\Delta_N)}^{-1}.
\end{equation*}
Since the spectrum of $\Delta_N$ is contained in $(-\infty,0]$,
the functional calculus (see e.g.~\cite[Chap. VII]{werner}) evaluating the exponential function at $t\,\Delta_N$, $t>0$, yields a strongly continuous contraction semigroup of symmetric operators 
\begin{equation*}
	\{\,T_t^N\,:=\,\exp(t\,\Delta_N)\,|\, t>0\,\}
\end{equation*}
on $H_N$. For $\alpha>0$ we write $\E_\alpha^N(u,v)=\E^N(u,v)+\alpha\langle u,v\rangle_N$ for $u,v\in\dom(\E^N)$.
Then $\E^N_\alpha$ defines a scalar product which makes $(\dom(\E^N),\E^N_\alpha)$ a Hilbert space.
The induced norm ${(\E_\alpha^N)}^{1/2}$ is equivalent to ${(\E_1^N)}^{1/2}$ for $\alpha>0$.
We shorten the notation a bit. If ${(u_N)}_{N\in\overline\N}\in\mathcal H$ is strongly (or weakly) convergent, we write `$u_N\overset{\textnormal{s.}}{\underset{N}{\longrightarrow}} u_\infty$'
(respectively `$u_N\overset{\textnormal{w.}}{\underset{N}{\rightharpoonup}} u_\infty$').
Analogously, we write `$L_N\overset{\textnormal{s.}}{\underset{N}{\longrightarrow}}L_\infty$' if ${(L_N)}_{N\in\overline\N}\in\mathcal L(\mathcal H)$ is strongly convergent.

\begin{thm}\label{thm:mokush}
	The following are equivalent.
	\begin{thmlist}
		\item \label{thm:mokushi}$G_\alpha^N\overset{\textnormal{s.}}{\underset{N}{\longrightarrow}} G_\alpha^\infty$ for $\alpha>0$.
		\item \label{thm:mokushii}
		(a) Let ${(u_N)}_{N\in\overline\N}\in \mathcal H$ and $\alpha>0$. Then, $u_N\overset{\textnormal{s.}}{\underset{N}{\longrightarrow}} u_\infty$ implies
		\begin{equation*}
			\lim_{N\to\infty} \,\E^N_\alpha(G_\alpha^Nu_N,G_\alpha^Nu_N)\,=\,\E^\infty_\alpha(G^\infty_\alpha u_\infty, G^\infty_\alpha u_\infty).
		\end{equation*}
		(b) For every $u\in\dom(\E^\infty)$ there is $u_N\in\dom(\E^N)$ for $N\in\N$ such that  
		$u_N\overset{\textnormal{s.}}{\underset{N}{\longrightarrow}} u$ and
		\begin{equation*}
			\lim_{N\to\infty}\E^N(u_N,u_N)\,=\,\E^\infty(u,u).
		\end{equation*}
		\item \label{eqn:m1m2m}
		(a) Let ${(u_N)}_{N\in\overline\N}\in\mathcal H$. Then, $u_N\overset{\textnormal{w.}}{\underset{N}{\rightharpoonup}} u_\infty$ implies
		\begin{equation*}
			\E^\infty(u_\infty,u_\infty)\,\leq\,\liminf_{N\to\infty}\E^N(u_N,u_N).
		\end{equation*}
		The inequality has to be read in the sense, that in 
		case $\#N$ with $u_N\in\dom(\E^N)$ is infinite and accounts for a finite right hand side, then $u_\infty\in\dom(\E^\infty)$ and the stated inequality holds true.\\
		(b) There is a dense linear subspace $V\subset (\dom(\E^N),\E^N_1)$ such that
		for every $u\in V$ there exists $u_N\in\dom(\E^N)$ for $N\in\N$ with  $u_N\overset{\textnormal{s.}}{\underset{N}{\longrightarrow}} u$ and
		\begin{equation*}
			\lim_{N\to\infty}\E^N(u_N,u_N)\,=\,\E^\infty(u,u).
		\end{equation*}
		\item \label{eqn:am1nd}
		(a) There exists $\alpha>0$ such that for every $\varphi\in\mathcal C$ and every weak accumulation point $u$ of ${(G_\alpha^N\Psi_N\varphi)}_{N\in\N}$ it holds
		\begin{equation*}
			\E^\infty(u,u)\,\leq\,\liminf_{k\to\infty}\E^{N_k}(G_\alpha^{N_k}\Psi_{N_k}\varphi,G_\alpha^{N_k}\Psi_{N_k}\varphi)
		\end{equation*}
		in case $G_\alpha^{N_k}\Psi_{N_k}\varphi\overset{\textnormal{w.}}{\underset{k}{\rightharpoonup}} u$ is a corresponding weakly convergent subsection.\\
		(b) Property (b) of \cref{eqn:m1m2m} holds true.
		\item\label{thm:mokushv} $\Delta_N^p\,T_t^N\overset{\textnormal{s.}}{\underset{N}{\longrightarrow}}\Delta_\infty^p\,T_t^\infty$ for $t>0$ and $p\geq0$.
	\end{thmlist}
\end{thm}
\begin{proof}
	Assume (i). Let $\alpha>0$. Property (a) of (ii) is a direct consequence. Then, the linear maps
	\begin{equation}\label{eqn:hei}
		G_\alpha^\infty(\mathcal C)\ni u\,\mapsto\,G_\alpha^N\Psi_N(\alpha-\Delta_\infty)u, \quad N\in\overline\N,
	\end{equation}
	make the sequence $(\dom(\E^N),\E^N_\alpha)_{N\in\N}$ a convergent sequence of Hilbert spaces on their own right, with asymptotic space 
	$(\dom(\E^\infty),\E^\infty_\alpha)$, since they satisfy the asymptotic equations
	\begin{align}\label{eqn:asympeq}
		G_\alpha^N\Psi_\infty(\alpha-\Delta_\infty)u\,&=\,u,\nonumber\\
		\lim_{N\to\infty}\,\E^N_\alpha(\,G_\alpha^N\Psi_N(\alpha-\Delta_\infty)u\,,\,G_\alpha^N\Psi_N(\alpha-\Delta_\infty)u\,)\,&=\,\E^\infty_\alpha(u,u)
	\end{align}
	for $u\in G_\alpha^\infty(\mathcal C)$.
	We hint at the fact that $G_\alpha^\infty(\mathcal C)$ 
	is a dense linear subspace of $(\dom(\E^\infty),\E^\infty_\alpha)$, because $v\in\dom(\E^\infty)$ such that
	\begin{equation*}
		0\,=\,\E_\alpha^\infty(v,G_\alpha^\infty\varphi)\,=\,\langle v,\varphi\rangle_\infty\qquad\textnormal{for all }\varphi\in\mathcal C
	\end{equation*}
	implies $v=0$. For short we set
	\begin{equation}\label{eqn:heii}
		\mathcal H^\E\,:=\,\prod_{N\in\overline\N}(\dom(\E^N),\E^N_\alpha).
	\end{equation}
	In the following lines, we exploit the interplay between the asymptotic of the identity
	\begin{equation}\label{eqn:lookouti}
		\E^N_\alpha(u_N,G_\alpha^Nv_N)\,=\,{\langle u_N,v_N\rangle}_N,\quad N\in\N,
	\end{equation}
	for $N\to\infty$ and the identity 
	\begin{equation}\label{eqn:lookoutii}
		\E^\infty_\alpha(u_\infty,G_\alpha^\infty v_\infty)\,=\,{\langle u_\infty,v_\infty\rangle}_\infty
	\end{equation}
	for suitable choices of ${(u_N)}_{N\in\overline \N}\in\mathcal H^\E$ and ${(v_N)}_{N\in\overline \N}\in\mathcal H$.
	First we set $v_N:=\Psi_N(\alpha-\Delta_\infty)w$ for $N\in\overline\N$ and $w\in G_\alpha^\infty(\mathcal C)$ and deduce via \wecrit{} that 
	\begin{equation}\label{eqn:firstint}
		u_N\overset{\textnormal{w.}}{\underset{N}{\rightharpoonup}} u_\infty\textnormal{ in the sense of }\mathcal H^\E\;\textnormal{  implies  }\;
		u_N\overset{\textnormal{w.}}{\underset{N}{\rightharpoonup}} u_\infty\textnormal{ in the sense of }\mathcal H
	\end{equation}
	for any ${(u_N)}_{N\in\overline \N}\in\mathcal H^\E$.
	Consequently, 
	\begin{equation*}
		v_N\overset{\textnormal{w.}}{\underset{N}{\rightharpoonup}} v_\infty\textnormal{ in the sense of }\mathcal H\;\textnormal{  implies  }\; 
		G_\alpha^Nv_N\overset{\textnormal{w.}}{\underset{N}{\rightharpoonup}} G_\alpha^\infty v_\infty\textnormal{ in the sense of }\mathcal H^\E
	\end{equation*}
	for any ${(v_N)}_{N\in\overline \N}\in\mathcal H$, where
	we make use of \cref{rem:adjoint}, \wesecomp{} and \cref{eqn:firstint}.
	Then, again looking at \cref{eqn:lookouti,eqn:lookoutii}, we deduce from \ksii{} that 
	\begin{equation*}
		u_N\overset{\textnormal{s.}}{\underset{N}{\longrightarrow}} u_\infty\textnormal{ in the sense of }\mathcal H^\E\;\textnormal{  implies  }\;
		u_N\overset{\textnormal{s.}}{\underset{N}{\longrightarrow}} u_\infty\textnormal{ in the sense of }\mathcal H
	\end{equation*} 
	for any ${(u_N)}_{N\in\overline \N}\in\mathcal H^\E$.
	So, the desired Property (b) follows from \ksi and the implication of \cref{thm:mokushii} by \cref{thm:mokushi} is shown.
	
	Assume (ii) now and let $\alpha>0$. Again defining linear maps as in \cref{eqn:hei} 
	we perceive $(\dom(\E^N),\E^N_\alpha)_{N\in\N}$ as a convergent sequence of Hilbert spaces with limiting space $(\dom(\E^\infty),\E^\infty_\alpha)$,
	since Property (a) of (ii) ensures the validity of the asymptotic equation \cref{eqn:asympeq}.
	In the same way as we did before, we argue by comparing the asymptotic of \cref{eqn:lookouti} with \cref{eqn:lookoutii} that
	\begin{equation}\label{eqn:secint}
		u_N\overset{\textnormal{w.}}{\underset{N}{\rightharpoonup}} u_\infty\textnormal{ in the sense of }\mathcal H^\E\;\textnormal{  implies  }\;
		u_N\overset{\textnormal{w.}}{\underset{N}{\rightharpoonup}} u_\infty\textnormal{ in the sense of }\mathcal H
	\end{equation}
	for any ${(u_N)}_{N\in\overline \N}\in\mathcal H^\E$. We now consider an arbitrary element ${(u_N)}_{N\in\overline \N}\in\mathcal H$
	and increasing positive integers $N_k\in\N$, $k\in\N$.
	Due to \cref{eqn:secint}, \lose{} and \wesecomp{} we have $u_\infty\in\dom(\E^\infty)$ with
	\begin{equation}\label{eqn:lima}
		\E_\alpha^\infty(u_\infty,u_\infty)\,\leq\,\liminf_{k\to\infty}\E_\alpha^{N_k}(u_{N_k},u_{N_k})
	\end{equation}
	whenever $u_{N_k}\in\dom(\E^{N_k})$ for $k\in\N$ with $\sup_k\E_\alpha^{N_k}(u_{N_k},u_{N_k})$ $<\infty$ 
	and $u_{N_k}\overset{\textnormal{w.}}{\underset{k}{\rightharpoonup}} u_\infty$ in the sense of $\prod_k H_{N_k}$.
	Property (a) of \cref{eqn:m1m2m} now follows considering arbitrarily small $\alpha>0$ in \cref{eqn:lima}.
	
	The implication of (iii) by (iv) is clear.
	We assume (iv) now and fix $\alpha>0$ for which Property (a) holds true. Let $\varphi\in\mathcal C$.
	We can choose $N_k\in\overline \N$ for $k\in\overline\N$, strictly increasing in $k$, such that both,
	\begin{equation*}
		\limsup_{N\to\infty}\,\E^N_\alpha(G_\alpha^N\Psi_N\varphi,G_\alpha^N\Psi_N\varphi)\,=\,\lim_k\,\E^{N_k}_\alpha(G_\alpha^{N_k}\Psi_{N_k}\varphi,G_\alpha^{N_k}\Psi_{N_k}\varphi)
	\end{equation*}
	and there is a weak accumulation point $u$ with $G_\alpha^{N_k}\Psi_{N_k}\varphi\overset{\textnormal{w.}}{\underset{k}{\rightharpoonup}} u$. 
	In case 
	\begin{equation*}
		\limsup_{N\to\infty}\,\E^N_\alpha(G_\alpha^N\Psi_N\varphi,G_\alpha^N\Psi_N\varphi)\,>0
	\end{equation*}
	it holds
	\begin{multline}\label{eqn:limsup}
		\limsup_{N\to\infty}\,\E^N_\alpha(G_\alpha^N\Psi_N\varphi,G_\alpha^N\Psi_N\varphi)^{1/2}\,=
		\,\lim_k\,\frac{\langle G_\alpha^{N_k}\Psi_{N_k}\varphi,\Psi_{N_k}\varphi\rangle_{N_k}}{\E^{N_k}_\alpha(G_\alpha^{N_k}\Psi_{N_k}\varphi,G_\alpha^{N_k}\Psi_{N_k}\varphi)^{1/2}}\\
		\leq\,\frac{\langle u,\varphi\rangle_\infty}{\E^\infty_\alpha(u,u)^{1/2}}
		\,=\,\frac{\E_\alpha^\infty( u,G_\alpha^\infty\varphi)}{\E^\infty_\alpha(u,u)^{1/2}}
		\,\leq\,\E_\alpha^\infty(G_\alpha^\infty\varphi,G_\alpha^\infty\varphi)^{1/2}.
	\end{multline}
	Otherwise, the analogue of \cref{eqn:limsup} is automatically fulfilled.
	Property (b) of (iv) allows us to define linear maps $\tilde\Psi_N:\dom(\E^\infty)\supset V\to\dom(\E^N)$ for $N\in\overline\N$ such that both,
	\[\tilde\Psi_N u\overset{\textnormal{s.}}{\underset{N}{\longrightarrow}}u=:\tilde\Psi_\infty u\] in the sense of $\mathcal H$ and
	\begin{equation*}
		\lim_N\E_\alpha^N(\tilde\Psi_N u,\tilde\Psi_Nu)\,=\,\E_\alpha^\infty(u,u)
	\end{equation*}
	for $u\in V$. In this way we can understand $(\dom(\E^N),\E^N_\alpha)$, $N\in\N$, as a sequence of converging Hilbert spaces with asymptotic space $(\dom(\E^\infty),\E^\infty_\alpha)$.
	In the emerging notion of convergence for elements of $\mathcal H^\E$ (defined as in \cref{eqn:heii})
	it holds $G_\alpha^N\Psi_N\varphi\overset{\textnormal{w.}}{\underset{N}{\rightharpoonup}}G_\alpha^\infty\varphi$ 
	for $\varphi\in\mathcal C$, because 
	\begin{equation}\label{eqn:cocrit}
		\lim_{N\to\infty}\,\E^N_\alpha(G_\alpha^N\Psi_N\varphi,\tilde\Psi_N u)\,=\,\lim_{N\to\infty}{\langle\Psi_N\varphi,\tilde\Psi_Nu\rangle}_N
		\,=\,{\langle \varphi,u\rangle}_\infty\, =\,\E^\infty_\alpha(G_\alpha^\infty\varphi, u)
	\end{equation}
	for $u\in V$.
	Hence, by \lose{} and \cref{eqn:limsup} we have
	\begin{equation*}
		\lim_N\E^N_\alpha(G^N_\alpha\Psi_N\varphi,G^N_\alpha\Psi_N\varphi)\,=\,\E^\infty_\alpha(G^\infty_\alpha \varphi,G^\infty_\alpha \varphi).
	\end{equation*}
	for $\varphi\in\mathcal C$.
	We now deduce that the family of maps
	\begin{equation*}
		\tilde\Psi_N': G_\alpha^\infty(\mathcal C)\ni u\,\mapsto\,G_\alpha^N\Psi_N(\alpha-\Delta_N)u,\quad N\in\N,
	\end{equation*}
	as have already been regarded in \cref{eqn:hei}, fulfil the asymptotic equations \cref{eqn:asympeq}.
	Moreover, ${(\tilde\Psi_N')}_{N\in\overline\N}$ generates the same notion of convergence for elements of $\mathcal H^\E$ 
	as do ${(\tilde\Psi_N)}_{N\in\overline\N}$, due to \cref{eqn:cocrit} and \cocrit{}.
	Let ${(v_N)}_{N\in\overline\N}\in\mathcal H$ be weakly convergent. First, we argue that $G_\alpha^Nv_N\overset{\textnormal{w.}}{\underset{N}{\rightharpoonup}}G_\alpha^\infty v_\infty$
	in the sense of $\mathcal H^\E$, because 
	\begin{equation*}
		\lim_{N\to\infty}\,\E_\alpha^N(G_\alpha^Nv_N,\tilde\Psi_N u)\,=\,\lim_{N\to\infty}\,{\langle v_N,\tilde\Psi_N u\rangle}_N
		\,=\,\,{\langle v_\infty, u\rangle}_\infty\,=\,\E_\alpha^\infty(G_\alpha^\infty v_\infty, u)
	\end{equation*}
	for $u\in V$. Then, we obtain $G_\alpha^Nv_N\overset{\textnormal{w.}}{\underset{N}{\rightharpoonup}}G_\alpha^\infty v_\infty$ in the sense of $\mathcal H$, because
	\begin{multline*}
		\lim_{N\to\infty}\,{\langle G_\alpha ^Nv_N,\Psi_N\varphi\rangle}_N\,=\,\lim_{N\to\infty}\,\E_\alpha^N(\,G_\alpha^Nv_N\,,\,\tilde\Psi'_N\,G_\alpha^N\varphi\,)
		\\=\,\E_\alpha^\infty(G_\alpha^\infty v_\infty, G_\alpha^\infty\varphi)\,=\,\,{\langle v_\infty, \varphi\rangle}_\infty
	\end{multline*}
	for $\varphi\in\mathcal C$.
	Now, $G_\alpha^N\overset{\textnormal{s.}}{\underset{N}{\longrightarrow}}G_\alpha^\infty$ is a consequence of \cref{rem:adjoint} and the self-adjointness of $G_\alpha^N$ for $N\in\overline\N$.
	
	Let $C_0{((-\infty,0])}$ denote the space of continuous function vanishing at $-\infty$ and set
	\begin{equation*}
		\mathcal A\,:=\,\{\,f\in C_0{((-\infty,0])}\,|\,f(\Delta_N)\overset{\textnormal{s.}}{\underset{N}{\longrightarrow}}f(\Delta_\infty)\,\}.
	\end{equation*}
	We now come to the final step of this proof. The observation, which completes the implications of (v) by (iv) and also includes the implication of (i) by (v),
	reads as follows:
	If $\mathcal A$ separates points, i.e.~for $t,s\in(-\infty,0]$ with $t\neq s$ there are $f,g\in\mathcal A$ with $f(t)\neq g(s)$, $f(t)\neq 0$, $g(s)\neq 0$,
	then $\mathcal A=C_0((-\infty])$. The observation is an application of the extended Stone-Weierstra{\ss} theorem, as formulated in \cite[Chap. 7]{simmons},
	since $\mathcal A$ is a closed subalgebra of $(C_0((-\infty]),\|\cdot\|_\infty)$. The latter fact can be verified easily via the formula
	\begin{equation*}
		(fg)\,(\Delta_N)\,=\,f(\Delta_N)\,\circ\, g(\Delta_N)
	\end{equation*}
	and the estimate
	\begin{equation*}
		|\,\langle f(\Delta_N)u_N,v_N\rangle_N\,|\,=\,\|f\|_\infty\,{\|u_N\|}_N\,{\|v_N\|}_N
	\end{equation*}
	for $f,g\in C_0((-\infty])$, ${(u_N)}_{N\in\overline \N},{(v_N)}_{N\in\overline \N}\in\mathcal H$ and $N\in\overline\N$.
	This concludes the proof.
\end{proof}

Until the end of \cref{sec:lowse}, let $E$ be a Polish space.
We denote the space of bounded, measurable functions from $E\to\R$ by $\mb(E)$
and the space of bounded, continuous functions from $E\to\R$ by $\cb(E)$. 
Let ${(\mu_N)}_{N\in\N}$ be a sequence of weakly convergent Probability measures on $E$ with limit $\mu_\infty$, i.e.
\begin{equation*}
	\lim_{N\to\infty}\int\limits_E f\de\mu_N\,=\,\int\limits_E f\de\mu_\infty\quad\textnormal{for }f\in\cb(E).
\end{equation*}
We moreover assume for the topological support of the measures
\begin{equation*}
	\supp[\mu_N]\,\subset\,\supp[\mu_\infty]
\end{equation*}
for $N\in\N$. This ensures that the map $\Psi_N$ which sends the $\mu$-class of a bounded, continuous function to its $\mu_N$-class
is well-defined on the linear subspace $\mathcal C:=\cb(E)\cap L^2(E,\mu_\infty)\subset L^2(E,\mu_\infty)$.
Since the asymptotic inequalities \cref{eqn:initop} are fulfilled we are dealing with a sequence
${(L^2(E,\mu_N))}_{N\in\N}$ of converging Hilbert spaces with asymptotic space $L^2(E,\mu_\infty)$. 
Finding a strongly convergent minorante and majorante can be a suitable way of proving that a section of non-negative measurable functions
is strongly convergent and identifying its asymptotic element.

\begin{lem}\label{lem:wemeco}
Let $g_N,f_N^m,F^m_N\in\mb(E)$ for $m\in\N$, $N\in\overline\N$. We assume
\begin{equation*}
	0\,\leq\, f_N^m(x)\,\leq\, g_N(x)\,\leq\, F^m_N(x)
\end{equation*}
for $x\in E$, $N\in\overline\N$ and also the strong convergence of
\begin{equation*}
	\lim_{m\to\infty}f^m_\infty=g_\infty\quad\textnormal{as well as}\quad\lim_{m\to\infty}F^m_\infty=g_\infty
\end{equation*}
in $L^2(E,\mu_\infty)$. The following statement regarding elements of $\prod_{N\in\overline\N}L^2(E,\mu_N)$ holds true: If
\begin{equation*}
	f_N^m\overset{\textnormal{s.}}{\underset{N}{\longrightarrow}}f_\infty^m\quad\textnormal{as well as}\quad 
	F^m_N\overset{\textnormal{s.}}{\underset{N}{\longrightarrow}} F^m_\infty
\end{equation*}
for every $m\in\N$, then also
$g_N\overset{\textnormal{s.}}{\underset{N}{\longrightarrow}}g_\infty$.
\end{lem}
\begin{proof}
Let $g_N, f_N^m,F^m_N$ for $m\in\N$, $N\in\overline\N$ be as in the assumptions of this lemma. In the next steps
$g_N\overset{\textnormal{w.}}{\underset{N}{\rightharpoonup}}  g_\infty$ is shown.
In view of \wesecomp{} and \extopii{} 
we may w.o.l.g.~assume that there exists $h\in L^2(E,\mu_\infty)$ such that
$g_N\overset{\textnormal{w.}}{\underset{N}{\rightharpoonup}} h$. Let $\varphi:E\to[0,\infty)$ be a bounded, continuous function and $m\in\N$. The inequality
\begin{equation*}
	0\leq \int\limits_E\varphi \,f_N^m\de\mu_N\,\leq \int\limits_E\varphi \,g_N\de\mu_N\,\leq\,\int\limits_E\varphi\, F^m_N\de\mu_N
\end{equation*}
for $N\in\N$ leads to the asymptotic inequality 
\begin{equation*}
	0\leq \int\limits_E\varphi\, f_\infty^m\de\mu_\infty\,\leq \int\limits_E\varphi \,h\de\mu_\infty\,\leq\, \int\limits_E\varphi \,F_\infty^m\de\mu_\infty
\end{equation*}
in the limit $N\to\infty$. Hence, $f^m_\infty(x)\leq h(x)\leq F^m_\infty(x)$ holds for $\mu_\infty$-a.e. $x\in E$.
Now passing to the limit $m\to\infty$ implies $g_\infty(x)=h(x)$ for $\mu_\infty$-a.e.~$x\in E$ and hence 
$g_N\overset{\textnormal{w.}}{\underset{N}{\rightharpoonup}} g_\infty$.

To verify the strong convergence we look at the inequality
\begin{equation*}
	0\leq \int\limits_E{(f^m_N)}^2\de\mu_N\,\leq \int\limits_Eg_N^2\de\mu_N\,\leq\,\int\limits_E{(F_N^m)}^2\de\mu_N
\end{equation*}
and by passing to the limit $N\to\infty$ observe that
\begin{equation*}
	\int\limits_E{(f^m_\infty)}^2\de\mu_\infty\,\leq\, \liminf_{N\to\infty}\int\limits_Eg_N^2\de\mu_N\quad\textnormal{as well as}\quad
	\limsup_{N\to\infty}\int\limits_Eg_N^2\de\mu_N\,\leq\,\int\limits_E{(F_\infty^m)}^2\de\mu_\infty.
\end{equation*}
Now, passing to the limit $m\to\infty$ yields
\begin{equation*}
	\limsup_{N\to\infty}\int\limits_Eg_N^2\de\mu_N\,\leq\,\int\limits_Eg_\infty^2\de\mu_\infty\,\leq\,\liminf_{N\to\infty}\int\limits_Eg_N^2\de\mu_N.
\end{equation*}
This concludes the proof.
\end{proof}

\subsection{Convergence of superposed standard gradient forms}\label{sec:lowsii}
Here, we assume that the state space $E$ is given as the product $E=S\times \R^d$, where $d\in\N$ and $S$ is a Polish space.
Denote by $\pi_1:E\to S$ the projection onto the first coordinate.
For $N\in\overline \N$ we define $m_s^N$ as the conditional distribution of $\mu_N$ given $\pi_1=s$ for $s\in S$.
This means by definition that $S\ni s\mapsto m_s^N(V)\in[0,1]$ is measurable for $V\in\mathcal B(\R^d)$ and $\mu_N$ is the superposition of $m_s^N$, $s\in S$, w.r.t.~the image measure $\nu_N$
of $\mu_N$ under $\pi_1$. The equations
\begin{align}\label{eqn:desdef}
\nu_N(\pi_1(A))\,&=\,\mu_N(A),\nonumber\\
\int\limits_S\int\limits_{\R^d} \eins_A(s,x)\de m_s^N(x)\de\nu_N(s)\,&=\,\mu_N(A),\qquad A\in\mathcal B(E),
\end{align}
equivalently characterize the resulting disintegration of $\mu_N$ along $\pi_1$.
The existence and uniqueness of the conditional densities is ensured by a general disintegration theorem as stated in \cite[Theorems 10.2.1 \& 10.2.2]{dudley}.
For simplicity we equivalently write $\mu$ for $\mu_\infty$, $\nu$ for $\nu_\infty$ and $m_s^\infty$ for $m_s$ if $s\in S$.
At the heart of \cref{thm:m1} is the superposition of
standard gradient forms on $L^2(m_s^N)$ w.r.t.~the mixing measure $\de \nu_N(s)$. The bilinear forms of \cref{eqn:dirint} are now lifted to the $L^2$ setting.
As in \cref{sec:mokushi} we want to work with closed forms. That is why in \cref{cond:hamza} we assume Hamza's condition for closability for each $N\in\overline\N$.
The theorem represents the main result of this paper in the abstract setting. 
Mosco convergence for $N\to\infty$ is obtained under some constraints on the conditional distributions. These are stated in \cref{cond:mucken} in terms of the quantities
$C_{r}^{\kappa(s,\cdot)}(m_s^N)$ and $\delta_{r}^{\kappa(s,\cdot)}(m_s^N)$,
depending on $r$ and $\kappa$, as defined in \cref{sec:diri} by \cref{eqn:delta}, respectively \cref{eqn:bic}.

\begin{cond}\label{cond:hamza}
Let $\mu_N$ be a probability measures on $E=S\times\R^d$ for $N\in\overline\N$.
We consider the disintegration according to \cref{eqn:desdef}. 
For $N\in\overline\N$ the family $m^N_s$, $s\in S$, is assumed to meet Hamza's condition in $\nu_N$ -a.e.~sense. This means that $m^N_s$ is absolutely continuous 
w.r.t.~the Lebesgue measure and its density $\varrho_N(s,\cdot)$ fulfills
\begin{equation*}
	\int\limits_S m^N_s\Big(\Big\{x\in\R^d\,\Big|\,\int\limits_{x+[-\varepsilon,\varepsilon]^d}\varrho_N^{-1}(s,y)\de y <\infty
	\textnormal{ for some }\varepsilon>0\Big\}\Big)\de\nu_N(s)\,=\,1.
\end{equation*}
\end{cond}
\newcommand{\hamza}{\cref{cond:hamza}}

\cref{thm:m1} identifies the Mosco limit for a sequence of Dirichlet forms.
Let $N\in\overline \N$. 
\cref{cond:hamza} says that for $\nu_N$ -a.e.~$s\in S$ there is an open set $U^{\varrho_N}_s\subset \R^d$ such that $x\mapsto \varrho_N^{-1}(s,x)$ is locally $\de x$ -integrable on $U^{\varrho_N}_s$ and 
$\varrho_N(s,x)=0$ holds $\de x$ -a.e.~on $\R^d\setminus U^{\varrho_N}_s$. By the Cauchy-Schwarz inequality
\begin{equation}\label{eqn:locemb}
L^2(S\times \R^d,\mu_N)\,\hookrightarrow\,L_\textnormal{loc}^1\big(\{\,(s,x)\,|\,s\in S,x\in U^{\varrho_N}_s\,\},\nu_N\times \de x\big)
\end{equation}
is continuously embedded.
Thanks to \cref{sec:est} we define a pre-domain
$\dom_{\textnormal{pre}}(\E^N)$ comprising elements of $u,v\in L^2(E,\mu_N)$ with representatives $\tilde u,\tilde v$ from $\dom(\int^\oplus E^{\varrho_N(s,\cdot)}\de\nu_N(s))$,
and a symmetric, non-negative bilinear form 
\begin{equation}\label{eqn:predi}
\E^N(u,v)\,:=\,\big(\textstyle{\int^\oplus} E^{\varrho_N(s,\cdot)}\de\nu_N(s)\big)\,(\tilde u,\tilde v),\qquad u,v\in\dom_{\textnormal{pre}}(\E).
\end{equation}
Due to \cref{cond:hamza} the form $(\E^N,\dom_{\textnormal{pre}}(\E^N))$ is well-defined and closable on $L^2(E,\mu_N)$ following
standard arguments, as disclosed in \cite[Chapter 3]{fuku94} and \cite[Chapter 2]{maro} and  - in particularly within the context of superposition of forms - also by \cite{albe}.
Its smallest closed extension on $L^2(E,\mu_N)$ is denoted by $(\E^N,\dom(\E^N))$.
If ${(u_m)}_{m\in\N}\subset \dom_\textnormal{pre}(\E^N)$ is 
an approximating Cauchy sequence for an element $u$ in $(\dom(\E^N),\E^N_1)$ and $i\in\{1,\dots,d\}$, then 
the family of functions
$(s,x)\mapsto \partial_{x_i}u_m(s,x)$, $m\in\N$, form a Cauchy sequence in $L^2(E,\mu_N)$.
From \cref{eqn:locemb} we deduce $u(s,\cdot)\in H^{1,1}_\textnormal{loc}(U^{\varrho_N}_s)$ for $\nu_N$ -a.e.~$s\in S$ and
\begin{equation}\label{eqn:conf}
\E^N(u,u)\,=\,\sum_{i=1}^d\,\int\limits_S\int\limits_{\R^d}|\partial_{x_i}u(s,x)|^2\de m_s^N(x)\de\nu_N(s)
=\,\sum_{i=1}^d\,\int\limits_{S\times\R^d}|\partial_{x_i}u(s,x)|^2\de \mu_N(s,x).
\end{equation}
The contraction property,
\begin{gather*}
0\lor(\eins\land u)\in\dom_\textnormal{pre}(\E^N)\quad\textnormal{for}\quad u\in\dom_\textnormal{pre}(\E^N)\\
\textnormal{with}\quad\E^N(\,0\lor(\eins\land u)\,,\,0\lor(\eins\land u)\,)\,\leq\,\E^N(u,u),
\end{gather*}
is inherited by $(\E^N,\dom(\E^N))$, which makes it a Dirichlet form. 
As is quite common, `$f\lor g$' (or `$f\land g$') denotes the maximum (respectively the minimum) of two measurable functions or the $\mu_N$-classes of such.
Consequently, the associated strongly continuous contraction resolvent ${(G_\alpha^N)}_{\alpha>0}$ is sub-Markovian, i.e.~
\begin{equation*}
0\,\leq\,\alpha\,G^N_\alpha u\,(\cdot)\,\leq 1\,\mu_N\textnormal{ -a.e.}\quad\textnormal{if}\quad0\,\leq\, u(\cdot)\,\leq\,1\,\mu_N\textnormal{ -a.e.}
\end{equation*}
for $u\in L^2(E,\mu_N)$ and $\alpha>0$.
A profound survey about the concept of Markovianity is provided by \cite[Chapters 1 \& 2]{maro}.
Replacing $E^{\varrho_N(s,\cdot)}$ by $E^{\kappa\varrho_N(s,\cdot)}$ in \cref{eqn:predi} we define a Dirichlet form $(\E^{N,\kappa},\dom(\E^{N,\kappa}))$ 
on $L^2(E,\kappa\mu_N)$ in an analogous way for  
a measurable function $\kappa:S\times\R^d\to[0,1]$. We note that $\E^{N,\kappa}$ is dominated by $\E^N$ 
in the sense that the natural linear inclusion $L^2(E,\mu_N)\to L^2(E,\kappa\mu_N)$, which sends an element $u\in L^2(E,\mu_N)$ to the $\kappa\mu_N$-class of a representative $\tilde u$,
restricts to a linear map $(\dom(\E^N),\E^N_1)\to(\dom(\E^{N,\kappa}),\E^{N,\kappa}_1)$
with operator norm smaller or equal $1$.

As a preparation for the proof of \cref{thm:m1} we extend the scope of \cref{lem:polygon} to the whole of $\dom(\E^N)$ for $N\in\N$.
This can be achieved in a straight-forward way. 
\begin{lem}\label{rem:polygon}
Let $N\in\N$, $\kappa:S\times \R^d\to [0,1]$ be a measurable function and $h:(s,x)\mapsto f(s)\,g(x)$ for $f\in\mb(S)$ and $g\in\cic(\R^d)$.
For $u\in\dom (\E^N)$ with $-1\leq u(\cdot)\leq 1$ and $r\in(0,\infty)$ there exists $\lambda\in\mathscr L_r$
such that each of the inequalities holds true.\vspace{1ex}
\begin{thmlist}[wide, labelindent=0pt]
	\item\label{rem:polygon0} $|\lambda^\alpha(\,\cdot\,)|\,\leq\,1$
	where $\lambda^\alpha$ is the coefficient of $\lambda$ with index $\alpha\in r\,\Z^d$.\vspace{1ex}
	\item\label{rem:polygon1} $\displaystyle\Big|\int\limits_{S\times\R^d}h\,(u-\lambda)\,\kappa\de \mu_N\,\Big|$
	$\,\displaystyle\leq\,\|f\|_{L^\infty(\nu_N)}\,\big(\omega_r^g+\|g\|_\infty\,\|\delta^{\kappa(s,\cdot)}_r(m^N_s)\|_{L^2(\nu_N)}\big)\,.$\vspace{1ex}
	\item\label{rem:polygon2} 
	$\displaystyle \E^{N,\kappa}(h,u-\lambda)$
	$\,\displaystyle\leq\,\sqrt{d}\,\|f\|_{L^\infty(\nu_N)}\,\big(\omega_r^{\nabla g}+D^g\,\|\delta^{\kappa(s,\cdot)}_r(m^N_s)\|_{L^2(\nu_N)}\big)\,
	{\E^N(u,u)}^{\frac{1}{2}}\,.$\vspace{1ex}
	\item\label{rem:polygon3} 
	$\displaystyle\E^{N,\kappa}(\lambda,\lambda)$
	$\,\displaystyle\leq\,\|C^{\kappa(s,\cdot)}_r(m_s^N)\|_{L^\infty(\nu_N)}\,\E^{N}(u,u).$\vspace{1ex}
\end{thmlist}
\end{lem}
\begin{proof}
Let $u\in\dom(\E^N)$ with $-1\leq u(\cdot)\leq 1$ and $r\in(0,\infty)$. 
Since $(\E^N,\dom_\textnormal{pre}(\E^N))$ is a pre-Dirichlet form and $\dom_{\textnormal{pre}}(\E^N)\subset (\dom(\E^N),\E^N_1)$ densely, 
there exists $u_m\in \dom_\textnormal{pre}(\E)$ for $m\in\N$ such that $-1\leq u_m(\cdot)\leq 1$ holds $\mu_N$ -a.e.~and $\lim_{m\to\infty}\E^N_1(u_m-u,u_m-u)=0$.
We apply \cref{lem:polygon} on a suitable representative $\tilde u_m$ of $u_m$ for $m\in\N$, which returns us an approximation $\lambda_m\in\mathscr L_r$, say
\begin{equation*}
	\lambda_m(s,x)\,=\,\sum_{\alpha\in r\Z^d}\lambda_m^\alpha(s)\,\chi_r^\alpha(x),\quad s\in S,x\in\R^d,
\end{equation*}
according to \cref{eqn:polygon0,eqn:polygon1,eqn:polygon2,eqn:polygon3}. Now we use the weak sequential compactness of bounded sets in $L^2(S,\nu_N)$.
Repeatedly dropping to a suitable subsequence and forming a diagonal sequence, we may w.l.o.g.~assume 
the existence of measurable functions $\lambda^\alpha$, $\alpha\in r\,\Z^d$, with $-1\leq\lambda^\alpha(\cdot)\leq 1$ and 
\begin{equation*}
	\lim_{m\to\infty}\int\limits_S(\lambda_m^\alpha-\lambda^\alpha)\, h\de\nu_N=0\qquad\textnormal{for}\quad h\in L^2(\nu_N).
\end{equation*}

Let now ${(R_\alpha)}_{\alpha>0}$ denote the resolvent of $(\E^{N,\kappa},\dom(\E^{N,\kappa}))$.
Exploiting the weak sequential compactness of bounded sets in $(\dom(\E^{N,\kappa}),\E^{N,\kappa}_1)$
and the energy bounds for ${(\lambda_m)}_m$ due to \cref{eqn:polygon3} we deduce the weak convergence of ${(\lambda_m)}_{m\in\N}$ in $(\dom(\E^{N,\kappa}),\E^{N,\kappa}_1)$.
Indeed, every weak accumulation point $w$ in $(\dom(\E^{N,\kappa}),\E^{N,\kappa}_1)$ must coincide with the function
\begin{equation*}
	\lambda:(s,x)\,\mapsto\,\sum_{\alpha\in r\Z^d}\lambda_m^\alpha(s)\,\chi_r^\alpha(x),
\end{equation*}
because for $v\in \mb(E)$ and each bounded, measurable set $K\subset \R^d$ we have
\begin{multline*}
	\int\limits_{S\times K} \lambda\,v\,\kappa\de\mu_N\,=\,\sum_{\alpha\in r\Z^d}\int\limits_S\lambda^\alpha(s)\,\int\limits_{K}\chi_r^\alpha\,v(s,\cdot)\,\kappa(s,\cdot)\de m_s^N\de\nu_N(s)\\
	=\,\lim_{m\to\infty}\sum_{\alpha\in r\Z^d}\int\limits_S\lambda_m^\alpha(s)\,\int\limits_{K}\chi_r^\alpha\,v(s,\cdot)\,\kappa(s,\cdot)\de m_s^N\de\nu_N(s)
	=\,\lim_{m\to\infty}\,\int\limits_{S\times K}\lambda_m\,v\,\kappa\de\mu_N\\=
	\,\lim_{m\to\infty}\,\E^{N,\kappa}_1(\lambda_m,R_1 (\eins_{S\times K}\,v))\,=\,\E^{N,\kappa}_1(w,R_1 (\eins_{S\times K}\,v))\,=\,\int\limits_{S\times K} w\,v\,\kappa\de\mu_N.
\end{multline*}
In particular, we have
\begin{equation*}
	\E^{N,\kappa}(\lambda,\lambda)\leq\liminf_{m\to\infty}\,\E^{N,\kappa}(\lambda_m,\lambda_m).
\end{equation*}

Now, the claimed estimates of (ii) to (iv) regarding $\lambda$ and $u$ emerge from
their analogues of \cref{eqn:polygon1,eqn:polygon2,eqn:polygon3} for their approximations $\lambda_m$ and $\tilde u_m$, when considering the limit $m\to\infty$.
\end{proof}

We denote by $\dom_\textnormal{min}(\E^\infty)$ the subspace of $\dom(\E^\infty)$ which is the topological closure in $(\dom(\E^\infty),\E^{\infty}_1)$ of the set
comprising all elements with representative in 
\begin{equation*}
	\cb(S)\otimes \cic(\R^d):=\textnormal{span}\Big(\Big\{\,S\times\R^d\ni(s,x)\,\mapsto f(s)\,g(x)\,\Big|\, f\in\cb(S),g\in\cic(\R^d)\,\Big\}\Big).
\end{equation*}
The strongly continuous contraction resolvent of the Dirichlet form $(\E^\infty,\dom_\textnormal{min}(\E^\infty))$ on $L^2(E,\mu_\infty)$ is denoted by ${(G_\alpha)}_{\alpha>0}$.

Analogously, we define the Dirichlet form $(\E^{\kappa,\infty},\dom_\textnormal{min}(\E^{\kappa,\infty}))$ on $L^2(E,\kappa\mu_\infty)$ for a measurable function $\kappa:S\times\R^d\to[0,1]$.
Again, we remark that $\E^{\infty,\kappa}$ is dominated by $\E^{\infty}$,
meaning the natural linear inclusion $L^2(E,\mu_\infty)\to L^2(E,\kappa\mu_\infty)$ restricts to a map $(\dom(\E^\infty),\E^\infty_1)\to(\dom(\E^{\infty,\kappa}),\E^{\infty,\kappa}_1)$
with operator norm smaller or equal $1$.
For short we equivalently write $\E$ for $\E^\infty$.

\begin{cond}\label{cond:mucken}
Let ${(\mu_N)}_{N\in\N}$ be a sequence of weakly converging probability measures on $E=S\times\R^d$ with limit $\mu:=\mu_\infty$,
hence we denote $\nu:=\nu_\infty$, $m_s:=m_s^\infty$ and $\E:=\E^\infty$.
For their disintegrations according to \cref{eqn:desdef} we assume the following.
\begin{thmlist}[wide, labelindent=0pt]
	\item \label{cond:desint}\vspace{1ex}
	$\displaystyle
	\lim_{N\to\infty}\,\int\limits_S\Big|\int\limits_{\R^d} g(s,\cdot)\de m_s^N\,\Big|^2\de\nu_N(s)\,=\,\int\limits_S\Big|\int\limits_{\R^d} g(s,\cdot)\de m_s\,\Big|^2\de\nu(s)$
	\qquad for $g\in\cb(E)$.\vspace{1ex}
	\item \label{cond:liminf}\label{lem:dommin} There exists an at most countable set $\mathscr U$ of continuous functions $E\to[0,1]$ such that
	\begin{equation}\label{eqn:liminf}
		\lim_{m\in\N}\,\sup_{N\in\N}\|\delta_{1/m}^{\kappa(s,\cdot)}(m_s^N)\|_{L^2(\nu_N)}=0\quad\textnormal{and}\quad
		\limsup_{m\in\N}\,\sup_{N\in\N}\|C_{1/m}^{\kappa(s,\cdot)}(m_s^N)\|_{L^\infty(\nu_N)}\,<\,\infty.
	\end{equation}
	for each element $\kappa\in\mathscr U$. In addition to that,
	\begin{gather}
		\sup_{\kappa\in\mathscr U}\kappa(s,x) \,=\, 1\textnormal{ for }\mu\textnormal{-a.e. }(s,x)\in S\times\R^d\textnormal{ and}\nonumber\\
		\dom_\textnormal{min}(\E)\,=\,\Big\{\,u\in L^2(E,\mu)\,\Big|\,u\in\bigcap_{\kappa\in\mathscr U}\dom_\textnormal{min}(\E^{\infty,\kappa})
		\textnormal{ with }\sup_{\kappa\in\mathscr U}\E^{\infty,\kappa}(u,u)<\infty\,\Big\}.\label{eqn:dommin}
	\end{gather}	    
\end{thmlist}
\end{cond}
\newcommand{\desint}{\cref{cond:desint}}
\newcommand{\mucken}{\cref{cond:liminf}}
We comment on the first item of the condition. A discussion about the second item follows after \cref{thm:m1}.
\begin{rem}\label{rem:onmucki}
\label{rem:desrem}Let $g\in\cb(E)$ and $f_N(s)$ $:=\int_{\R^d} g(s,\cdot)\de m_s^N$ for $s\in S$ and $N\in\overline \N$.
\desint{}  is equivalent to $f_N\overset{s.}{\underset{N}{\longrightarrow}}f_\infty$ referring to
\sloppy $\prod_{N\in\overline\N}L^2(S,\nu_N)$, since $f_N\overset{w.}{\underset{N}{\rightharpoonup}}f_\infty$ is already implied by the weak measure convergence of $\mu_N$ towards $\mu$.
\end{rem}
\newcommand{\desrem}{\cref{rem:desrem}}

Under \cref{cond:mucken,cond:hamza} we obtain our main abstract result on Mosco convergence.
The observation of $\mathscr L_r\subset\dom_\textnormal{min}(\E)$ for $r\in(0,\infty)$ is vital for the proof. 
Te matter basically boils down to a generally known result about the minimal domain of a gradient-type Dirichlet form containing the Lipschitz continuous functions.
We manifest this fact in a lemma before we state the theorem.

\begin{lem}\label{lem:lipmin}
$\mathscr L_r\subset\dom_\textnormal{min}(\E)$ for $r\in(0,\infty)$. 	
\end{lem}
\begin{proof}
The proof works is a double application of \cite[Lemma 2.12]{maro}. It provides us with a sufficient criterion for an element $u\in L^2(E,\mu)$ to be a member of $\dom_{\textnormal{min}}(\E)$:
The existence of a sequence ${(u_k)}_{k\in\N}$ $\subset\dom_{\textnormal{min}}(\E)$ such that
\begin{equation}\label{eqn:rolem}
	\lim_{k\to\infty}\,\int\limits_E|u_k-u|^2\de\mu\,=\,0\qquad\textnormal{and}\qquad\sup_{k\in\N}\E(u_k,u_k)<\infty.
\end{equation}
Let $r\in(0,\infty)$, $\alpha\in r\,\Z^d$ and $f\in\mb(S)$.
Further, let ${(f_k)}_{k\in\N}\subset\cb(S)$ be an approximation for $f$ such that $f_k\to f$ in $\nu$ -a.e.~sense and $\sup_k\|f_k\|_\infty=:C<\infty$.
We choose a non-negative function $\varphi_k\in\cic((-1/k,1/k))$ with $\int_{\R}\varphi_k\de x=1$ and set 
\[u_k:S\times\R^d\ni (s,x)\,\mapsto\,f_k(s)\,(\varphi_k*\chi^\alpha_r)(x)\]
for $k\in\N$, the symbol `$*$' denoting the convolution.
We note that $\chi_r^\alpha$ is globally Lipschitz continuous with constant smaller equal $\sqrt{2}/r$ by virtue of \cref{thm:partuniii}.
Then, \[u:S\times\R^d\ni (s,x)\,\mapsto\,f(s)\,\chi_r^\alpha(x)\] defines an element of $\dom_{\textnormal{min}}(\E)$, because \cref{eqn:rolem} 
is verified with Lebesgue's dominated convergence and the estimate
\begin{multline*}
	\sup_{(s,x)\in E}|\partial_{x_i} u_k(s,x)|\,\leq\,C\,\sup_{b\in\R}\,\Big|\,\frac{1}{b}\,
	\int\limits_{\R^d}\varphi_k(y)\,\chi_r^\alpha(x-y)\de y-\int\limits_{\R^d}\varphi_k(y)\,\chi_r^\alpha(x+b\,\univ_i-y)\de y\,\Big|\\
	=\,C\,\sup_{b\in\R}\,\int\limits_{\R^d}\varphi_k(y)\,\Big|\,\frac{\chi_r^\alpha(x-y)-\chi_r^\alpha(x+b\,\univ_i-y)}{b}\,\Big|\,\de y\,\leq\,\frac{\sqrt{2}\,C}{r}
\end{multline*}
for $k\in\N$.
Let $M\in(0,\infty)$ and $\lambda^\alpha:S\to[-M,M]$ measurable functions for $\alpha\in r\,\Z^d$. The claim now follows, since again by Lebesgue's dominated convergence 
\begin{equation*}
	\lim_{k\to\infty}\,\int\limits_E\Big|\,\sum_{\alpha\in r\Z^d\setminus [-k,k]^d}\lambda^\alpha(s)\,\chi_r^\alpha(x)\,\Big|^2\de\mu(s,x)\,=\,0,
\end{equation*}
while we estimate 
\begin{equation*}
	\sup_{s\in S}\,\Big\|\sum_{\alpha\in r\Z^d\cap[-k,k]^d}\lambda^\alpha(s)\,\partial_i\chi_r^\alpha\,\Big\|_{L^\infty(\R^d,m_s)}
	\leq\,\frac{2\,M}{r}
\end{equation*}
for $k\in\N$ and $i=1,\dots,d$ with \cref{thm:partuniii}.
\end{proof}

\begin{thm}\label{thm:m1}
	Let \cref{cond:mucken,cond:hamza} be fulfilled. 
	$(\E^N,\dom(\E^N))$ converges to $(\E,\dom_\textnormal{min}(\E))$ in the sense of Mosco.
\end{thm}
\begin{proof} We want to verify \amind. Let $u,v\in$ $\cb(S)\otimes \cic(\R^d)$. Due to the weak convergence of measures we have
	\begin{align}
		\lim_{N\to\infty}\E^{N,\kappa}(u,v)\,&=\,\lim_{N\to\infty}\,\sum_{i=1}^d\,\int\limits_{S\times\R^d}\partial_{x_i}u(s,x)\,\partial_{x_i}v(s,x)\,\kappa(s,x)\de \mu_N(s,x)\nonumber\\
		&=\,\sum_{i=1}^d\,\int\limits_{S\times\R^d}\partial_{x_i}u(s,x)\,\partial_{x_i}v(s,x)\,\kappa(s,x)\de \mu(s,x)\,=\,\E^{\infty,\kappa}(u,v)\label{eqn:rev3}
	\end{align}
	as well as
	\begin{equation}
		\lim_{N\to\infty}\,\int\limits_{S\times\R^d}u\,v\,\kappa\de \mu_N\,=\,\int\limits_{S\times\R^d}u\,v\,\kappa\de \mu.\label{eqn:rev4}
	\end{equation}
	So, with the choice $\kappa=\eins_E$ Property (b) is satisfied. 
	
	Property (a) is left to proof. We fix $f\in\cb(S\times\R^d)$. Since \cref{cond:mucken} is stable under dropping to subsequences, it suffices to show a modified version of \amind{} (a):
	\begin{equation*}
		u^*\in\dom_{\textnormal{min}}(\E)\qquad\textnormal{and}\qquad\E(u^*,u^*)\,\leq\,\liminf_{N\to\infty}\E^N(G_1^Nf,G_1^Nf)
	\end{equation*}
	under the condition that 
	\begin{equation*}
		u^*\in L^2(E,\mu)\qquad\textnormal{with}\qquad G_1^Nf\overset{\textnormal{w.}}{\underset{N}{\rightharpoonup}}u^*\quad\textnormal{in the sense of }\prod_{N\in\overline\N}L^2(E,\mu_N).
	\end{equation*}

	In case $f\equiv 0$ there is nothing to do. Otherwise we can equivalently consider $f/\|f\|_\infty$, respectively $u^*/\|f\|_\infty$, instead of $f$ and $u^*$.
	So, we assume $-1\leq f(\cdot)\leq 1$, and hence $-1\leq G_1 ^Nf(\cdot)\leq 1$ for $N\in\N$ due to the sub-Markovian property.
	The strategy to obtain the desired statement is to interpret ${({G_1^N}f)}_{N\in\N}$ as an element from the product space of Hilbert spaces
	\begin{align*}
		&\mathcal H^{\E,\kappa,\beta}=\prod_{N\in\N}\dom(\E^{N,\kappa}),\\
		&\quad\textnormal{where }\dom(\E^{N,\kappa})\textnormal{ is equipped with the scalar product }\E_\beta^{N,\kappa}\textnormal{ for }N\in\N,
	\end{align*}
	with arbitrary $\beta>0$ and $\kappa\in\mathscr U$.
	For such $\kappa$ and $\beta$ the Hilbert spaces $(\dom(\E^{N,\kappa}),\E_\beta^{N,\kappa})$, $N\in\N$, converge to the asymptotic space 
	$(\dom_\textnormal{min}(\E^{\infty,\kappa}),\E_\beta^{\infty,\kappa})$ due to \cref{eqn:rev3,eqn:rev4}.

	The main part of this proof is dedicated to show 
	\begin{equation}\label{eqn:dedicate}
		u^*\in\dom_{\textnormal{min}}(\E^{\infty,\kappa})\qquad\textnormal{and}\qquad
		G_1^Nf\overset{\textnormal{w.}}{\underset{N}{\rightharpoonup}}u^*\textnormal{ in the sense of }\mathcal H^{\E,\kappa,\beta}
	\end{equation}
	for every $\kappa\in\mathscr U$ and $\beta>0$.
	Once this is achieved, the modified version of Property (a) of \amind{} follows easily from \cref{eqn:dedicate} under usage of \lose{}. Indeed, in that case
	\begin{align*}
		\sup_{\kappa\in\mathscr U}\,\E_\beta^{\infty,\kappa}(u^*,u^*)\,&\leq\,\sup_{\kappa\in\mathscr U}\,\liminf_{N\to\infty} \E^{N,\kappa}_\beta(G_1^Nf,G_1 ^Nf)\\
		&\leq\,\liminf_{N\to\infty} \E^{N}_\beta(G_1^Nf,G_1 ^Nf).
	\end{align*}
	Now, for one thing \cref{lem:dommin} implies $u^*\in\dom_{\textnormal{min}}(\E)$, and for the other  
	\begin{equation*}
		\E(u^*,u^*)\,\leq\,\liminf_{N\to\infty}\E^N(G_1^Nf,G_1^Nf)
	\end{equation*}
	considering the limit $\beta\to0$.

	So, let $\beta>0$ and $\kappa\in\mathscr U$ be fixed. Only \cref{eqn:dedicate} is left to show.
	By virtue of \wesecomp{} and \extopii{} we may w.l.o.g.~assume that there exists $w^*_\kappa\in\dom_\textnormal{min}(\E^{\infty,\kappa})$ such that
	\begin{equation}\label{eqn:heki}
		G_1^Nf\overset{\textnormal{w.}}{\underset{N}{\rightharpoonup}}w_\kappa^*\quad\textnormal{in the sense of }\mathcal H^{\E,\kappa,\beta}.
	\end{equation}

		First, we prove a related statement considering only sections with elements from $\mathscr L_r$ for fixed $r>0$.
		Let $u_N\in\mathscr L_r$ such that
		\begin{equation*}
			u_N(s,x)\,=\,\sum_{\alpha\in r\Z^d} q_N^\alpha(s)\,\chi_r^{\alpha}(x),\quad s\in S,\,x\in \R^d,
		\end{equation*}
		with measurable coefficient functions $q_N^\alpha:S\to[-1,1]$ for $\alpha\in r\Z^d$ and $N\in\N$. Moreover, let
		\begin{align*}
			&u^{**}\in L^2(\kappa \mu):&&u_N\overset{\textnormal{w.}}{\underset{N}{\rightharpoonup}}u^{**}\quad\textnormal{in the sense of }\prod_{N\in\overline \N}L^2(E,\kappa\mu_N)\textnormal{ and}\\
			&w^{**}\in\dom_\textnormal{min}(\E^{\infty,\kappa}):&&u_N\overset{\textnormal{w.}}{\underset{N}{\rightharpoonup}}w^{**}\quad\textnormal{in the sense of }\mathcal H^{\E,\kappa,\beta}.
		\end{align*}

		The goal is now to show the identity $u^{**}=w^{**}$ in $L^2(E,\kappa\mu)$.
		By repeated usage of \wesecomp{} we can - after dropping to suitable diagonal subsequence - w.l.o.g.~assume the existence of
		\begin{equation*}
			q^\alpha\in L^2(S,\nu),\,\alpha\in r\Z^d:\quad  q^\alpha_N\overset{\textnormal{w.}}{\underset{N}{\rightharpoonup}} 
			q^\alpha\quad\textnormal{in the sense of }\prod_{N\in\overline \N}L^2(S,\nu_N).
		\end{equation*}
		Due to \desint{} and \cref{rem:onmucki} we deduce
		\begin{equation*}
			u^{**}(s,x)=\sum_{\alpha\in r\Z^d} q^\alpha(s)\,\chi_r^{\alpha}(x)\quad\textnormal{for }\kappa\mu\textnormal{ -a.e.~}s\in S,\,x\in\R^d,
		\end{equation*}
		as
		\begin{multline*}
			\int\limits_{S\times\R^d}u^{**}(s,x)\,g(s)\,\varphi(x)\,\kappa(s,x)\de\mu(s,x)\\
			\begin{aligned}&=\,
				\lim_{N\to\infty}\,\sum_{\alpha\in r\Z^d}\int\limits_Sq_N^\alpha(s)\,g(s)\int\limits_{\R^d}\chi_r^\alpha\,\varphi\,\kappa(s,\cdot)\de m_s^N\de\nu_N(s)\\
				&=\,\sum_{\alpha\in r\Z^d}\int\limits_Sq^\alpha(s)\,g(s)\int\limits_{\R^d}\chi_r^\alpha\,\varphi\,\kappa(s,\cdot)\de m_s\de\nu(s).
			\end{aligned}
		\end{multline*}
		holds true for $g\in\cb (S)$ and $\varphi\in\cc(\R^d)$. We note that the summation over $\alpha$ in this equation is actually a finite sum.

		Now, we set
		\begin{equation*}
			U\:=\,\bigcup_{T\in\mathcal T_r}\accentset{\circ}{D_T}.
		\end{equation*}
		Since $U$ is an open set in $\R^d$ and  $m_s (\R^d\setminus U)$ $=0$ holds for every $s\in S$, the
		linear span of the product indicator functions  from the family
		\begin{equation*}
			\Big\{\,S\times \R^d\ni (s,x)\,\mapsto\, \eins_{A}(s)\,\eins_K(x)\,\Big|\,A\in\mathcal B(S), K\subset U\textnormal{ and } K\textnormal{ is a compact set}\,\Big\}
		\end{equation*}
		is dense in $L^2(E,\kappa \mu)$.
		Then, by approximation, a dense subspace of $L^2(E,\kappa\mu)$ is also induced by $\cb(S)\otimes\cc(U)$, i.e.~the linear span of all functions
		$(s,x)\,\mapsto\, g(s)\,\varphi(x)$, where $g\in C_b(S)$ and $\varphi$ is a continuous function with compact support contained in $U$.
		
		Let $i=1,\dots,d$. Under use of \desint{}, \cref{rem:onmucki} and \partuniii{} we obtain
		\begin{multline*}
			\lim_{N\to\infty}\,\sum_{\alpha\in r\Z^d}\,\int\limits_Sq_N^\alpha(s)\,g(s)\int\limits_{\R^d}\partial_i\chi_r^\alpha\,\varphi\,\kappa(s,\cdot)\de m_s^N\de\nu_N(s)\\
			=\,\lim_{N\to\infty}\,\frac{1}{r}\sum_{\alpha\in r\Z^d}\,\sum_{\substack{T\in\mathscr T_r:\\T(\sigma^{-1}_T(i)-1)=\alpha}}\,
			\int\limits_S(\,q^{\alpha+r\univ_i}_N(s)-q^{\alpha}_N(s)\,)\,g(s)\int\limits_{\R^d}\eins_{D_T}\,\varphi\,\kappa(s,\cdot)\de m_s^N\de\nu_N(s)\\
			=\,\frac{1}{r}\sum_{\alpha\in r\Z^d}\,\sum_{\substack{T\in\mathscr T_r:\\T(\sigma^{-1}_T(i)-1)=\alpha}}\,
			\int\limits_S(\,q^{\alpha+r\univ_i}(s)-q^{\alpha}(s)\,)\,g(s)\int\limits_{\R^d}\eins_{D_T}\,\varphi\,\kappa(s,\cdot)\de m_s\de\nu(s)\\
			=\,\sum_{\alpha\in r\Z^d}\int\limits_Sq^\alpha(s)\,g(s)\int\limits_{\R^d}\partial_i\chi_r^\alpha\,\varphi\,\kappa(s,\cdot)\de m_s\de\nu(s)
		\end{multline*}
		for $g\in\cb (S)$ and $\varphi\in\cc(U)$. 
		Again, the summation over $\alpha$ in this equation is actually a finite sum.
		So, referring to $\prod_{N\in\overline\N}L^2(E,\kappa\mu_N)$, the function
		\begin{equation*}
			S\times \R^d\ni (s,x)\,\mapsto\,\sum_{\alpha\in r\Z^d} q^\alpha(s)\,\partial_i\chi_r^\alpha(x)
		\end{equation*}
		is the asymptotic element of 
		\begin{equation*}
			S\times \R^d\ni (s,x)\,\mapsto\,\sum_{\alpha\in r\Z^d} q_N^\alpha(s)\,\partial_i\chi_r^\alpha(x)
		\end{equation*}
		in the terminology of a weakly convergent section.  
		Now, to prove the identity $u^{**}=w^{**}$, let $v\in$ $\cb(S)\otimes \cic(\R^d)$. Then
		\begin{align*}
			\E^{\infty,\kappa}_\beta(u^{**},v)\,&=\,\sum_{i=1}^d\,\int\limits_{S\times\R^d}\sum_{\alpha\in r\Z^d}q^\alpha(s)\,\partial_i\chi_r^{\alpha}(x)\,\partial^x_iv(s,x)\,\kappa(s,x)\de\mu(s,x)\\
			&\begin{aligned}
				&\qquad\qquad\qquad\qquad+\beta\int\limits_{S\times\R^d} \sum_{\alpha\in r\Z^d}q^\alpha(s)\,\chi_r^{\alpha}(x)\,v (s,x)\,\kappa(s,x)\de\mu(s,x)\\
				&=\,\sum_{i=1}^d\,\lim_N\,\int\limits_{S\times\R^d}\sum_{\alpha\in r\Z^d}q_N^\alpha(s)\,\partial_i\chi_r^{\alpha}(x)\,\partial_i^xv(s,x)\,\kappa(s,x)\de\mu_N(s,x)\\
				&\qquad\qquad\qquad\qquad+\beta\,\lim_N\,\int\limits_{S\times\R^d}\sum_{\alpha\in r\Z^d} q_N^\alpha(s)\,\chi_r^{\alpha}(x)\,v (s,x)\,\kappa(s,x)\de\mu_N(s,x)\\
				&=\,\lim_N\,\E_\beta^{N,\kappa}(u_N,v)
				\,=\,\E_\beta^{\infty,\kappa}(w^{**},v)\,.
			\end{aligned}
		\end{align*}
		This yields the identity $u^{**}=w^{**}$ in $L^2(E,\kappa\mu)$.
		
		To bridge the gap and come back to the problem of \cref{eqn:dedicate} we choose an approximation $\lambda_r^N\in\mathscr L_r$ for $G_1^Nf$ according to \cref{rem:polygon} for $N\in\N$ and $r>0$.
		By \cref{rem:polygon3} we estimate
		\begin{equation}\label{eqn:rbound}
			\sup_{N\in\N}\,\E^{N,\kappa}_\beta(\lambda_r^N,\lambda_r^N)\,\leq\,\sup_{N\in\N}\,\Big(\,\beta\,\|\lambda_r^N\|_\infty
			+\|C^{\kappa(s,\cdot)}_r(m^N_s)\|_{L^\infty(\nu_N)}\,\E^N(u_N,u_N)\,\Big)
		\end{equation}
		for $r>0$. For every $m\in\N$ the right hand side of \cref{eqn:rbound} takes a finite value w.r.t.~the choice $r_m:=1/m$.
		So, by repeated usage of \wesecomp{} and
		dropping to a suitable diagonal sequence, we obtain $N_k\in\overline\N$, strictly increasing in $k\in\overline\N$, such that there exists
		\begin{align}
			&\lambda_m^*\in \dom_{\textnormal{min}}(\E^{\infty,\kappa}):
			&&\lambda_{r_m}^{N_k}\overset{\textnormal{w.}}{\underset{k}{\rightharpoonup}}\lambda^*_m\quad\textnormal{in the sense of }\prod_{N\in\overline \N}L^2(E,\kappa\mu_N)\textnormal{ and}\nonumber\\
			& &&\lambda_{r_m}^{N_k}\overset{\textnormal{w.}}{\underset{k}{\rightharpoonup}}\lambda^*_m\quad\textnormal{in the sense of }\mathcal H^{\E,\kappa,\beta}\label{eqn:hekii}
		\end{align}
		for every $m\in\N$. Moreover, 
		\begin{multline}\label{eqn:lambound}
			\limsup_{m\to\infty}\,\E^{\infty,\kappa}_\beta(\lambda^*_m,\lambda^*_m)\,\leq\,\limsup_{m\to\infty}\,\liminf_{k\to\infty}\,\E^{{N_k},\kappa}_\beta(\lambda_{r_m}^{N_k},\lambda_{r_m}^{N_k})\\
			\leq\,\limsup_{m\in\N}\sup_{N\in\N}\,\Big(\,\beta\,\|\lambda_{1/m}^N\|_\infty
			+\|C^{\kappa(s,\cdot)}_{1/m}(m^N_s)\|_{L^\infty(\nu_N)}\,\E^N(u_N,u_N)\,\Big)\,<\,\infty
		\end{multline}
		because of \lose{}, \cref{eqn:rbound}, \cref{rem:polygon0} and \mucken{}.
		For $\varphi\in\cb(S)\times\cic(\R^d)$ we estimate 
		\begin{multline*}
			\Big|\int\limits_{S\times\R^d}\varphi\,(\lambda^*_m-u^*)\,\kappa\de\mu\Big|
			\,\leq\,\Big|\int\limits_{S\times\R^d}\varphi\,\lambda^{N_k}_{r_m}\,\kappa\de\mu_{N_k}-\int\limits_{S\times\R^d}\varphi\,\lambda^*_m\,\kappa\de\mu\,\Big|\\
			+\Big|\int\limits_{S\times\R^d}\varphi\,(G_1^{N_k}f-\lambda^{N_k}_{r_m})\,\kappa\de\mu_{N_k}\,\Big|
			+\Big|\int\limits_{S\times\R^d}\varphi\,(G_1^{N_k}f)\,\kappa\de\mu_{N_k}-\int\limits_{S\times\R^d}\varphi\,u^*\,\kappa\de\mu\,\Big|
		\end{multline*}
		and
		\begin{multline*}
			\Big|\E_\beta^{\infty,\kappa}(\varphi,\lambda_m^*-w_\kappa^*)\,\Big|
			\,\leq\,\Big|\,\E_\beta^{N_k,\kappa}(\varphi,\lambda_{r_m}^{N_k})-\E_1^{\infty,\kappa}(\varphi,\lambda_m^*)\,\Big|\\
			+\Big|\,\E_\beta^{N_k,\kappa}(\varphi,G_1^{N_k}f-\lambda_{r_m}^{N_k})\,\Big|
			+\Big|\,\E_\beta^{N_k,\kappa}(\varphi,G_1^{N_k}f)-\E_1^{\infty,\kappa}(\varphi,w_\kappa^*)\,\Big|\,.
		\end{multline*}
		In both estimates the second summand of the right hand side becomes arbitrarily small if $m$ is large enough, independent of $k$, due to 
		\cref{rem:polygon1,rem:polygon2} in combination with \mucken{}.
		So, we can first choose $m$ large enough, and then $k$, depending on $m$, to make also the first and third summand arbitrarily small, by  virtue of \cref{eqn:heki,eqn:hekii}.
		An $\varepsilon/3$ argument yields
		\begin{center}
			$\lambda_m^*\overset{m\to\infty}{\rightharpoonup}u^*$ weakly in $L^2(\kappa\mu)$\quad and\quad
			$\lambda_m^*\overset{m\to\infty}{\rightharpoonup}w_\kappa^*$ weakly in $(\dom_{\textnormal{min}}(\E^{\infty,\kappa}),\E^{\infty,\kappa}_\beta)$
		\end{center}
		in view of \cref{eqn:lambound}. Denoting the resolvent of $(\E^{\infty,\kappa},\dom_{\textnormal{min}}(\E^{\infty,\kappa}))$ by ${(R_\alpha)}_{\alpha>0}$
		the identity $u^*=w_\kappa^*$ in $\kappa\mu$ -a.e.~sense now follows from the equation
		\begin{multline*}
			\int\limits_E u^*\,v\,\kappa\de\mu\,=\,\lim_m\,\int\limits_E\lambda_m^*\,v\,\kappa\de\mu\,=
			\lim_m\,\E_\beta^{\infty,\kappa}(\lambda_m^*,R_\beta v)\\=\,\E_\beta^{\infty,\kappa}(w_\kappa^*,R_\beta v)
			\,=\,\int\limits_E w_\kappa^*\,v\,\kappa\de\mu
		\end{multline*}
		for $v\in L^2(\kappa\mu)$. Now, \cref{eqn:dedicate} is shown. 
		This concludes the proof.
	\end{proof}
	
	It is the last item listed in \cref{cond:mucken} through which the analysis of \cref{sec:tria} supports the proof of \cref{thm:m1}.
	We append a discussion about it here.
	Firstly, we ask about the role of $\kappa$. If \cref{eqn:liminf} holds for the choice $\kappa=\eins_E$, then we just pick $\mathscr U=\{\eins_E\}$
	and nothing needs to be proven concerning \cref{eqn:dommin}. 
	However, the option to consider different $\kappa$ provides the chance to make \cref{eqn:liminf} potentially weaker, hence easier to be verified.
	Such a procedure is legitimate as long as the family $\mathscr U$ is still large enough in a sense specified by \cref{eqn:dommin}.
	Beyond that, the next lemma gives a sufficient criterion under which \mucken{} still holds if the measure $\mu_N$ is perturbed by a weight
	function $g_N\in\mb(E)$ for $N\in\overline\N$. 
	The lemma addresses \mucken{} as an individual property, which a countable family of finite measures may have or may not have,
	and is not concerned with any other properties of that family, such as weak convergence, etc.
	\begin{lem}\label{lem:mopert}
		Let $0<c_1<c_2<\infty$ be constants and $g_N:S\times\R^d\to [c_1,c_2]$ for $N\in\overline\N$ be a function
		which meets at least one of the following three properties.
		\begin{enumerate}[wide, labelindent=0pt]
			\item For $s\in S$ the function $g_N(s,\cdot)$ is Lipschitz continuous on $\R^d$ with
			\begin{equation*}
				|\,g_N(s,x)-g_N(s,y)\,|\,\leq\,C_{\textnormal{Lip},N}(s)\,\sqrt{(x_1-y_1)^2+\dots+(x_d-y_d)^2}
			\end{equation*}
			for $x,y\in\R^d$, where the family $C_{\textnormal{Lip},N}(s)\in(0,\infty)$ meet
			\begin{equation*}
				\sup_N\,{\|\,C_{\textnormal{Lip},N}(\cdot)\,\|}_{L^2(\nu_N)}\,<\infty.
			\end{equation*}
			\item 
			$g_N(s,x_1,\dots,x_d)\leq g_N(s,y_1,\dots,y_d)$ for $s\in S$, $x,y\in\R^d$ with $x_1\leq y_1$, $\dots$, $x_d\leq y_d$.
			\item 
			$g_N(s,x_1,\dots,x_d)\geq g_N(s,y_1,\dots,y_d)$ for $s\in S$, $x,y\in\R^d$ with $x_1\leq y_1$, $\dots$, $x_d\leq y_d$.
		\end{enumerate}
		
		Then, the family $g_N\de\mu_N$, $N\in\overline\N$, meets \mucken{}.
	\end{lem}
	
	\begin{proof}
		Let the family $\mathscr U$ be the one which is suitable to verify that \mucken{} is met by ${(\mu_N)}_{N\in\overline\N}$
		and their disintegration measures $\nu_N$ and $\de m_s^N$ $=\varrho_N(s,\cdot)\de x$ with $s\in S$ for $N\in\N$ from \cref{eqn:desdef}. 
		As to \cref{eqn:dommin} there is nothing to show since the domain of the perturbed forms coincide with the unperturbed domains.
		We deal with the verification of \cref{eqn:liminf} in the following.
		Let $N\in\N$ be fixed. 
		
		The relevant densities in the perturbed case are given by
		\begin{equation*}
			\tilde g_N(s,x)\,:=\,\frac{g_N(s,x)}{w_N(s)}\qquad\text{with}\qquad
			w_N(s)\,:=\, \int\limits_{\R^d} g_N(s,y)\de m_s^N(y)
		\end{equation*}
		for $s\in S$ and $x\in\R^d$. Now \cref{eqn:desdef} holds if we replace $\mu_N$ by $g_N\mu_N$, $\nu_N$ by $w_N\nu_N$ and 
		$m_s^N$ by $\tilde g_N(s,\cdot)$ $m_s^N$ for $s\in S$. We observe that if $g_N$ satisfies either (i), (ii) or (iii) from the assumptions,
		then so does $\tilde g_N$ respectively. In the first part of this proof,
		we obtain a general estimate and it doesn't matter which of the three properties it is.

		Let $s\in S$, $r>0$, $\kappa\in\mathscr U$.
		We derive of an estimate for $\delta_r^{\kappa(s,\cdot)}(\,\tilde g_N(s,\cdot)\de m_s^N\,)$.
		To do so, we first use the characterization of \cref{rem:blt} and then apply the inequality  
		\begin{multline*}
			|\,(\tilde g_N\,\kappa\,\varrho_N)(s,x)-(\tilde g_N\,\kappa\,\varrho_N)(s,y)\,|\,
			\\\leq\,\tilde g_N(s,x)\,|\,(\kappa\,\varrho_N)(s,x)-(\kappa\,\varrho_N)(s,y)\,|
			+(\kappa\,\varrho_N)(s,y)\,|\,\tilde g_N(s,x)-\tilde g_N(s,y)\,|
		\end{multline*}
		for $x,y\in\R^d$ together with $c_1/c_2\leq \tilde g_N(s,\cdot)\leq c_2/c_1$.
		The supremum is taken over all primal functions $\varphi$, $\eta$ from $\mathscr C$ and the sum runs over $\alpha\in r\,\Z^d$ in the next lines.
		\begin{align*}
			&{\delta_r^{\kappa(s,\cdot)}(\,\tilde g_N(s,\cdot)\de m_s^N\,)}^2\\
			&=\,\sup_{\varphi,\eta}\int\limits_{\R^d}\Big|\,\sum_{\alpha}\int\limits_{\R^d}
			\frac{|\,(\tilde g_N\,\kappa\,\varrho_N)\,(s,y)-(\tilde g_N\,\kappa\,\varrho_N)\,(s,x)\,|}{r^d\,\varrho_N(s,y)}
			\eta_r^\alpha(x)\de x\,\varphi_r^\alpha(y)\,\Big|^2\frac{\varrho_N(s,y)}{\tilde g_N(s,y)}\de y\\
			&\leq\,
				\frac{2\,c_2^3}{c_1^3}\,\sup_{\varphi,\eta}\int\limits_{\R^d}\Big|\int\limits_{\R^d}
				\sum_{\alpha}\,\frac{|\,(\kappa\,\varrho_N)\,(s,y)-(\kappa\,\varrho_N)\,(s,x)\,|}{r^d\,\varrho_N(s,y)}\,\eta_r^\alpha(x)\,\varphi_r^\alpha(y)\de x\,\Big|^2\varrho_N(s,y)\de y\\
				&+\frac{2\,c_2}{c_1}\,\sup_{\varphi,\eta}\int\limits_{\R^d}\Big|\,r^{-d}\,\sum_{\alpha}\int\limits_{\R^d}
				\kappa(s,y)\,|\,\tilde g_N(s,y)-\tilde g_N(s,x)\,|\,
				\eta_r^\alpha(x)\de x\,\varphi_r^\alpha(y)\,\Big|^2\varrho_N(s,y)\de y\\
			&\leq\frac{2\,c_2^3}{c_1^3}\,{\delta_r^{\kappa(s,\cdot)}(m_s^N)}^2\\
				&+\,\frac{2\cdot 9^dc_2}{c_1}\,\sup_{\varphi,\eta}\,r^{-d}\sum_{\alpha}\,\int\limits_{\R^d\times\R^d}
				\kappa(s,y)^2\,|\,\tilde g_N(s,y)-\tilde g_N(s,x)\,|^2\,
				\eta_r^\alpha(x)\de x\,\varphi_r^\alpha(y)\,\varrho_N(y)\de y.
		\end{align*}
		For the estimate leading to the last term we used that for $\varphi\in\mathscr C$ and $\alpha,\beta\in r\,\Z^d$ it holds $\varphi_r^\alpha\,\varphi_r^\beta\equiv 0$
		unless $\beta$ is contained in the set $\alpha+[-4r,4r]^d$. So, for any family $h_\alpha:\R^d\to[0,\infty)^d$, $\alpha\in r\,\Z^d$, we have 
		\begin{align*}
			\Big|\,\sum_{\alpha\in r\Z^d}h_\alpha(y)\,\varphi_r^\alpha(y)\,\Big|^2
			\,&=\,\sum_{\alpha,\beta\in r\Z^d}h_\alpha(y)\,\varphi_r^\alpha(y)\,h_\beta(y)\,\varphi_r^\beta(y)\\
			&\leq\sum_{\substack{\alpha,\beta\in r\Z^d\\\beta\in\alpha+ [-4r,4r]^d}}\Big(\,
			\frac{1}{2}\,\big(h_\alpha(y)\,\varphi_r^\alpha(y)\big)^2
			+\frac{1}{2}\,\big(h_\beta(y)\,\varphi_r^\beta(y)\big)^2\,\Big)\\
			&\leq\,9^d\,\sum_{\alpha\in r\Z^d}h_\alpha(y)^2\,\varphi_r^\alpha(y).
		\end{align*}
		at any point $y\in\R^d$. So, for the estimate in question, we just have to choose
		\begin{equation*}
			h_\alpha(y)\,:=\,r^{-d}\int\limits_{\R^d}\kappa(s,y)\,|\,\tilde g_N(s,y)-\tilde g_N(s,x)\,|\,\eta_r^\alpha(x)\de x,\qquad\alpha\in r\,\Z^d,\,y\in\R^d,
		\end{equation*}
		for given $\eta\in\mathscr C$ and then apply Jensen's inequality and Fubini's theorem.
		
		We now fix $\varphi,\eta\in\mathscr C$ and tackle the term
		\begin{equation*}
			(\star)\,:=\,r^{-d}\sum_{\alpha\in r\Z^d}\,\int\limits_{\R^d}\int\limits_{\R^d}\kappa(s,y)^2\,|\,\tilde g_N(s,y)-\tilde g_N(s,x)\,|^2\,
			\eta_r^\alpha(x)\de x\,\varphi_r^\alpha(y)\,\varrho_N(s,y)\de y,
		\end{equation*}
		which appears in the estimate for $\delta_r^{\kappa(s,\cdot)}(\,\tilde g_N(s,\cdot)\de m_s^N\,)$ above.

		First, we look at the significantly easier case where  property (i) of the assumptions of this lemma is satisfied. 
		Since both, $\varphi_r^\alpha$ and $\eta_r^\alpha$, are supported on $\alpha+[-2r,2r]^d$ and the Lipschitz constant of
		$\tilde g_N(s,\cdot)$ is smaller equal $C_{\textnormal{Lip},N}(s)/c_1$, we have
		\begin{align*}
			(\star)\,&\leq\, \big(\,4\,r\,\sqrt d\,C_{\textnormal{Lip},N}(s)/c_1\,\big)^2\,r^{-d}\sum_{\alpha\in r\Z^d}\,\int\limits_{\R^d}\int\limits_{\R^d}
			\kappa(s,y)^2\,\eta_r^\alpha(x)\de x\,\varphi_r^\alpha(y)\,\varrho_N(s,y)\de y\\
			&=\,\frac{4^2\,r^2\,d}{c_1^2}\,C_{\textnormal{Lip},N}(s)^2\,\sum_{\alpha\in r\Z^d}\,\int\limits_{\R^d}\kappa(s,y)^2\,\varphi_r^\alpha(y)\,\varrho_N(s,y)\de y\\
			&=\,\frac{4^2\,r^2\,d}{c_1^2}\,C_{\textnormal{Lip},N}(s)^2\,\int\limits_{\R^d}\kappa(s,y)^2\,\varrho_N(s,y)\de y
			\,\leq\,\frac{4^2\,r^2\,d}{c_1^2}\,C_{\textnormal{Lip},N}(s)^2.
		\end{align*}
		If we plug in this estimate for $(\star)$ into the initial bound for $\delta_r^{\kappa(s,\cdot)}(\,\tilde g_N(s,\cdot)\de m_s^N\,)$ above
		and use the triangular inequality for the norm of $L^2(w_N\nu_N)$, then we arrive at
		\begin{multline}\label{eqn:delesti}
			{\|\,\delta_r^{\kappa(s,\cdot)}(\,\tilde g_N(s,\cdot)\de m_s^N\,)\,\|}_{L^2(w_N\nu_N)}\\
			\leq\,c_2^2\,{\Big(\,\frac{2}{ c_1^3}\,\Big)}^{\frac{1}{2}}\,{\|\,\delta_r^{\kappa(s,\cdot)}(m_s^N)\,\|}_{L^2(\nu_N)}
			+4\,r\,c_2\,{\Big(\,2\,9^d\,\frac{d}{c_1^3}\,\Big)}^{\frac{1}{2}}\,{\|\,C_{\textnormal{Lip},N}(s)\,\|}_{L^2(\nu_N)}.
		\end{multline}

		We address the case, where either property (ii) or (iii) of the assumptions of this lemma is satisfied, in which a similar bound as in \cref{eqn:delesti} can be obtained.
		We claim that in this case 
		\begin{equation}\label{eqn:starestii}
			(\star)\,\leq\,4\,\frac{c_2^2}{c_1^2}\,\delta_{2r}^{\kappa(s,\cdot)}(m_s^N).
		\end{equation}
		If true, plugging in the estimate for $(\star)$ into the initial bound for $\delta_r^{\kappa(s,\cdot)}(\,\tilde g_N(s,\cdot)\de m_s^N\,)$ above,
		using the triangular inequality for the norm of $L^2(w_N\nu_N)$ and then Jensen's inequality, we arrive at
		\begin{multline}\label{eqn:delestii}
			{\|\,\delta_r^{\kappa(s,\cdot)}(\,\tilde g_N(s,\cdot)\de m_s^N\,)\,\|}_{L^2(w_N\nu_N)}\\
			\leq\,c_2^2\,{\Big(\,\frac{2}{ c_1^3}\,\Big)}^{\frac{1}{2}}\,{\|\,\delta_r^{\kappa(s,\cdot)}(m_s^N)\,\|}_{L^2(\nu_N)}
			+2\,c_2^2\,{\Big(\,9^d\,\frac{2}{c_1^3}\,{\|\,\delta_{2r}^{\kappa(s,\cdot)}(\varrho_N\de x)\,\|}_{L^2(\nu_N)}\,\Big)}^{\frac{1}{2}}.
		\end{multline}
		We only write down the proof of \cref{eqn:starestii} in the case, where property (ii) is satisfied, since the case of (iii) works analogous. The point $s\in S$ is fixed. 
		We set $\rho$ $:=\varrho_N(s,\cdot)$, $\tau$ $:=\kappa(s,\cdot)^2$ and $f$ $:=\tilde g(s,\cdot)$ for short.
		In the following lines we first exploit the monotonicity of $f$, 
		then shift the index $\alpha$ of the sum, before we use linearity of the integral and translation invariance of the Lebesgue measure.
		Again we observe that both, $\varphi_r^\alpha$ and $\eta_r^\alpha$, are supported on $\alpha+[-2r,2r]^d$ and \hypertarget{eqn:atarestiia}{calculate} for $(\star)$
		\begin{multline*}
			\,r^{-d}\sum_{\alpha\in r\Z^d}\,\int\limits_{\R^d}\int\limits_{\R^d}|\,f(x)-f(y)\,|^2\,\tau(y)\,\eta_r^\alpha(x)\de x\,\varphi_r^\alpha(y)\,\rho(y)\de y\\
			\begin{aligned}
				&\leq\,\sum_{\alpha\in r\Z^d}\,\int\limits_{\R^d}\big(f(\alpha+2\,r\,\univ)-f(\alpha-2\,r\,\univ)\big)^2\varphi_r^\alpha(y)\,\tau(y)\,\rho(y)\,\de \\
				&\leq\,2\,\frac{c_2}{c_1}\,\sum_{\alpha\in r\Z^d}\,\int\limits_{\R^d}\big(f(\alpha+2\,r\,\univ)-f(\alpha-2\,r\,\univ)\big)\varphi_r^\alpha(y)\,\tau(y)\,\rho(y)\,\de y\\
				&=\,2\,\frac{c_2}{c_1}\,\sum_{\alpha\in r\Z^d}f(\alpha)\int\limits_{\R^d}\big(\,\varphi_r^{\alpha-2r\univ}(y)-\varphi_r^{\alpha+2r\univ}(y)\,\big)\,\tau(y)\,\rho(y)\,\de y\\
				&=\,2\,\frac{c_2}{c_1}\,\sum_{\alpha\in r\Z^d}f(\alpha)\int\limits_{\R^d}\varphi_r^{\alpha}(y)\,\big(\,(\tau\,\varrho)\,(y-2\,r\,\univ)
				-(\tau\,\rho)\,(y+2\,r\,\univ)\,\big)\de y\\
				&\leq\,2\,\frac{c_2^2}{c_1^2}\,\sum_{\alpha\in r\Z^d}\,\int\limits_{\R^d}\varphi_r^{\alpha}(y)\,|\,(\tau\,\varrho)\,(y-2\,r\,\univ)
				-(\tau\,\rho)\,(y+2\,r\,\univ)\,|\de y
			\end{aligned}\\
				=\,2\,\frac{c_2^2}{c_1^2}\,\int\limits_{\R^d}|\,(\tau\,\rho)\,(y-2\,r\,\univ)\,
				-(\tau\,\rho)\,(y+2\,r\,\univ)\,|\de y.\numberthis
		\end{multline*}
		\newcounter{eqnAtarestiiA}
		\setcounter{eqnAtarestiiA}{\value{equation}}
		\newcommand{\eqnAtarestiiA}{Equations \hyperlink{eqn:atarestiia}{(\thec.\theeqnAtarestiiA)}}
		Note that neither the function $f$, nor the primal functions $\varphi$ or $\eta$, appear in the latter expression.
		To go on with the estimate for $(\star)$, we define $\hat\varphi$ $:=\eins_{[0,1)^d}$.
		Recalling the perturbation operator from \cref{eqn:phetaper} we now split
		\begin{multline}
			\int\limits_{\R^d}|\,(\tau\,\rho)\,(y-2\,r\,\univ)\,
			-(\tau\,\rho)\,(y+2\,r\,\univ)\,|\de y
			\,\leq\,\int\limits_{\R^d}\Big|(\tau\,\rho)\,(y-2\,r\,\univ)-I^{\hat\varphi,\hat\varphi}_{2r}(\tau\,\rho)(y)\,\Big|\de y
			\\
			+\int\limits_{\R^d}\Big|\,I^{\hat\varphi,\hat\varphi}_{2r}(\tau\,\rho)(y)-(\tau\,\varrho)\,(y+2\,r\,\univ)\,\Big|\de y.
		\end{multline}
		We can show that each of the two summands is bounded from above by $\delta_r^\tau(\rho\de x)$.
		Since the argumentation is analogous for the two summands, we restrict ourselves to put it here for only one of them.
		In the following lines let $\hat\eta$ $:=\eins_{[-1,0)^d}$.
		First we use the translation invariance of the Lebesgue measure,
		then a straight-forward calculation using the elementary properties of primal functions \hypertarget{eqn:atarestiic}{yields}
		\begin{align}
			&\int\limits_{\R^d}\Big|(\tau\,\rho)\,(y-2\,r\,\univ)-I^{\hat\varphi,\hat\varphi}_{2r}(\tau\,\rho)(y)\,\Big|\de y\nonumber\\
			&\qquad=\,\int\limits_{\R^d}\Big|(\tau\,\varrho)\,(y)-
			\sum_{\alpha\in {(2r)}\Z^d}(2r)^{-d}\int\limits_{\R^d}\hat\varphi_{2r}^\alpha(x)\,\tau(x)\,\rho(x)\de x\,\hat\varphi_{2r}^{\alpha-2r\univ}(y)\Big|\de y\nonumber\\
			&\qquad=\,\int\limits_{\R^d}\Big|(\tau\,\varrho)\,(y)-
			\sum_{\alpha\in {(2r)}\Z^d}(2r)^{-d}\int\limits_{\R^d}\hat\varphi_{2r}^\alpha(x)\,\tau(x)\,\rho(x)\de x\,\hat\eta_{2r}^{\alpha}(y)\Big|\de y\nonumber\\
			&\qquad\leq\,\int\limits_{\R^d}\sum_{\alpha\in {(2r)}\Z^d}(2r)^{-d}\int\limits_{\R^d}|\,(\tau\,\rho)\,(y)-(\tau\,\rho)\,(x)\,|\,
			\hat\varphi_{2r}^\alpha(x)\de x\,\hat\eta_{2r}^{\alpha}(y)\de y\nonumber\\
			&\qquad=\,\int\limits_{\R^d}\eins_{\R^d}(y)\,R^{\hat\eta,\hat\varphi}_{2r}(\tau\,\rho)(y)\de y\,\leq \,\delta_{2r}^\tau(\rho\de x).
		\end{align}
		\eqnAtarestiiA
		\refstepcounter{eqnAtarestiiA}
		\refstepcounter{eqnAtarestiiA}to \hyperlink{eqn:atarestiic}{(\thec.\theeqnAtarestiiA)}
		provide the proof of \cref{eqn:starestii} and we go on to make the final remarks which are necessary to finish the proof of this lemma.
		
		We observe that $c_1/c_2\leq \tilde g_N(s,\cdot)\leq c_2/c_1$ implies
		\begin{equation*}
			C_r^{\kappa(s,\cdot)}(\,\tilde g_N(s,\cdot) \de m_s^N\,)\,\leq\,\frac{c_2^2}{c_1^2}\,C_r^{\kappa(s,\cdot)}(m_s^N)
		\end{equation*}
		for $s\in S$.
		Finally, due to \cref{eqn:delesti}, respectively \cref{eqn:delestii}, 
		the family $\mathscr U$ is suitable to provide \mucken{} for the sequence ${(g_N\de\mu_N)}_N$. 	
		This concludes the proof.
	\end{proof}

	\section{Application to infinite-dimensional problems and a first example}\label{sec:appli}
	
	\subsection{Mosco convergence of standard gradient forms on Fr{\'e}chet spaces}\label{sec:stagra}
	
	This section deals with gradient type Dirichlet forms on a locally convex, real topological vector space $E$,
	which is also assumed to be a Polish space. Hence, $E$ is a separable Fr{\'e}chet space. Its topological dual space is denoted by $E'$.
	We define the linear space 
	\begin{equation*}
		\fcb\,:=\,\Big\{\,f\circ\,(l_1,\cdots,l_m)\,\Big|\,m\in\N,\,f\in\cb^\infty(\R^m),\,l_1,\dots,l_m\in E'\,\Big\}
	\end{equation*}
	of cylindrical smooth functions on $E$.
	Let $(H,|\cdot|,\langle\cdot,\cdot\rangle)\subset E$ be a Hilbert space which is densely embedded in $E$. The gradient $\nabla F(z)$ 
	of a cylindrical smooth function $F:E\to\R$ at a point $z\in E$ denotes the unique element in $H$
	which (via the Riesz isomorphism) represents the linear functional
	\begin{equation*}
		H\ni h\mapsto \frac{\partial F}{\partial h}(z)\,:=\,\frac{\de F(z+t\,h)}{\de t}\Big|_{t=0}.
	\end{equation*}
	The right hand side $\frac{\partial F}{\partial h}(z)$ is the G{\^a}teaux derivative of $u$ in the direction $h$ at $z$. 
	Identifying $H$ with its dual via the Riesz isomorphism we get
	\begin{equation*}
		E'\,\subset\, H'\,=\, H\,\subset\, E.
	\end{equation*}
	For $F=f\circ\,(l_1,\cdots,l_m)$ with $f\in\cb^\infty(\R^m)$, $l_1,\dots,l_m\in E'$, $m\in\N$, the directional derivative 
	at a point $z\in E$ in a direction $h\in H$ then reads
	\begin{equation*}
		\langle\nabla F(z),h\rangle\,=\,\sum_{i=1}^m\partial_if(l_1(z),\dots ,l_m(z))\,\langle l_i,h\rangle.
	\end{equation*}
	The norm of the gradient can be estimated from above by
	\begin{align*}
		|\nabla F(z)|\,\leq\,\sup_{\substack{h\in H\\|h|=1}}\,\langle\nabla F(z),h\rangle\,
		\leq\,\sum_{i=1}^m\partial_if(l_1(z),\dots ,l_m(z))\,|l_i|\,\leq\,\sup_{1\leq i\leq m}\,\|\partial_i f\|_\infty\,\sum_{i=1}^m|l_i|.
	\end{align*}
	For a subset $A\subset E'$ we specify a linear subspace of $\fcb$ by writing
	\begin{equation*}
		\fcb(A)\,:=\,\Big\{\,f\circ\,(l_1,\cdots,l_m)\,\Big|\,m\in\N,\,f\in\cb^\infty(\R^m),\,l_1,\dots,l_m\in A\,\Big\}.
	\end{equation*}
	Let ${(\mu_N)}_{N\in\N}$ be a sequence of weakly converging probability measure on $E$ with limit $\mu_\infty$.
	The minimal gradient form on $E$ is a Dirichlet form $(\E^N,\dom(\E^N))$ on $L^2(E,\mu_N)$ for given $N\in\overline\N$.
	It arises from taking the closure in $L^2(E,\mu_N)$ of the form 
	\begin{equation}\label{eqn:mingrad}
		\E^N(u,v)\,=\,\int\limits_E\langle \nabla \tilde u,\nabla \tilde v\rangle\de\mu_N,\qquad u,v\in\dom_\textnormal{pre}(\E^N),
	\end{equation}
	with pre-domain 
	\begin{equation*}
		\dom_\textnormal{pre}(\E^N)\,:=\,\Big\{\,u\in L^2(E,\mu_N)\,\Big|\, u(\cdot)=\tilde u(\cdot)\,\mu_N\textnormal{ -a.e. for some }\tilde u\in\fcb\,\Big\},
	\end{equation*}
	always assuming this procedure and the assignment of \cref{eqn:mingrad} is well-defined.
	The gradient forms, as defined here, are extensively studied in \cite{albe, rocla}.
	We find a criterion under which the minimal gradient forms on $E$ with varying reference measure $\mu_N$ converge as $N\to\infty$.
	The well-definedness and closability of the respective forms is implied, alongside the Mosco convergence of their closures, 
	by the conditions listed in \cref{cor:mosco}. These focus on certain `component forms'. 
	We look at two different ways how the component forms can be defined.
	In one case an orthonormal basis $\eta_1,\eta_2,\dots$ of $H$ is selected.
	In the other case a set $K_d:=\{\xi_i^d|i=1,\dots,d\}\subset H$ of linearly independent vectors is chosen for each $d\in\N$.
	In the second case, we further assume 
	\begin{equation}\label{eqn:piconv}
		\sup_{d\in\N}\,\sum_{i=1}^d\langle h,\xi^d_i\rangle^2\,=\,|h|^2,\qquad h\in H.
	\end{equation}
	For $d\in\N$ we fix linear spaces $S_{\eta_d}\subset E$ and $S_{K_d}\subset E$, which are closed complementing subspaces of span$(\{\eta_d\})$, respectively of span$(K_d)$, in $E$.
	In other words, we decompose $E$ into the direct sum
	\begin{equation}\label{eqn:iso1}
		E\,=\,S_{\eta_d}\oplus\textnormal{span}(\{\eta_d\})\,\simeq\,S_{\eta_d}\oplus\R,
	\end{equation}
	respectively
	\begin{equation}\label{eqn:iso2}
		E\,=\,S_{K_d}\oplus\textnormal{span}(K_d)\,\simeq\,S_{K_d}\oplus\R^d
	\end{equation}
	for $d\in\N$. Let $\pi_{\eta_d}:E\to\R$ and $\pi_{K_d}:E\to\R^d$
	denote the second components of the isomorphisms behind \cref{eqn:iso1}, respectively \cref{eqn:iso2}.
	W.l.o.g., we may assume that $\pi_{\eta_d}:E\to\R$, $\pi_{K_d}:E\to\R^d$ are surjective linear mappings such that
	$\pi_{\eta_d}\eta_d$ equals $1\in\R$ and $\pi_{K_d}\xi_i^d$ equals $\univ_i$ (the $i$-th unit vector of a Euclidean space), while
	$S_{\eta_d}=\textnormal{Ker}(\pi_{\eta_d})$ and $S_{K_d}= \textnormal{Ker}(\pi_{K_d})$.
	Hence, we consider
	\begin{equation}\label{eqn:jeta}
		J_{\eta_d}:E\ni h\,\mapsto\,(\,h-(\pi_{\eta_d}h)\,\eta_d\,,\,\pi_{\eta_d} h\,)\,\in\,S_{\eta_d}\times\R,
	\end{equation}
	respectively
	\begin{equation}\label{eqn:xixi}
		J_{K_d}:E\ni h\,\mapsto\,\Big(\,h-\sum_{i=1}^d(\univ_i^\textnormal{T}\,\pi_{K_d}h)\,\xi^d_i\,,\,\pi_{K_d} h\,\Big)\,\in\,S_{K_d}\times\R^d.
	\end{equation}
	Clearly, $J_{\eta_d}^{-1}$ and $J_{K_d}^{-1}$ are continuous and so are $J_{\eta_d}$ and $J_{K_d}$ by the open mapping theorem.

	For every $d\in\N$ the criteria of \cref{cond:hamza,cond:mucken} are now imposed on the family, indexed by $N\in\overline\N$, which is obtained by taking the image of $\mu_N$
	under the maps of \cref{eqn:jeta}, respectively \cref{eqn:xixi}. 
	We recall the family of Dirichlet forms, indexed by $N\in\overline\N$, constructed in \cref{sec:lowsii}, subsequent to \cref{cond:mucken}. 
	The starting point in \cref{sec:lowsii} has been a sequence of weakly convergent probability measures on a product of an abstract Polish space $S$ and a 
	finite-dimensional Euclidean space.
	For each $d\in\N$ we now look at the family ${(\mu_N\circ J_{K_d}^{-1})}_{N\in\overline\N}$ with state space $S_{K_d}\times \R^d$
	and denote the corresponding family of forms, defined as in \cref{sec:lowsii}, by $(\bar\E^{N,K_d},\dom(\bar\E^{N,K_d}))$, $N\in\overline\N$.
	Accordingly, $(\bar\E^{N,K_d},\dom(\bar\E^{N,K_d}))$ is then a Dirichlet form on $L^2(S_{K_d}\times\R^d,\mu_N\circ J_{K_d}^{-1})$ for $N\in\overline\N$.
	Next, we consider the image forms under the inverse of the map in \cref{eqn:jeta}.
	For $d\in\N$, $N\in\overline\N$ and $u,v\in\dom(\E^{N,K_d}):=$ $\{\,w\in L^2(E,\mu_N)\,|\,w\circ J_{K_d}^{-1}\in\dom(\bar\E^{N,K_d})\,\}$ we define
	\begin{align*}
		\E^{N,K_d}(u,v)\,&:=\,\bar\E^{N,K_d}(u\circ J_{K_d}^{-1},v\circ J_{K_d}^{-1})\\
		&=\,\sum_{i=1}^d\,\int\limits_{S_{K_d}\times \R^d} \partial_{x_i}(u\circ J_{K_d}^{-1})\,\partial_{x_i}(v\circ J_{K_d}^{-1})\de (\mu_N\circ J_{K_d}^{-1})
	\end{align*}
	(confer with \cref{eqn:conf}). In the same way $(\bar\E^{N,\eta_d},\dom(\bar\E^{N,\eta_d}))$ 
	shall be defined as a  Dirichlet form on $L^2(S_{\eta_d}\times\R,\mu_N\circ J_{\eta_d}^{-1})$ for $N\in\overline\N$ and $d\in\N$. Again, for
	$u,v\in\dom(\E^{N,\eta_d}):=$ $\{\,w\in L^2(E,\mu_N)\,|\,w\circ J_{\eta_d}^{-1}\in\dom(\bar\E^{N,\eta_d})\,\}$ we set
	\begin{align*}
		\E^{N,\eta_d}(u,v)\,&:=\,\bar\E^{N,\eta_d}(u\circ J_{\eta_d}^{-1},v\circ J_{\eta_d}^{-1})\\
		&=\,\int\limits_{S_{\eta_d}\times \R} \partial_x(u\circ J_{\eta_d}^{-1})\,\partial_x(v\circ J_{\eta_d}^{-1})\de (\mu_N\circ J_{\eta_d}^{-1}).
	\end{align*}
	Furthermore, we define the Dirichlet forms $\sup_d\E^{N,K_d}$ and $\sum_i\E^{N,\eta_i}$ on $L^2(E,\mu_N)$ for $N\in\overline\N$. Their domains read
	\begin{align*}
		\dom(\,{\textstyle\sup_d}\E^{N,K_d}\,)\,:=\,\Big\{u\in\bigcap_{d\in\N}\dom(\E^{N,K_d})\, \Big|\,\sup_{d\in\N}\E^{N,K_d}(u,u)<\infty\,\Big\},
	\end{align*}
	respectively
	\begin{equation*}
		\dom(\,{\textstyle\sum_i}\E^{N,\eta_i}\,)\,:=\,\Big\{u\in\bigcap_{i\in\N}\dom(\E^{N,\eta_i})\, \Big|\,\sum_{i=1}^\infty\E^{N,\eta_i}(u,u)<\infty\,\Big\}.
	\end{equation*}
	
	\begin{rem}
		Let $N\in\overline\N$.
		\begin{thmlist}[wide, labelindent=0pt]
			\item \label{rem:m1} We assume that the family $\mu_N\circ J_{\eta_i}^{-1}$, $N\in\overline \N$, satisfies \cref{cond:hamza,cond:mucken} for $i\in\N$.
			Since $J_{\eta_i}^{-1}(s,x)=s+x\,\eta_i$ for $s\in S_{\eta_i}$, $x\in\R$ and $i\in\N$, the G{\^a}teaux derivative of $\tilde u\in\fcb$ at
			the point $J_{\eta_i}^{-1}(s,x)$ in the direction $\eta_i$ calculates as
			\begin{equation*}
				\frac{\partial\tilde u}{\partial \eta_i}(J_{\eta_i}^{-1}(s,x))\,=\,
				\frac{\de \tilde u(s+x\,\eta_i+t\,\eta_i)}{\de t}\Big|_{t=0}\,=\,\partial_x(\tilde u\circ J_{\eta_i}^{-1})(s,x).
			\end{equation*}
			Hence,
			\begin{equation*}
				|\nabla \tilde u|^2(z)\,=\,\sum_{i=1}^\infty\,\langle \eta_i,\nabla\tilde u(z)\rangle^2\,
				=\,\sum_{i=1}^\infty|\,\partial_x(\tilde u\circ J_{\eta_i}^{-1})\,(J_{\eta_i}z)\,|^2
			\end{equation*}
			for $z\in E$. Clearly, $\dom_\textnormal{pre}(\E^N)\subset\dom(\,{\textstyle\sum_i}\E^{N,\eta_i}\,)$ and moreover
			\begin{equation*}
				\int\limits_{E}|\nabla \tilde u|^2\de\mu_N\,=\,\sup_{d\in\N}\,\sum_{i=1}^d\,\int\limits_{E}|\,\partial_x(\tilde u\circ J_{\eta_i}^{-1})\,(J_{\eta_i}z)\,|^2\de\mu_N(z)
				\,=\,\sum_{i=1}^\infty\E^{N,\eta_i}(u,u)
			\end{equation*}
			for $u\in\dom_\textnormal{pre}(\E^N)$ (with representative $\tilde u\in\fcb$). 
			This, in turn, implies that the form in \cref{eqn:mingrad} is well-defined and closable on $L^2(E,\mu_N)$ (with closure $(\E^N,\dom(\E^N))$) and that
			$({\textstyle\sum_i}\E^{N,\eta_i},\dom(\,{\textstyle\sum_i}\E^{N,\eta_i}\,))$ is an extension of $(\E^N,\dom(\E^N))$, i.e.~
			$\dom(\E^N)\subset\dom(\,{\textstyle\sum_i}\E^{N,\eta_i}\,)$ and
			\begin{equation*}
				\E^N(u,u)\,=\,\sum_{i=1}^\infty\E^{N,\eta_i}(u,u)
			\end{equation*}
			for $u\in\dom(\E^N)$.
			
			\item  The other case behaves analogously. 
			We assume that the family $\mu_N\circ J_{K_d}^{-1}$, $N\in\overline \N$, satisfies \cref{cond:hamza,cond:mucken} for $d\in\N$.
			We have $J_{K_d}^{-1}(s,x)=s+x_1\,\xi_1^d+\dots+x_d\,\xi_d^d$ for $s\in S_{K_d}$, $x\in\R^d$ and $d\in\N$. 
			If $\tilde u\in\fcb$ and $1\leq i\leq d$, then
			\begin{equation*}
				\frac{\partial\tilde u}{\partial \xi_i^d}(J_{K_d}^{-1}(s,x))\,=\,
				\frac{\de \tilde u(s+x_1\,\xi_1^d+\dots+x_d\,\xi_d^d+t\,\xi_i^d)}{\de t}\Big|_{t=0}\,=\,\partial_{x_i}(\tilde u\circ J_{K_d}^{-1})(s,x).
			\end{equation*}
			Using \cref{eqn:piconv} we conclude
			\begin{equation*}
				|\nabla \tilde u|^2(z)\,=\,\sup_{d\in\N}\,\sum_{i=1}^d\,\langle \xi_i^d,\nabla\tilde u(z)\rangle^2\,
				=\,\sup_{d\in\N}\,\sum_{i=1}^d|\,\partial_{x_i}(\tilde u\circ J_{K_d}^{-1})\,(J_{K_d}z)\,|^2
			\end{equation*}
			for $z\in E$. Clearly, $\dom_\textnormal{pre}(\E^N)\subset \dom(\,{\textstyle\sup_d}\E^{N,K_d}\,)$ and moreover
			\begin{equation*}
				\int\limits_{E}|\nabla \tilde u|^2\de\mu_N\,=\,\sup_{d\in\N}\,\sum_{i=1}^d\,\int\limits_{E}|\,\partial_{x_i}(\tilde u\circ J_{K_d}^{-1})\,(J_{K_d}z)\,|^2\de\mu_N(z)
				\,=\,\sup_{d\in\N}\,\E^{N,K_d}(u,u)
			\end{equation*}
			for $u\in\dom_\textnormal{pre}(\E^N)$ (with representative $\tilde u\in\fcb$). 
			This, in turn, implies that the form in \cref{eqn:mingrad} is well-defined and closable on $L^2(E,\mu_N)$ (with closure $(\E^N,\dom(\E^N))$) and that
			$({\textstyle\sup_d}\E^{N,K_d},\dom(\,{\textstyle\sup_d}\E^{N,K_d}\,))$ is an extension of $(\E^N,\dom(\E^N))$.
		\end{thmlist}
	\end{rem}
	
	We now assume $\supp[\mu_N]\subset\supp[\mu_\infty]$ for $N\in\N$, as in \cref{sec:mokushi}, and understand ${(L^2(E,\mu_N))}_{N}$ as a sequence of converging Hilbert spaces
	with asymptotic space $L^2(E,\mu_\infty)$.
	
	\begin{thm}\label{cor:mosco}
		The sequence ${(\E^N)}_{N\in\N}$ converges to $\E^\infty$ in the sense of Mosco if one of the following two conditions is fulfilled:
		\begin{thmlist}[wide, labelindent=0pt]
			\item \label{cor:mosci} The family $\mu_N\circ J_{\eta_i}^{-1}$, $N\in\overline \N$, satisfy \cref{cond:hamza,cond:mucken} for $i\in\N$ and
			\begin{equation*}
				\dom(\textstyle{\sum_i\E^{\infty,\eta_i}}\big)=\dom(\E^\infty).
			\end{equation*}
			\item We assume \cref{eqn:piconv} and 
			\begin{equation*}
				\dom(\textstyle{\sup_{d\in\N}}\E^{\infty,K_d})=\dom(\E^\infty).
			\end{equation*}
			Moreover, the family $\mu_N\circ J_{K_d}^{-1}$, $N\in\overline \N$, satisfy \cref{cond:hamza,cond:mucken} for $d\in\N$.
		\end{thmlist}
	\end{thm}

	\begin{proof}
		The proof of (i) and (ii) work analogously.
		We only write down the proof of (i) here. This is accomplished by verifying both conditions of \cref{eqn:m1m2m}. 
		
		We start with Property (a). 
		Let ${(u_N)}_{N\in\overline\N}\in \prod_{N\in\overline\N}L^2(E,\mu_N)$ with $u_N\overset{\textnormal{w.}}{\underset{N}{\rightharpoonup}}u_\infty$.
		Then $u_N\circ J_{\eta_i}^{-1}\overset{\textnormal{w.}}{\underset{N}{\rightharpoonup}}u_\infty\circ J_{\eta_i}^{-1}$
		(in the sense of $\prod_{N\in\overline\N}L^2(E,\mu_N\circ J_{\eta_i}^{-1})$) for $i\in\N$, since $J_{\eta_i}$ is a topological homeomorphism.
		Let $d\in\N$. In the following estimate we make a multiple use of \cref{thm:m1}, apply Fatou's lemma and then \cref{rem:m1}.
		We have $u_\infty\circ J_{\eta_i}^{-1}\in\dom(\bar\E^{\infty,\eta_i})$ for $i\in\{1,\dots,d\}$ and
		\begin{multline}\label{eqn:propa}
			\sum_{i=1}^d\bar\E^{\infty,\eta_i}(u_\infty\circ J_{\eta_i}^{-1},u_\infty\circ J_{\eta_i}^{-1})\,
			\leq\,\sum_{i=1}^d\,\liminf_{N\to\infty}\,\bar\E^{N,\eta_i}(u_N\circ J_{\eta_i}^{-1},u_N\circ J_{\eta_i}^{-1})\\
			=\,\sum_{i=1}^d\,\liminf_{N\to\infty}\,\E^{N,\eta_i}(u_N,u_N)
			\,\leq\,\liminf_{N\to\infty}\,\sum_{i=1}^\infty \E^{N,\eta_i}(u_N,u_N)\,=\,\liminf_{N\to\infty}\,\E^N(u_N,u_N)
		\end{multline}
		under the condition that $u_N\in\dom(\E^N)$ for infinitely many $N$ and that the right hand side of \cref{eqn:propa} is finite.
		Now Property (a) follows from the assumption, since the choice of $d\in \N$ is arbitrary.
		
		We address Property (b). For $\tilde u=f\circ\,(l_1,\cdots,l_m)$ with $f\in\cb^\infty(\R^m)$, $l_1,\dots,l_m\in E'$, $m\in\N$ we calculate
		\begin{multline*}
			\lim_{N\to\infty}\,\int\limits_E|\nabla \tilde u|^2\de\mu_N
			\\=\lim_{N\to\infty}\,\int\limits_E\sum_{i,j=1}^m\partial_if(l_1(z),\cdots,l_m(z))\,\partial_jf(l_1(z),\cdots,l_m(z))\,\langle l_i, l_j\rangle\de\mu_N(z)\\
			=\,\int\limits_E\sum_{i,j=1}^m\partial_if(l_1(z),\cdots,l_m(z))\,\partial_jf(l_1(z),\cdots,l_m(z))\,\langle l_i, l_j\rangle\de\mu(z)
			=\,\int\limits_E|\nabla \tilde u|^2\de\mu.
		\end{multline*}
		Hence
		$\lim_{N}\E^N(u,u)\,=\,\E^\infty(u,u)$
		for $u\in\dom_\textnormal{pre}(\E^\infty)$, which concludes the proof.
	\end{proof}
	
	\subsection{A Gaussian measure and orthogonal projections: An example with a non-convex perturbing potential}\label{sec:gauss}
	
	In the final part of our survey, we present a frame in which the abstract assumptions of \cref{cond:hamza,cond:mucken} systematically hold
	and a convergence result in infinite dimension can be retrieved from \cref{cor:mosco}.
	We start with a finite measure space $(\Omega,\mathcal A,\lambda)$ and a non-degenerate, mean zero Gaussian measure $\tilde \mu$ on
	the state space $E=L^2(\Omega,\lambda)$. Norm and scalar product on $E$ are denoted by $|\cdot|$, respectively $\langle\cdot,\cdot\rangle$.
	We further assume that $E$ is separable.
	Let $E$ be densely embedded into another real Hilbert space $(\tilde H,(\cdot,\cdot)_{\tilde H})$ and $(A,\dom(A))$ be a self-adjoint operator on $\tilde H$
	such that \vspace{1ex}
	\begin{center}
		\begin{tabular}{l}
			$\bullet$ $A$ has pure point spectrum contained in the interval $[c,\infty)$ for some $c>0$,\vspace{1ex}\\
			$\bullet$ $E=\dom(\sqrt{A})$ with $(\sqrt{A}\,h,\sqrt{A}\,h)_{\tilde H}=\langle h,h\rangle$ for $h\in E$,\vspace{1ex}\\
			$\bullet$ the restriction $(\cdot,\cdot)_{\tilde H}|_{E\times E}$ is given by the covariance of $\tilde\mu$.\vspace{1ex}
		\end{tabular}
	\end{center}
	The last point says that for $h,k\in E$ the dual pairing w.r.t.~the inner product of $\tilde H$ reads
	\begin{equation*}
		(h,k)_{\tilde H}\,=\,\int\limits_E\langle h,z\rangle\,\langle k,z\rangle\de\tilde \mu(z).
	\end{equation*}
	It shall be noted that the listed properties do not restrict the class of Gaussian measures on $E$ which have mean zero and an injective covariance operator .
	We look at an increasing family ${(V_N)}_{N\in\overline\N}$ of closed subspaces of $E$ with $V_i\subset V_j$ if $1\leq i\leq j\leq\N$ and $V_\infty=E$.
	The images of $\tilde\mu$ under the orthogonal projections $P_N:E\to V_N$, $N\in\overline\N$, then serve as reference measures in our setting.
	To define a perturbation of these reference measures we consider a function $f$ from $\R$ to $\R$ with bounded variation.
	We have to assume another property which links $f$, $\tilde\mu$ and $\lambda$. As stated in \cref{cond:lvl} below,
	for any number in $\R$ at which $f$ is discontinuous the
	corresponding level set of $\tilde\mu$ is almost surely $\lambda$ -negligible.
	The set of real numbers at which $f$ is discontinuous is denoted by $U_f$.
	\begin{cond}\label{cond:lvl}
		Let $f:\R\to\R$ be a function of bounded variation such that
		\begin{equation*}
			\lambda(\{h=a\})\,=\,0\quad\textnormal{for }\tilde\mu\textnormal{ -a.e. }h\in E,\textnormal{ if }a\in U_f.
		\end{equation*}
	\end{cond}
	Now we define a perturbing potential by
	\begin{equation*}
		Q_f:E\ni h\,\mapsto\,\int\limits_\Omega f(h(\omega))\de\lambda(\omega)\,\in[\,-\|f\|_\infty\,\lambda(\Omega)\,,\,\|f\|_\infty\,\lambda(\Omega)\,].
	\end{equation*}
	A word should be said concerning the measurability of $Q_f$. By Lebesgue's dominated convergence $Q_f$ is continuous for $f\in\cb(\R)$.
	For $a\in\R$, choosing a monotone increasing sequence ${(f_m)}_{m\in\N}\subset\cb(\R)$ with  $\sup_{m}f_m(x)=\eins_{(a,\infty)}(x)$, $x\in\R$,
	yields the measurability of $Q_{\eins_{(a,\infty)}}$ by the monotone convergence theorem. By a monotone class argument the set 
	$\{\,g:\R\to\R\,|\,Q_g$ is measurable $\}$ contains all bounded functions which are $\mathcal B(\R)$ -mea\-surable.
	
	\begin{lem} \label{lem:wemecap}
		The weighted sequence of image measures ${(\,e^{-Q_f}(\tilde\mu\circ P_N^{-1})\,)}_{N\in\N}$ converges weakly towards $e^{-Q_f}\tilde\mu$, i.e.
		\begin{equation}\label{eqn:wemecapi}
			\lim_{N\to\infty}\,\int\limits_{E}g\circ P_N\,\exp(-Q_f\circ P_N)\de\tilde\mu \,=\,\int\limits_Eg\,\exp(-Q_f)\de\tilde\mu
		\end{equation}
		for $g\in\cb(E)$. Moreover,
		\begin{equation}\label{eqn:wemecapii}
			\lim_{N\to\infty}\,\int\limits_{E}g\,\exp(-Q_f\circ P_N)\de\tilde\mu \,=\,\int\limits_Eg\,\exp(-Q_f)\de\tilde\mu.
		\end{equation}
	\end{lem}
	\begin{proof}
		The lemma is an application of \cref{lem:wemeco}. The function $f$ can be approximated with sequences 
		${(f_m^\textnormal{min})}_{m\in\N},{(f_m^\textnormal{maj})}_{m\in\N}\subset\cb(\R)$ in the following sense. The inequality
		\begin{equation*}
			-\|f\|_\infty\,\leq\,f_m^\textnormal{min}(x)\,\leq\,f(x)\,\leq\,f_m^\textnormal{maj}(x)\,\leq\,\|f\|_\infty
		\end{equation*}
		holds for $m\in\N$ and $x\in\R$. Furthermore, if $x\in\R\setminus U_f$, then
		\begin{equation*}
			\lim_{m\to\infty} f_m^\textnormal{min}(x)\,=\,f(x)\,=\,\lim_{m\to\infty}f_m^\textnormal{maj}(x).
		\end{equation*}
		Such an approximation can be obtained for any bounded function on $\R$, e.g.~using the one-dimensional tent functions and setting
		\begin{equation*}
			f_m^\textnormal{min}\,=\,\sum_{\alpha\in (1/m)\Z}\,\Big(\inf_{y\in[\alpha-\frac{1}{m},\alpha+\frac{1}{m}]}f(y)\Big)\,\chi_{1/m}^\alpha
		\end{equation*}
		and
		\begin{equation*}
			f_m^\textnormal{maj}\,=\,\sum_{\alpha\in (1/m)\Z}\,\Big(\sup_{y\in[\alpha-\frac{1}{m},\alpha+\frac{1}{m}]}f(y)\Big)\,\chi_{1/m}^\alpha.
		\end{equation*}
		
		The set $U_f$ is at most countable, because $f$ is of bounded variation. Hence, 
		under \cref{cond:lvl},
		there exists a $\tilde\mu$ -nullset $\mathcal N\subset E$ such that $\lambda(\{\,\omega\,|\,h(\omega)\in U_f\,\})=0$ holds true for $h\in E\setminus\mathcal N$.
		By Lebesgue's dominated convergence, $\lim_m Q_{f^\textnormal{min}_m}(h)\,=\,Q_f(h)$ as well as $\lim_m Q_{f^\textnormal{maj}_m}(h)\,=\,Q_f(h)$ for $h\in E\setminus \mathcal N$.
		A second use of Lebesgue's dominated convergence yields the strong convergences, $\lim_m \exp(-Q_{f^\textnormal{min}_m})=\exp(-Q_f)$ and
		$\lim_m \exp(-Q_{f^\textnormal{maj}_m})=\exp(-Q_f)$, in $L^2(E,\tilde \mu)$.
		
		Since $\exp(-Q_{f^\textnormal{maj}_m}),\exp(-Q_{f^\textnormal{min}_m})\in\cb(E)$ for $m\in\N$ and
		\begin{equation*}
			\exp(-Q_{f^\textnormal{maj}_m}(h))\,\leq\,\exp(-Q_f(h))\,\leq\,\exp(-Q_{f^\textnormal{min}_m}(h))
		\end{equation*}
		for $h\in E$, we can apply \cref{lem:wemeco}  
		in the frame of $\prod_{N\in\overline\N}L^2(E,\tilde\mu\circ P_N^{-1})$. \cref{eqn:wemecapi} is proven.
		
		As to the second claim of this lemma, we want to obtain a convergence result in $L^2(E,\tilde\mu)$, so we apply \cref{lem:wemeco}
		in the frame of $\prod_{N\in\overline\N}L^2(E,\tilde\mu)$. On the one hand,
		\begin{equation*}
			\exp(-Q_{f^\textnormal{maj}_m}(P_Nh))\,\leq\,\exp(-Q_f(P_Nh))\,\leq\,\exp(-Q_{f^\textnormal{min}_m}(P_Nh))
		\end{equation*}
		for $N,m\in\N$ and $h\in E$. 
		On the other,
		\begin{gather*}
			\lim_{N\to\infty}\int\limits_E|\,\exp(-Q_{f^\textnormal{maj}_m}\circ P_N)-\exp(-Q_{f^\textnormal{maj}_m})\,|^2\de\tilde\mu\,=\,0,\\
			\lim_{N\to\infty}\int\limits_E|\,\exp(-Q_{f^\textnormal{min}_m}\circ P_N)-\exp(-Q_{f^\textnormal{min}_m})\,|^2\de\tilde\mu\,=\,0
		\end{gather*}
		for $m\in\N$ follows by Lebesgue's dominated convergence since $\exp(-Q_{f^\textnormal{maj}_m})$, $\exp(-Q_{f^\textnormal{min}_m})$ $\in\cb(E)$.
		Now, \cref{eqn:wemecapii} is a consequence of \cref{lem:wemeco}. and the proof is competed.
	\end{proof}
	
	As to the relevant partition functions we have
	\begin{equation*}
		Z_N\,:=\,\int\limits_E\exp(-Q_{f}\circ P_N)\de\tilde\mu\,\overset{N\to\infty}{\longrightarrow}\,\int\limits_E\exp(-Q_{f})\de\tilde\mu\,=:\,Z_\infty\in(0,\infty).
	\end{equation*}
	Since ${(Z_N)}_N$ is a convergent sequence of real numbers, we don't include it into the analysis below to shorten notation.
	The weighted measure $\exp(-Q_f\circ P_N)\tilde\mu$ is denoted by $\mu_N$ for $N\in\overline\N$.
	We define the relevant Dirichlet forms for the concluding results of this article. 
	Let $(\E^N,\dom(\E^N))$ denote of the smallest closed extension on $L^2(E,\mu_N)$ of the form in \cref{eqn:mingrad} for $\N\in\overline\N$, 
	i.e.~the minimal gradient form which have been analysed in \cref{sec:stagra}
	(we are now in the special case where $H=E=L^2(\Omega,\lambda)$). 
	For $N\in\N$ we also consider another Dirichlet form, which represents a similar yet slightly different point of view.
	We want to consider $Q_f$ as a perturbing potential for $\tilde\mu\circ P_N^{-1}$. Such an approach is taken in \cite{boune} with fixed examples
	for the respective choices of an $L^2$ space $E$, a Gaussian measure $\tilde\mu$ on $E$ and  increasing subspaces ${(V_N)}_{N}$. 
	There, the focus lies on the law of a Brownian bridge from $0$ to $0$ with state $E=L^2((0,1),\de x)$ and 
	subspaces $V_N$, which are the linear span of indicator functions $\eins_{[2^{-N}(i-1),2^{-N}i)}$, $i=1,\dots, 2^N$, $N\in\N$.
	In \cref{thm:bougen} we generalize \cite[Theorem 5.6]{boune} to our more abstract setting.
	$(\tilde \E^N,\dom(\tilde\E^N))$ denotes the smallest closed extension on $L^2(V_N,e^{-Q_f}\tilde\mu\circ P_N^{-1})$ of 
	\begin{equation*}
		\tilde\E^N(u,v)\,:=\,\int\limits_{V_N}\langle\nabla\tilde u,\nabla\tilde v\rangle\,\exp(-Q_f)\de(\tilde\mu\circ P_N^{-1})
	\end{equation*}
	with pre-domain 
	\begin{align*}
		\Big\{\,u\in L^2(V_N,e^{-Q_f}\tilde\mu\circ P_N^{-1})\,\Big|\,u(\cdot)=\tilde u(\cdot)\,(\tilde\mu\circ P_N^{-1})\textnormal{ -a.e. on }V_N
		\textnormal{for some }\tilde u\in\fcb(V_N)\,\Big\}.
	\end{align*}

	\begin{prop}\label{prop:only}
		Let $f$ be as in \cref{cond:lvl}. We consider the converging Hilbert spaces of $L^2(E,\mu_N)$, $N\in\N$, with limit $L^2(E,e^{-Q_f}\tilde \mu)$.
		
		${(\E^N,\dom(\E^N))}_N$ converges to $(\E^\infty,\dom(\E^\infty))$ in the sense of Mosco.
	\end{prop}
	\begin{proof}
		We verify the assumptions of \cref{cor:mosci}. To do so,
		we choose eigenvectors $\eta_1,\eta_2,\dots$ of $A$ which form an orthonormal basis of $E$. 
		Let $d\in\N$ and
		\begin{equation*}
			\pi_{\eta_d}:E\ni h\,\mapsto\,
			\,\langle \eta_d,h\rangle\,\in\R
		\end{equation*}
		We set $S_{\eta_d}:=\textnormal{Ker}(\pi_{\eta_d})$ and recall $J_{\eta_d}$ from \cref{eqn:jeta}.
		We define $(\E^{N,\eta_d},\dom(\E^{N,\eta_d}))$ for $d\in\N$ and $N\in\overline\N$ as in \cref{sec:stagra}.
		At first, we argue briefly why the assumptions of \cref{cor:mosci} are fulfilled in the trivial the case $f\equiv 0$, i.e.~$\mu_N=\tilde\mu$ for $N\in\overline\N$.
		Then, we can generalize using perturbation methods from \cref{sec:lowse}, in particular \cref{lem:wemeco,lem:mopert}.
		
		Assume $f\equiv0$. Let $d\in\N$. \cref{cond:hamza,cond:mucken} for the family ${(\mu_N\circ J_{\eta_d}^{-1})}_{N\in\overline\N}$ can be checked easily
		with the disintegration formula given in \cite[Proposition 5.5]{albe}:
		\begin{equation}\label{eqn:normaldis}
			\mu_N\circ J_{\eta_d}^{-1}(A)\,=\,\tilde\mu\circ J_{\eta_d}^{-1}(A)\,
			=\,\frac{1}{\sqrt{2\pi\sigma^2}}\,\int\limits_{S_{\eta_d}}\int\limits_{\R} \eins_A(s,x)\,e^{-x^2/(2\sigma^2)}\de x\de \nu_{\eta_d} (s)
		\end{equation}
		for $A\in\mathcal B(S_{\eta_d}\times\R)$, where $\sigma^2={\langle A\eta_d,\eta_d\rangle}^{-1}$ and
		$\nu_{\eta_d}$ is the image of $\tilde\mu$ under $E\ni h\mapsto h-(\pi_{\eta_d}h)\,\eta_d\in S_{\eta_d}$ (i.e.~under the first component of $J_{\eta_d}$).
		The issue of form domains can be settled with \cite[Proposition 3.2]{rocla}.
		To apply the latter result we have to do a remark on the existence of a gradient for $u\in \dom(\textstyle{\sum_i\E^{\infty,\eta_i}}\big)$.
		Indeed, each element $u\in \dom(\textstyle{\sum_i\E^{\infty,\eta_i}}\big)$
		can be assigned a gradient $\nabla u$, a $\tilde\mu$ -class of measurable maps $E\to E$ with $\langle\nabla u,\nabla u\rangle \in L^1(E,\tilde\mu)$.
		It is given by 
		\begin{align}
			&\nabla u: E\ni h\,\mapsto\,\sum_{i=1}^\infty \frac{\partial u}{\partial \eta_i}(h)\,\eta_i,\qquad\textnormal{where, for }i\in\N\textnormal{, we define}\label{eqn:weakgrad}\\ 
			&\frac{\partial u}{\partial \eta_i}(J_{\eta_i}^{-1}(s,x))\,:=\, (u\circ J_{\eta_i}^{-1}\,(s,\,\cdot\,))'(x)
			\quad\textnormal{ as a }\de\nu_{\eta_i}\times\de x\textnormal{ -class on }S_{\eta_i}\times \R.\label{eqn:weakpar}
		\end{align}
		The right hand side of \cref{eqn:weakpar} is well-defined, since $u\circ J_{\eta_i}^{-1}(s,\cdot)\in H^{1,1}_\textnormal{loc}(\R)$ for $\nu_{\eta_i}$~-a.e. $s\in S_{\eta_i}$
		(see \cref{eqn:conf}). The assignment of \cref{eqn:weakgrad} thus
		extends the gradient in the sense of G{\^a}teaux derivatives, which has been defined at the beginning of \cref{sec:stagra} for cylindrical smooth functions.
		Moreover,
		\begin{equation*}
			\int\limits_E\langle\nabla u,\nabla u\rangle\de\tilde\mu\,=\,\sum_{i=1}^\infty \E^{\infty,\eta_i}(u,u)
		\end{equation*}
		for $u\in\dom(\sum_i\E^{\infty,\eta_i})$.
		The action of the gradient of $u\in\dom(\sum_i\E^{\infty,\eta_i})$
		on an element $h\in$ span$(\{\eta_1,\eta_2,\dots\})$ can be interpreted as a (weak) directional derivative because of the chain rule, as follows.
		For $d\leq m$, $x_1,\dots,x_m\in\R$ and $z\in E$ we have 
		\begin{align*}
			u\Big(z+\sum_{i=1}^mx_i\,\eta_i\Big)\,&=\,u\circ J_{\eta_d}^{-1}\circ J_{\eta_d}\Big(z+\sum_{i=1}^mx_i\,\eta_i\Big)
			\\&=\,u\circ J_{\eta_d}^{-1}\Big(\,z-\langle\eta_d,z\rangle\,\eta_d+\sum_{\substack{i=1\\i\neq d}}^mx_i\,\eta_i\,\,,\,\,\langle\eta_d,z\rangle+x_d\,\Big).
		\end{align*}
		Hence, in accordance with \cref{eqn:weakpar}, $\R^m\ni x\mapsto u(z+x_1\eta_1+\dots+x_m\eta_m)$ is an element in $H^{1,1}_\textnormal{loc}(\R^m)$ for $\tilde\mu$ -a.e.~$z\in E$,
		whose $d$-th partial derivative reads
		\begin{equation*}
			\partial_{x_d} u\Big(z+\sum_{i=1}^mx_i\,\eta_i\Big)\,=\,\frac{\partial u}{\partial \eta_d}\Big(z+\sum_{i=1}^mx_i\,\eta_i\Big)\qquad\de x\textnormal{ -a.e.}
		\end{equation*}
		This, in turn, implies that $\R\ni s\mapsto u(z+sa_1\eta_1+\dots+sa_m\eta_m)$ is an element in $H^{1,1}_\textnormal{loc}(\R)$ for $a\in\R^m$
		and $\tilde\mu$ -a.e.~$z\in E$ with
		\begin{equation}\label{eqn:gradpre}
			\Big(u\Big(z+\,\,\cdot\,\,\sum_{i=1}^ma_i\,\eta_i\Big)\Big)'(s)\,=\,
			\sum_{d=1}^ma_d\,\frac{\partial u}{\partial \eta_d}\Big(z+s\sum_{i=1}^ma_i\,\eta_i\Big)\qquad\de s \textnormal{ -a.e.}
		\end{equation}
		Let $h\in$ span$(\{\eta_1,\eta_2,\dots\})$, i.e.~$h=\sum_{i=1}^ma_i\eta_i$ for some $m\in\N$ with $a_i=\langle \eta_i,h\rangle$ for $i=1,\dots,m$. We obtain
		\begin{equation}\label{eqn:gradcar}
			\langle \nabla u(z+s\,h),h\rangle\,=\,(u(z+\,\cdot\,h))'(s)\qquad\tilde\mu(\de z)\times\de s\textnormal{ -a.e.}
		\end{equation}
		by \cref{eqn:weakgrad} and \cref{eqn:gradpre}. With the existence of a gradient for elements of $\dom(\sum_i\E^{\infty,\eta_i})$, which has the property of \cref{eqn:gradcar},
		we see that the uniqueness result provided in \cite[Proposition 3.2]{rocla} is a stronger statement than the last assumption of \cref{cor:mosci}, the equality of domains
		\begin{equation*}
			\dom(\textstyle{\sum_i\E^{\infty,\eta_i}}\big)=\dom(\E^\infty).
		\end{equation*}
		This concludes our discussion about the case $f\equiv 0$.

		We now turn the attention to a non-trivial choice for $f$, in accordance with \cref{cond:lvl}. Since $\exp(-Q_f)$ is bounded uniformly on $E$ from below and above by positive numbers,
		we only have to care about \cref{cond:desint,cond:liminf}. The first one is handled via \cref{lem:wemeco}. The proper tool to tackle the second is \cref{lem:mopert}.
		Let's start with the verification of \cref{cond:desint} regarding the family $\mu_N\circ J_{\eta_d}^{-1}$, $N\in\overline\N$, where $d\in\N$ is fixed. 
		Taking into account the perturbing potential, the disintegration, which results from \cref{eqn:normaldis}, following the scheme of \cref{eqn:desdef} is given by
		\begin{equation*}
			\mu_N\circ J_{\eta_d}^{-1}(A)
			\,=\,\int\limits_{S_{\eta_d}}\int\limits_{\R} \eins_A(s,x)\,\frac{1}{z_s^N}\,e^{-Q_f\circ P_N\circ J_{\eta_d}^{-1}(s,x)}\,e^{-x^2/2\sigma^2}\de x\,z_s^N\de \nu_{\eta_d} (s)
		\end{equation*}
		for $A\in\mathcal B(S_{\eta_d}\times\R)$ and $N\in\overline\N$, where 
		\begin{equation*}
			z_s^N:=\int_\R e^{-Q_f\circ P_N\circ J_{\eta_d}^{-1}(s,x)}\,e^{-x^2/2\sigma^2}\de x.
		\end{equation*}
		Let $g\in\cb(S_{\eta_d}\times\R)$. We have to show
		\begin{multline}\label{eqn:finali}
			\lim_{N\to\infty}\,\int\limits_{S_{\eta_d}}\Big|\,\frac{1}{z_s^N}\,\int\limits_\R g(s,x)\, e^{-Q_f\circ P_N\circ J_{\eta_d}^{-1}(s,x)}\,e^{-x^2/2\sigma^2}\de x\,\Big|^2 \,z_s^N\de\nu_{\eta_d}(s)\\
			=\int\limits_{S_{\eta_d}}\Big|\,\frac{1}{z_s^\infty}\,\int\limits_\R g(s,x)\, e^{-Q_f\circ J_{\eta_d}^{-1}(s,x)}\,e^{-x^2/2\sigma^2}\de x\,\Big|^2 \,z_s^\infty\de\nu_{\eta_d}(s)
		\end{multline}
		We make an observation based on \cref{rem:onmucki}.
		If \cref{eqn:finali} is true for two bounded, continuous functions $g_1$ and $g_2$, then \cref{eqn:finali} also holds for their sum $g_1+g_2$.
		So, we can w.l.o.g.~assume that $g$ is a non-negative function. We verify \cref{eqn:finali} by proving the strong convergence of
		$z_{(\cdot)}^N\overset{N}{\to}z_{(\cdot)}^\infty$ as well as
		\begin{equation*}
			\underbrace{\int\limits_\R g(\cdot,x)\, e^{-Q_f\circ P_N\circ J_{\eta_d}^{-1}(\cdot,x)}\,e^{-x^2/2\sigma^2}\de x}_{=:G_N(\cdot)}\,\overset{N}{\longrightarrow}\,
			\underbrace{\int\limits_\R g(\cdot,x)\,e^{-Q_f\circ J_{\eta_d}^{-1}(\cdot,x)}\,e^{-x^2/2\sigma^2}\de x}_{=:G_\infty(\cdot)}
		\end{equation*}
		in $L^2(S_{\eta_d},\nu_{\eta_d})$.
		Indeed, once this is accomplished, the strong convergence of $G_N^2(\cdot)/z_{(\cdot)}^N\overset{N}{\to} G_\infty^2(\cdot)/z_{(\cdot)}^\infty$ in $L^1(S_{\eta_d},\nu_{\eta_d})$
		follows, since the sequence ${(z_s^N})_{N\in\N}$ is bounded from below by a positive number uniformly in $s$ and 
		${(G_N(s))}_{N\in\N}$ is bounded from above uniformly in $s$. 
		To obtain a convergence result in $L^2(S_{\eta_d},\nu_{\eta_d})$ we apply \cref{lem:wemeco} in the frame of 
		$\prod_{N\in\overline\N}L^2(S_{\eta_d},\nu_{\eta_d})$. The required comparison functions are constructed similarly as in the proof of \cref{lem:wemecap}
		with the continuous approximations $-\|f\|_\infty\leq f_m^{\textnormal{min}}(\cdot)\leq f(\cdot)\leq f_m^{\textnormal{maj}}(\cdot)\leq\|f\|_\infty$, $m\in\N$.
		For each $m\in\N$ and $N\in\overline\N$ we have
		\begin{align*}
			G_N(s)\,&\geq\,\int\limits_\R g(s,x)\, e^{-Q_{f_m^{\textnormal{maj}}}\circ P_N\circ J_{\eta_d}^{-1}(s,x)}\,e^{-x^2/2\sigma^2}\de x\,=:\,G_{N,m}^-(s)\\
			\textnormal{and}\quad G_N(s)\,&\leq\,\int\limits_\R g(s,x)\, e^{-Q_{f_m^{\textnormal{min}}}\circ P_N\circ J_{\eta_d}^{-1}(s,x)}\,e^{-x^2/2\sigma^2}\de x=:G_{N,m}^+(s).
		\end{align*}
		The continuity of $Q_{f_m^{\textnormal{maj}}}$ and $Q_{f_m^{\textnormal{min}}}$ on $E$ implies $G_{N,m}^+\,\overset{\textnormal{s.}}{\underset{N}{\longrightarrow}}\,G^+_{\infty,m}$
		and $G_{N,m}\,\overset{\textnormal{s.}}{\underset{N}{\longrightarrow}}\,G_{\infty,m}$ on $L^2(S_{\eta_d},\nu_{\eta_d})$ for $m\in\N$ by a multiple use of Lebesgue's dominated convergence.
		
		We now argue why $\lim_m G^+_{\infty,m}=\lim_m G^-_{\infty,m}=G_\infty$ holds strongly in $L^2(S_{\eta_d},\nu_{\eta_d})$.
		The set $U_f$ is at most countable, because $f$ is of bounded variation. Hence, 
		there exists a $\tilde\mu$ -nullset $\mathcal N\subset E$ such that $\lambda(\{\,\omega\,|\,h(\omega)\in U_f\,\})=0$ holds true for $h\in E\setminus\mathcal N$
		under \cref{cond:lvl}.
		We set $\mathcal N_s:=J_{\eta_d}^{-1}(s,\cdot)(\mathcal N)\subset\R$ for $s\in S_{\eta_d}$. For $\nu_{\eta_d}$ -a.e.~$s\in S_{\eta_d}$ the set $\mathcal N_s$ is a Lebesgue nullset.
		By repeatedly using Lebesgue's dominated convergence we build an argumentation as follows. First, 
		$\lim_m Q_{f^\textnormal{min}_m}(h)\,=\,Q_f(h)$ as well as $\lim_m Q_{f^\textnormal{maj}_m}(h)\,=\,Q_f(h)$ for $h\in E\setminus \mathcal N$.
		Secondly, $\lim_m G^+_{\infty,m}(s)=\lim_m G^-_{\infty,m}(s)=G_\infty(s)$ for $\nu_{\eta_d}$ -a.e.~$s\in S_{\eta_d}$. Finally, 
		$\lim_m G^+_{\infty,m}=\lim_m G^-_{\infty,m}=G_\infty$ holds strongly in $L^2(S_{\eta_d},\nu_{\eta_d})$.
		
		So, $\lim_NG_N=G_\infty$ holds strongly in $L^2(S_{\eta_d},\nu_{\eta_d})$ by \cref{lem:wemeco}. The corresponding convergence of 
		$\lim_Nz_{(\cdot)}^N=z_{(\cdot)}^\infty$ is already implied, since it results from the case where $g\equiv\eins_{S_{\eta_d}\times\R}$.
		The verification of \cref{cond:desint} regarding the family $\mu_N\circ J_{\eta_d}^{-1}$, $N\in\overline\N$, is completed.
		
		We address \cref{cond:liminf} as the last step of this proof. Let $d\in\N$ be fixed. We want to apply the perturbation result of \cref{lem:mopert} to deal with the relevant densities
		$\exp(-Q_f\circ P_N\circ J_{\eta_d}^{-1})$, $N\in\overline\N$. To do so, we use the Jordan decomposition of the function $f$.
		Let TV$(f)\in[0,\infty)$ denote the total variation of $f$.
		There exist monotone increasing functions $f_1,f_2:\R\to [0,\textnormal{TV(f)}]$ such that $f=a+f_1-f_2$ for some constant $a\in\R$, see e.g.~\cite[Chapter 5]{royden}. 
		We define the functionals
		\begin{align*}
			R_N:E\ni h\,&\,\mapsto\,\phantom{-}\int\limits_{\Omega}\big(\eins_{\{P_N\eta_d\geq 0\}}(\omega)\, f_1(h(\omega))-\eins_{\{P_N\eta_d< 0\}}(\omega)\, f_2(h(\omega))\big)\de\lambda(\omega)\\
			\textnormal{and}\qquad T_N:E\ni h\,&\,\mapsto\,\phantom{-}\int\limits_\Omega\big(\eins_{\{P_N\eta_d< 0\}}(\omega)\, f_1(h(\omega))-\eins_{\{P_N\eta_d\geq 0\}}(\omega)\, f_2(h(\omega))\big)\de\lambda(\omega).
		\end{align*}
		for $N\in\overline\N$.
		If $s\in S_{\eta_d}$, $-\infty<x\leq y<\infty$ and $N\in\overline\N$, then 
		\begin{align*}
			\exp{(-R_N\circ P_N\circ J_{\eta_d}^{-1})}(s,x)\,&=\,\exp{(-R_N(P_Ns+x\,P_N\eta_d))}\\
			&\geq \,\exp{(-R_N(P_Ns+y\,P_N\eta_d))}\,=\,\exp{(-R_N\circ P_N\circ J_{\eta_d}^{-1})}(s,y)
		\end{align*}
		and further
		\begin{align*}
			\exp{(-T_N\circ P_N\circ J_{\eta_d}^{-1})}(s,x)\,&\leq\,\exp{(-T_N\circ P_N\circ J_{\eta_d}^{-1})}(s,y).\\
		\end{align*}
		Since we can write
		\begin{equation*}
			\exp(-Q_f\circ P_N\circ J_{\eta_d}^{-1})\,=\,\exp(-a\,\lambda(\Omega))\,\exp(-R_N\circ P_N\circ J_{\eta_d}^{-1})\,\exp(-T_N\circ P_N\circ J_{\eta_d}^{-1}),
		\end{equation*}
		the family $\mu_N\circ J_{\eta_d^{-1}}$, $N\in\overline\N$ satisfies \cref{cond:liminf} by a double application of \cref{lem:mopert}. This concludes the proof.
	\end{proof}

	\begin{thm}\label{thm:bougen}
		Let $f$ be as in \cref{cond:lvl}. We consider the converging Hilbert spaces of $L^2(V_N,e^{-Q_f}\tilde \mu\circ P_N^{-1})$, $N\in\N$, with limit $L^2(E,e^{-Q_f}\tilde \mu)$.
		
		${(\tilde\E^N,\dom(\tilde\E^N))}_N$ converges to $(\E^\infty,\dom(\E^\infty))$ in the sense of Mosco.
	\end{thm}
	\begin{proof}
		We proof Mosco convergence by verifying the two conditions of \mimii{}. We start with (a).
		Let $N\in\N$. If $v\in L^2(V_N,e^{-Q_f}\tilde\mu\circ P_N^{-1})$ is in the pre-domain of $\tilde\E^N$, then
		choosing a representative $\tilde v\in\fcb(V_N)$ of $v$ we have 
		\begin{multline*}
			\tilde\E^N(v,v)\,=\,\int\limits_E|\nabla\tilde v\,(P_N(h))|^2\de\mu_N(h)\,=\,\int\limits_E|\nabla\tilde v\,(h)|^2\de\mu_N(h)
			\\=\,\int\limits_E|\nabla (\tilde v\circ P_N)\,(h)|^2\de\mu_N(h)\,=\,\E^N(v\circ P_N,v\circ P_N).
		\end{multline*}
		Since the image form of $(\E^N,\dom(\E^N))$ under $P_N$ is a closed form on $L^2(V_N,e^{-Q_f}\tilde\mu\circ P_N^{-1})$,
		its domain $\{\,u\,|\,u\circ P_N\in \dom(\E^N)\,\}$ must contain the whole of $\dom (\tilde\E^N)$ and furthermore
		\begin{equation}\label{eqn:imdom}
			\tilde\E^N(v,v)\,=\,\E^N(v\circ P_N,v\circ P_N),\,\qquad v\in\dom(\tilde\E^N).
		\end{equation}
		
		Let $g_1,g_2\in\cb(E)$.
		The convergence $\lim_{N\to\infty}g_1\circ P_N=g_1$ holds strongly in $L^2(E,\tilde\mu)$. It follows from \cref{eqn:wemecapii} that
		\begin{equation*}
			\lim_{N\to\infty}\,\int\limits_{E}(g_1\circ P_N)\,g_2\de\mu_N \,=\,\int\limits_Eg_1\,g_2\,e^{-Q_f}\de\tilde\mu.
		\end{equation*}
		As a consequence , we have $g\circ P_N\,\overset{\textnormal{s.}}{{\underset{N}{\longrightarrow}}}\,g$ 
		for $g\in\cb(E)$ in the sense of $\prod_{N\in\overline\N}L^2(E,\mu_N)$.
		
		Let now ${(u_N)}_{N\in\overline\N}\in \prod_{N\in\overline\N} L^2(V_N,e^{-Q_f}\tilde\mu\circ P_N^{-1})$ be a weakly convergent section.
		\begin{equation*}
			\lim_{N\to\infty}\,\int\limits_E (u_N\circ P_N)\,(g\circ P_N)\de\mu_N\,=\,\lim_{N\to\infty}\,\int\limits_{V_N} u_N\,g\,e^{Q_f}\de(\tilde\mu\circ P_N^{-1})
			\,=\,\int\limits_E u_\infty\,g\,e^{-Q_f}\de\tilde\mu
		\end{equation*}
		for $g\in\cb(E)$ and hence $u_N\circ P_N\,\overset{\textnormal{w.}}{{\underset{N}{\rightharpoonup}}}\,u$ referring to $\prod_{N\in\overline\N}L^2(E,\mu_N)$
		by virtue of \wecrit{}. 
		Now \cref{prop:only}, \mimii{} and \cref{eqn:imdom} imply $u_\infty\in\dom(\E)$ with
		\begin{equation*}
			\E(u_\infty,u_\infty)\,\leq\,\liminf_{N\to\infty}\,\E^N(u_N\circ P_N,u_N\circ P_N)
			\,\leq\,\liminf_{N\to\infty}\tilde\E^N(u_N,u_N)
		\end{equation*}
		assuming that $u_N\in\dom(\tilde\E^N)$ for infinitely many $N$ and the right hand side of the inequality is finite.
		Property (a) is proven. 
		
		As to (b), let $u\in L^2(E,\mu_\infty)$ be in the pre-domain of $\E^\infty$ with representative $\tilde u\in\fcb$.
		Then, $\tilde u\circ P_N\in\fcb(V_N)$ for $N\in\N$ with $\nabla (\tilde u\circ P_N)\,(h)=P_N\nabla\tilde u\,(P_N h)$ for $h\in E$ by the chain rule.
		Let $u_N\in L^2(V_N,e^{-Q_f}\tilde\mu\circ P_N^{-1})$ denote the class of $\tilde u\circ P_N$. 
		The convergence $u_N\,\overset{\textnormal{s.}}{{\underset{N}{\longrightarrow}}}\,u$ in the sense of 
		$\prod_{N\in\overline\N}L^2(V_N,e^{-Q_f}\tilde \mu\circ P_N^{-1})$ is an immediate consequence of \cref{eqn:wemecapi}, 
		the equality $P_N=P_N^2$ and the transformation of integrals.
		Now (b) follows from
		\begin{multline*}
			\limsup_{N\to\infty}\,\tilde\E^N(u_N,u_N)\,=\,\limsup_{N\to\infty}\,\int\limits_{V_N}|P_N\nabla\tilde u|^2\,e^{-Q_f}\de(\tilde\mu\circ P_N^{-1})\\
			\leq\,\limsup_{N\to\infty}\,\int\limits_{V_N}|\nabla\tilde u|^2\,e^{-Q_f}\de(\tilde\mu\circ P_N^{-1})\,=\,\E^\infty(u,u)
		\end{multline*}
		together with Property (a). This concludes the proof.
	\end{proof}

\hypersetup{bookmarksdepth=-2}
\section*{Acknowledgements}
We gratefully acknowledge financial support by the DFG through the project GR 1809/14-1.
Sincere thanks to Dr.~Sema Co\c{s}kun for her kind help with the figures in this article.
We thank an anonymous referee whose comments greatly helped to improve the quality of the article. 

\hypersetup{bookmarksdepth}

\bibliographystyle{plain}
\bibliography{literaturs}
\end{document}